\documentclass[12pt,leqno]{article}
\usepackage[utf8]{inputenc}
\usepackage{tikz}
\usepackage[official]{eurosym}
\usepackage{mathtools,fullpage,amsmath,theorem}
\usepackage{epic,eepic}
\usepackage{amssymb,amsmath}
\usepackage{accents}

\usepackage{cases}

\def\Mfor{\mathbb{M}[\vec{U}]}
\def\Mforc{\mathbb{M}_{\langle\kappa_1,...\kappa_n\rangle}[\vec{U}]}
\def\MforI{\mathbb{M}_I[\vec{U}]}
\def\min{{\rm min}}
\def\max{{\rm max}}
\def\sup{{\rm sup}}
\def\otp{{\rm otp}}
\def\Lim{{\rm Lim}}
\def\Inx{{\rm Index}}
\def\succ{{\rm Succ}}

\def\llvdash{{\|\hskip-2pt \raise 3pt\hbox{\vrule
height 0.25pt width 0.4cm}}}

\def\l{{\langle}}
\def\r{{\rangle}}

\def\calP{\mathcal P}

\def\upr{\upharpoonright}

\def\oa{{\overline A^{\,\lower 7pt_{\hbox{$\scriptstyle\bet}}
\hbox{$\scriptstyle 0\tau$}}}}

\def\bet{\beta}

\def\llvdash{{\|\hskip-2pt \raise 3pt\hbox{\vrule height 0.25pt
width 0.4cm}}}

\newtheorem{theorem}{Theorem}[section]
\newtheorem{lemma}[theorem]{Lemma}
\newtheorem{corollary}[theorem]{Corollary}
\newtheorem{proposition}[theorem]{Proposition}
{\theorembodyfont{\rmfamily}
\newtheorem{definition}[theorem]{Definition}}
{\theorembodyfont{\rmfamily}
}
{\theorembodyfont{\rmfamily}
\newtheorem{claim}{Claim}}
{\theorembodyfont{\rmfamily}
}
{\theorembodyfont{\rmfamily}
}
 \DeclareMathOperator{\dom}{dom}

\newcommand{\pr}{\medskip\noindent\textit{Proof}. }

\newcommand{\lusim}[1]{\smash{\underset{\raisebox{1.2pt}[0cm][0cm]{$\sim$}}
{{#1}}}}

\def\dom{{\rm dom}}

\def\otp{{\rm otp}}

\def\Ult{{\rm Ult}}

\def\llvdash{{\|\hskip-2pt \raise 3pt\hbox{\vrule
height 0.25pt width 0.15cm}}}

\def\Vdashbks{\hbox{$\Vdash\!\!\!\!{\raise2pt\hbox
{$\scriptscriptstyle\backslash$}}$}}
\providecommand{\keywords}[1]
{
  \small	
  \textbf{\textit{Keywords---}} #1
}
\providecommand{\subjectclass}[1]
{
  \small	
  \textbf{\textit{Mathematical Subject Classification---}} #1
}

\title{Sets in Prikry and Magidor Generic Extensions}
\author{ Tom Benhamou and Moti Gitik\footnote{ The work of the second author was partially supported by ISF grant No.58/14.}}

\begin{document}

\maketitle
\begin{abstract}
We continue \cite{PrikryCase} and study sets in generic extensions by the Magidor forcing and by the Prikry forcing with non-normal ultrafilters.
  \end{abstract}
  
  .
  
  \subjectclass{03E40, 03E55, 03E35,03E99}.
  \section*{Introduction}
In the paper \cite{PrikryCase} by V. Kanovei, P. Koepke and the second author, subforcings of the Prikry forcing were characterized.
Here we extend the analyzes to the Magidor forcing, introduced by M. Magidor in \cite{ChangeCofinality} (more recent account can be found in  \cite{GitHandbook} or in a very recent
nice and detailed paper by G. Fuchs \cite{Fuchs}).
\\The following is shown:
\\
\\
\textbf{Theorem \ref{Main theorem}}\emph{ Let \(\vec{U}\) be a coherent sequence in \(V\),  \(\langle\kappa_1,...\kappa_n\rangle\) be a sequence such that \(o^{\vec{U}}(\kappa_i)<\min(\nu\mid 0<o^{\vec{U}}(\nu))\), let \(G\) be \(\Mforc\)-generic\footnote{\(\Mforc\) is Magidor forcing with the coherent sequence \(\vec{U}\) above a condition which has \(\langle\kappa_1,...,\kappa_n\rangle\) as its ordinal sequence.} and let \(A\in V[G]\) be a set of ordinals. Then there exists \(C'\subseteq C_G\) such that \(V[A]=V[C']\), where $C_G$ is the Magidor sequence of $G$.}
\\
\\
One of the main methods used in the proof was the construction of a forcing \(\MforI\in V\), which is a projection of Magidor forcing \(\Mfor\). This forcing is a Magidor type forcing which uses only measures from \(\vec{U}\) with index \(i\in I\). Moreover, \(\MforI\) adds a prescribed subsequence \(C_I:=(C_G)\restriction I\)  as a generic object, where \(I\subseteq \lambda_0\) is a set of indexes in \(\lambda_0=\otp(C_G)\). Hence, we may examine the intermediate extensions \(V\subseteq V[C_I]\subseteq V[C_G]\) as an iteration of two forcings, which resemble \(\Mfor\).\\
A consequence of this theorem is the classification of all complete subforcings of $\Mfor$, this result is stated in theorem \ref{classification}.

Another direction addressed in this paper, is an attempt to extend the results of \cite{PrikryCase} to Prikry forcings with non-normal ultrafilters.
 The full generalization is not possible.
Thus, for example, once $\kappa$ is a $\kappa-$compact cardinal, then the Prikry forcing with a non-normal ultrafilter over $\kappa$ can add a generic for any $\kappa-$distributive forcing of size $\kappa$, see \cite{Git-kcompact} for more on this. Here we show that even from a single measurable, one can produce counter examples to generalizations of \cite{PrikryCase} to non-normal ultrafilters.
\\Namely the following is proved:
\\
\\
\textbf{Theorem 7.1}\emph{ Suppose that $V$ satisfies GCH and $\kappa$ is a measurable cardinal. Then in a generic cofinality preserving extension there is a $\kappa-$complete ultrafilter $U$ over $\kappa$ such that the Prikry forcing with $U$ adds a Cohen subset to $\kappa$ over $V$.
In particular, this forcing has a non-trivial subforcing which preserves regularity of $\kappa$.}
\\
\\
However, if one restricts to the Prikry forcing with $P-$point ultrafilters, then the following holds:
\\
\\
\textbf{Theorem 6.7}
\emph{Let $\mathbb{U}=\l U_a \mid a \in [\kappa]^{<\omega}\r$ consists of P-point ultrafilters over $\kappa$.
Then for every new set of ordinals $A$ in $V^{P(\mathbb{U})}$, $\kappa$ has cofinality $\omega$ in $V[A]$, where $P(\mathbb{U})$ denotes the tree Prikry  forcing with $\mathbb{U}$.
}
\\
\\

The paper is organized as follows:
\begin{itemize}
    \item Sections $1-5$ present the results for Magidor forcing.
    \item Section $6$ presents the generalization of \cite{PrikryCase} to the tree Prikry forcings with  $P-$points.
    \item Section $7$ presents the proof for Theorem 7.1.
\end{itemize}

\newpage
\section*{Notations}
\begin{itemize}
\item \(V\) denotes the ground model.
\item For any set \(A\), \(V[A]\) denote the minimal model of ZFC containing \(V\) and \(\{A\}\).
    \item 
\(\prod_{j=1}^{n} A_j \) increasing sequences \(\langle a_1,...,a_n\rangle\) where \(a_i\in A_i\).
\item \(\overset{m}{\underset{i=1 }{\prod}}\overset{n}{\underset{j=1}{\prod}} A_{i,j} \) left-lexicographically increasing sequences (which is denoted by \(\leq_{LEX}\)).
\item \([\kappa]^{\alpha}\) increasing sequences of length \(\alpha\).
\item \([\kappa]^{<\omega}=\underset{n<\omega}{\bigcup}[\kappa]^n\).
\item \(^{\alpha}[\kappa]\) not necessarily increasing sequences, i.e functions with domain \(\alpha\) and range \(\kappa\).
\item \(^{\omega>}[\kappa]=\underset{n<\omega}{\bigcup} \ ^n[\kappa]\).
\item \(\langle\alpha,\beta\rangle\) an ordered pair of ordinals. \((\alpha,\beta)\) the interval between \(\alpha\) and \(\beta\).
\item \(\vec{\alpha}=\langle\alpha_1,...,\alpha_n\rangle\) , \(|\vec{\alpha}|=n \ , \ \vec{\alpha}\setminus\langle\alpha_i\rangle=\langle\alpha_1,...,\alpha_{i-1},\alpha_{i+1},...,\alpha_n\rangle\).
\item  For every \(\alpha<\beta\), The Cantor normal form (abbreviated C.N.F) equation is
\(\alpha+\omega^{\nu_1}+...+\omega^{\nu_m}=\beta\),  \(\nu_1\geq...\geq\nu_m\) are unique. If \(\alpha=0\) this is the C.N.F of \(\beta\), otherwise, this is the C.N.F difference of \(\alpha,\beta\).
\item The limit otder if $\alpha$, denoted by \(o_L(\alpha)=\gamma\), where \(\alpha=\omega^{\gamma_1}+...+\omega^{\gamma_n}+\omega^{\gamma}\) (C.N.F).
\item \(\Lim(A)=\{\alpha \in A \mid \sup(A\cap\alpha)=\alpha\}\).
\item \(\succ(A)=\{\alpha\in A\mid \sup(A\cap\alpha)<\alpha\}\).
\item \(\underset{i\in I}{\biguplus}A_i\) is the union of \(\{A_i \mid i\in I\}\) with the requirement that \(A_i\)'s are pairwise disjoint.
\item If \(f:A\rightarrow B\) is a function then for every \(A'\subseteq A\), \(B'\subseteq B\)
\begin{center}
\(f''A'=\{f(x)\mid x\in A'\}\) ,  \( f^{-1''}B'=\{x\in A\mid f(x)\in B'\}\).

\end{center}
\item Let \(B\subseteq \langle\alpha_\xi\mid \xi<\delta\rangle=A \) be sequences of ordinals, \begin{center}
 
\(\Inx(B,A)=\{\xi<\delta\mid \exists b\in B \ \alpha_\xi=b\}\).
\end{center}
\item Let \(\mathbb{P}\) be a forcing notion, \(\sigma\) a formula in the forcing language and \(p\in\mathbb{P}\). If \(\underaccent{\sim}{A}\) is a \(\mathbb{P}\)-name, then
\begin{center}
\(p\ ||\ \underaccent{\sim}{A}\) means "there is \(a\in V\) such that \(p\Vdash\overset{\vee}{a}=\underaccent{\sim}{A} \)".
\end{center}
\item Let \(p,q\in\mathbb{P}\) then
"\(p,q\) are compatible in \(\mathbb{P}\)" if there exists \(r\in\mathbb{P}\) such that \(p,q\leq_{\mathbb{P}} r\). Otherwise, if they are incompatible denote it by \( p\bot q\).
\item In any forcing notion, \(p\leq q\) means "\(q\) extends \textit{p}".
\item The notion of complete subforcing, complete embedding and projection is used as defined in \cite{Shelah}.
\end{itemize}
\newpage
\section{Magidor forcing}
 \begin{definition}
 A \textit{coherent sequence} is a sequence \\ \( \vec{U}=\langle U(\alpha,\beta) \mid \beta<o^{\vec{U}}(\alpha)\ ,\alpha\leq\kappa \rangle \) such that:
\begin{enumerate}
    \item \(U(\alpha,\beta)\) is a normal ultrafilter over \(\alpha\).
    \item Let \(j:V\rightarrow Ult(U(\alpha,\beta),V)\) be the corresponding elementary embedding , then $$j(\vec{U})\restriction\alpha=\vec{U}\restriction\langle\alpha,\beta\rangle$$
Where 
$$\vec{U}\restriction\alpha=\langle U(\gamma,\delta) \mid \delta<o^{\vec{U}}(\gamma) \ ,\gamma\leq\alpha \rangle$$ 
$$\vec{U}\restriction\langle\alpha,\beta\rangle=\langle U(\gamma,\delta) \mid (\delta<o^{\vec{U}}(\gamma) ,\ \gamma<\alpha) \vee (\delta<\beta ,\ \gamma=\alpha)\rangle$$

\end{enumerate}
 \end{definition}
\begin{flushright}
 \(\blacksquare\)
 \end{flushright}
Fix  a coherent sequence of ultrafilters \( \vec{U}\) with maximal measurable \( \kappa\). Assume that
 \(o^{\vec{U}}(\kappa)<\min(\nu\mid o^{\vec{U}}(\nu)>0):=\delta_0\). Let \(\alpha\leq\kappa\) with \(o^{\vec{U}}(\alpha)>0\), define 
 \begin{center}
  \(\bigcap U(\alpha,i)=\underset{i<o^{\vec{U}}(\alpha)}{\bigcap}U(\alpha,i)\)
 \end{center}
 We will follow the description of Magidor forcing as presented in \cite{GitHandbook}. \\
 \\
\begin{definition}\label{Definition of conditions in MU}
The \textit{Magidor forcing}, denoted by
\(\mathbb{M}[\vec{U}]\), consist of conditions \(p\) of the form
\(p=\langle t_1,...,t_n,\langle\kappa,B\rangle\rangle\).
 For every  \(1\leq i\leq n \), \(t_i\) is either an ordinal
 \(\kappa_i\) if \( o^{\vec{U}}(\kappa_i)=0\)
 or a pair \(\langle\kappa_i,B_i\rangle\)  if \ \(o^{\vec{U}}(\kappa_i)>0\) such that:
\begin{enumerate}
\item \(B\in \underset{\xi<o^{\vec{U}}(\kappa)}{\bigcap}U(\kappa,\xi)\), \ \ \(\min(B)>\kappa_n\).
    \item  for every  \(1\leq i\leq n \):
    \begin{enumerate}
    \item \(\langle\kappa_1,...,\kappa_n\rangle\in [\kappa]^{<\omega}\).
    \item \(B_i\in \underset{\xi<o^{\vec{U}}(\kappa_i)}{\bigcap}U(\kappa_i,\xi)\). 
    \item  \(\min(B_i)>\kappa_{i-1}\) \ \  (\(i>1\)).
    \end{enumerate}
\end{enumerate}

\end{definition}
\begin{flushright}
 \(\blacksquare\)
 \end{flushright}
 We shall adopt the following notations: 
 \begin{enumerate}
\item \(t_0=0 , t_{n+1}=\langle\kappa,B\rangle\).
\item \(o^{\vec{U}}(t_i)=o^{\vec{U}}(\kappa(t_i))\).
\item \(o^{\vec{U}}(t_i)>0\) then \(t_i=\langle\kappa_i,B_i\rangle=\langle\kappa(t_i),B(t_i)\rangle \).
\item\( o^{\vec{U}}(t_i)=0\) then \(t_i=\kappa_i=\kappa(t_i)\) .

\item \(\kappa(p)=\{\kappa(t_1),...,\kappa(t_n)\}\).
\item \(B(p)=\overset{n+1}{\underset{i=1}{\biguplus}}B(t_i)\).
  \end{enumerate}
 The ordinals \(\kappa_i\) are designated to form the eventual Magidor sequence and candidates for the sequence’s missing elements in the interval \((\kappa(t_{i-1}),\kappa(t_i))\)  (where \(t_0=0,\ \kappa(t_{n+1})=\kappa\)) are provided by the sets \(B(t_i)\).\\ 
 \\
\begin{definition}\label{Definition of MU order}
 Let \(p=\langle t_1,t_2,...,t_n,\langle\kappa,B\rangle\rangle,q=\langle s_1,...,s_m,\langle\kappa,C\rangle\rangle\in \mathbb{M}[\vec{U}]\) , define  the \textit{Magidor forcing order} by \(p \leq q\) ("q extends p") iff:
\begin{enumerate}
    \item \(n \leq m\).
    \item \( B \supseteq C\).
    \item \( \exists 1 \leq i_1 <...<i_n \leq m \) such that for every \(1 \leq j \leq m\):
    \begin{enumerate}
        \item If \( \exists 1\leq r\leq n\) such that \(i_r=j\) then \(\kappa(t_r)=\kappa( s_{i_r})\) and \(C(s_{i_r})\subseteq B(t_r)\).
        \item Otherwise \( \exists \ 1 \leq r \leq n+1\) such that \( i_{r-1}<j<i_{r}\) then 
        \begin{enumerate}
        \item \(\kappa(s_j) \in B(t_r)\).
        \item\(o^{\vec{U}}(s_j)<o^{\vec{U}}(t_{r})\).
        \item \(B(s_j)\subseteq B(t_r)\cap \kappa(s_j)\).
        \end{enumerate}
    \end{enumerate}
\end{enumerate}
The \textit{direct extension order} is defined by \( p \leq^{*} q\) iff:
\begin{enumerate}
    \item \( p \leq q\).
    \item \(n=m\).
\end{enumerate}

\end{definition}
\begin{flushright}
 \(\blacksquare\)
 \end{flushright}
 \underline{Remarks:} 
  \begin{enumerate}
 \item Let \( p=\langle t_1,...,t_n,\langle\kappa,B\rangle\rangle\). Assume we would like to add an element \(s_j\) to \(p\) between \(t_{r-1}\) and \(t_{r}\). It is possible only if \(o^{\vec{U}}(t_r)>0\).
 Moreover, let \(\xi=o^{\vec{U}}(s_j)\), then
 \begin{center}
  \(s_j\in \{\alpha\in B(t_r)\mid o^{\vec{U}}(\alpha)=\xi\}\)
 \end{center} 
 If \(s_j=\kappa(s_j)\) (i.e. \(\xi=0\)), then any \(s_j\) satisfying this requirement can be added. 
 If \(s_j=\langle\kappa(s_j),B(s_j)\rangle\) (i.e. \(\xi>0\)), Then according to definition \ref{Definition of MU order} (3.b.iii) \(s_j\) can be added iff 
 \begin{center}
         \(B(t_r)\cap \kappa(s_j)\in \underset{\xi'<\xi}{\bigcap}U(\kappa(s_j),\xi')\)
\end{center}
 \item If \(p=\langle t_1,...,t_n,\langle \kappa,B\rangle\rangle\in\Mfor\). Fix some \(1\leq j\leq n\) with \(o^{\vec{U}}(t_j)>0\). Then \(t_j\)  yields a Magidor forcing in the interval \((\kappa(t_{j-1}),\kappa(t_j))\) with the coherent sequence \(\vec{U}\restriction\kappa(t_j)\). \(t_j\) acts autonomously in the sense that the sequence produced by it is independent of how the sequence develops in other parts. This observation becomes handy when manipulating \(p\), since we can make local changes at \(t_j\) with no impact on other \(t_i\)'s.
  \end{enumerate}
 Let
  \(Y= \{ \alpha \leq \kappa  \mid  o^{\vec{U}}(\alpha)<\delta_0 \}\). From the coherency of \(\vec{U}\), it follows that \(Y \in \bigcap U(\kappa,i)\).
 For every \(\beta\in Y \) with \(o^{\vec{U}}(\beta)>0\), and \( i<\delta_0\) define
 \begin{center}
  \(Y(i)= \{ \alpha < \kappa  \mid  o^{\vec{U}}(\alpha) =i\} \) and \(Y[\beta]=\underset{ i<o^{\vec{U}}(\beta)}{\biguplus}Y(i)\)
   \end{center} 
 It follows that for every \(\beta\in Y\) and \(i<o^{\vec{U}}(\beta)\), 
 \(Y(i)\cap\beta \in 
U(\beta ,i)\).
To see this take \(\beta\leq\kappa\) in \(Y\) and \( j_{\beta,i}:V \rightarrow Ult(U(\beta ,i),V)\) . 
\begin{center}
\(Y(i)\cap\beta \in U(\beta ,i) \ \Leftrightarrow\) \( \beta \in j_{\beta,i}(Y(i)\cap\beta)\) 
\end{center}
By coherency, \(o^{j_{\beta ,i}(\vec{U}) }(\beta)=i \) and therefore 
\begin{center}
\( \beta\in j_{\beta, i}(Y(i)\cap\beta)= \{ \alpha<j_{\beta,i}(\beta) \mid o^{j_{\beta,i}(\vec{U})}(\alpha)=j_{\beta,i}(i)=i \}\).
\end{center}
Consequently, \( Y[\beta]\cap\beta\in \underset{i<o^{\vec{U}}(\beta)}{\bigcap}U(\beta,i)\).\\
 For \(  B \in \underset{i<o^{\vec{U}}(\beta)}{\bigcap}U(\beta,i)\) define recursively, \(B^{(0)}=B\),
 \begin{center}
  \( B^{(n+1)}=\{ \alpha \in B^{(n)} \mid (o^{\vec{U}}(\alpha)=0) \vee (B^{(n)} \cap \alpha \in \cap U(\alpha,i))  \}\) \\
   \end{center}
    Let \(B^{\star}=\underset{n<\omega}{\bigcap}B^{(n)}\)  it follows by induction that for all \( n<\omega\),
    \begin{center}
    \(B^{(n)}\in \underset{i<o^{\vec{U}}(\beta)}{\bigcap}U(\beta,i)\)
    \end{center}
    By \(\beta\)-completeness \( B^{\star}\in\underset{i<o^{\vec{U}}(\beta)}{\bigcap}U(\beta,i)\).
\(B^{\star}\) has the feature that ,
\begin{center}
\(\forall \alpha\in B^{\star}\) \(\alpha\cap B^{\star}\in \underset{i<o^{\vec{U}}(\alpha)}{\bigcap}U(\alpha,i) \)
\end{center}

 The previous paragraph indicates that by restricting to a dense subset of \(\mathbb{M}[\vec{U}]\) we can assume that given \(p=\langle t_1,t_2,...,t_n,\langle\kappa,B\rangle\rangle \in \mathbb{M}[\vec{U}]\), every choice of  ordinal in \(B(t_r)\) automatically satisfies the requirement that we discussed in remark (2).
 Formally, we work above \(\langle\langle\rangle ,\langle\kappa,Y\rangle\rangle\) and we directly-extend any \(p=\langle t_1,t_2,...,t_n,\langle\kappa,B\rangle\rangle\) as follows: \\For every \(1\leq r\leq n+1\) and \( i<o^{\vec{U}}(t_{r})\) define 
\begin{center}
 
\( B(t_{r},i):=Y(i) \cap B(t_r)^{\star} \in U(\kappa(t_r),i) \)
\end{center}
It follows that 
\begin{center}
\( B^{\star}(t_r):= \underset{i<o^{\vec{U}}(t_{r})}{\biguplus} B(t_r,i)\in\underset{i<o^{\vec{U}}(t_r)}{\bigcap}U(\kappa(t_r),i)\). 
\end{center}
 Shrink \(B(t_r)\) to \(B^{\star}(t_r)\) to obtain
 \begin{center}
  \(p\leq^*p^*=\langle t_1',...,t_n',\langle\kappa,B^{\star}\rangle\rangle\) \end{center}
  \begin{center}
  \(t'_r=\begin{cases}
t_r & o^{\vec{U}}(t_r)=0\\
\langle\kappa(t_r),B^{\star}(t_r)\rangle & otherwise \end{cases}\)
  \\
 \end{center}
 This dense subset also simplifies \(\leq\) to
 \begin{center}
  \(p\leq q\) iff \(\kappa(p)\subseteq\kappa(q)\) , \( B(p)\subseteq B(q)\)
 \end{center} 
When applying the revised approach regarding the large sets, it is apparent that \(B(t_r,i)\) provide candidates, precisely, for the \(i\)-limit indices in the final sequence \(C_{G}\) i.e. of indices \(\gamma\) such that \(o_L(\gamma)=i\) (for the definition of \(o_L(\gamma)\) see Notations). This is stated formally in proposition \ref{Order type of generic club}.  \\
 Recall that:
 \begin{itemize}
     \item \(\Mfor\) satisfies \(\kappa^+-c.c.\)
     \item Let \(p=\langle t_1,...,t_n,\langle\kappa,B\rangle\rangle\in\Mfor\)  and denote \(\nu=\kappa(t_j)\) where \(j\) is the minimal such that \(o^{\vec{U}}(t_j)>0\). Then above \(p\) there is \(\nu-_ {\leq^*}\)closure.
     \item \(\Mfor\) satisfies the Prikry condition.
\end{itemize}
 Let \(G\subseteq \Mfor\) be generic, define 
 \begin{center}
  \(C_{G}=\bigcup \{\kappa(p) \mid p\in G\}\)
 \end{center}

  We will abuse notation by considering \(C_G\) as a the canonical enumeration of the set \(C_G\). \(C_{G}\) is closed and unbounded in \(\kappa\). Therefore, the order type of \(C_{G}\) determines the cofinality of \(\kappa\) in \(V[G]\). The next propositions can be found in \cite{GitHandbook}.\\
 \\
\begin{proposition}\label{Reconstruction of generic}
 \textit{ Let \(G\subseteq\Mfor\) be generic. Then \(G\) can be reconstructed from \(C_{G}\) as follows}
     \begin{center}
      \(G=\{p\in\Mfor\mid (\kappa(p)\subseteq C_{G}) \wedge (C_{G}\setminus\kappa(p)\subseteq B(p))\}\)
     \end{center}
      \textit{Therefore \(V[G]=V[C_{G}]\).}
      \begin{flushright}
 \(\blacksquare\)
 \end{flushright}
      
\end{proposition}
 \begin{proposition}\label{Order type of generic club}
 \textit{ Let G be \(\Mfor\)-generic and \(C_{G}\) the corresponding Magidor sequence. Let \(\langle t_1,...,t_n,\langle\kappa,B\rangle\rangle\in G\), then }
    \begin{center}
    \( \otp( (\kappa(t_i),\kappa(t_{i+1}))\cap C_{G} )=\omega^{o^{\vec{U}}(\kappa(t_{i+1}))}\)
    \end{center}
    \textit{Thus if \(\kappa(t_{i+1})=C_G(\gamma)\) then \(o_L(\gamma)=o^{\vec{U}}(t_{i+1})\).}
     \begin{flushright}
 \(\blacksquare\)
 \end{flushright}
 \end{proposition}
\begin{corollary}\label{Cor of changing cofinality}
 If $o^{\vec{U}}(\kappa)<\delta_0$, then $cf^V(o^{\vec{U}}(\kappa))=cf^{V[G]}(o^{\vec{U}}(\kappa))$ and 
 \(cf^{V[G]}(\kappa)=cf(\omega^{o^{\vec{U}}(\kappa)}).\)
 \begin{flushright}
 \(\blacksquare\)
 \end{flushright}
 \end{corollary}
 Let \(p=\langle t_1,...,t_n,\langle\kappa,B\rangle\rangle\in G\). By  proposition \ref{Order type of generic club}, for each \(i\leq n\) one can determine the position of \(\kappa(t_i)\) in  \(C_{G}\).  Namely, \(C_{G}(\gamma(t_i,p))=\kappa(t_i) \) where 
   \begin{center}
   
 \( \gamma(t_i,p):=\underset{j\leq i}{\sum}\omega^{o^{\vec{U}}(t_j)}<\omega^{o^{\vec{U}}(\kappa)} \) \ \ (*) 
\end{center}
Addition and power are of ordinals. The equation (*) induces a C.N.F equation
\begin{center}
   
 \( \gamma=\sum_{r=1}^{m}\omega^{o^{\vec{U}}(t_{j_r})} \) \ \ (C.N.F)
 \end{center}
This indicates the close connection between Cantor normal form of the index \(\gamma\) in \(\otp(C_G)\) and the elements \(t_{j_1},...,t_{j_m}\) in $p$ which determines that \(\gamma(t_i,p)=\gamma\). Let \(q=\langle s_1,...,s_m,\langle\kappa,B'\rangle\rangle\) be another condition, by definition \ref{Definition of MU order} (3.b.ii), if \(s_j\) is an element of \(q\) which was added to \(p\) in the interval \((\kappa(t_r),\kappa(t_{r+1}))\) then \(o^{\vec{U}}(s_j)<o^{\vec{U}}(t_{r+1})\). Consequently, adding it to $p$, does not impact the Cantor normal form and
$\gamma(t_{r+1},p)=\gamma(t_{r+1},q)$.
      
\newpage
\section{Combinatorial properties }

 The combinatorial nature of \(\Mfor\) is most clearly depicted through the language of step-extensions as presented below. \\
\textbf{To perform a one step extension} of \(p=\langle t_1,t_2,...,t_n,\langle\kappa,B\rangle\rangle\):
\begin{enumerate}
\item Choose \(1 \leq r \leq n+1\) with \(0<o^{\vec{U}}(t_{r})\).
\item Choose \(i < o^{\vec{U}}(t_{r})\).
\item Choose an ordinal \(\alpha \in B(t_r,i)\).
\item Shrink the sets \(B(t_s,j)\) to \(C(t_s,j) \in U(t_s,j)\) for every \(1\leq s\leq n+1\) and set $$C(t_s)=\underset{j<o^{\vec{U}}(t_i)}{\biguplus}C_s(j)$$
\item For \(j<o^{\vec{U}}(\alpha)\) pick \(C(\alpha,j)\in U(\alpha,j)\),  \( C(\alpha,j)\subseteq B(t_r,j)\cap\alpha \) and set  $$C(\alpha)=\underset{j<o^{\vec{U}}(\alpha)}{\biguplus}C(\alpha,j)$$
\item Cut \(C(t_r)\) above \(\alpha \).

\end{enumerate}
Extend \(p\) to 
\begin{center}
\(p^{\frown}\langle\alpha,(C(t_s))^{n+1}_{s=1},C(\alpha)\rangle= \langle t'_1,...,t'_{i-1},\langle\alpha,C(\alpha)\rangle,t'_i,...,t'_n,\langle\kappa,C(t_{n+1})\rangle\rangle\)
\end{center}
\begin{center}
\(t'_r=\begin{cases}
t_r & o^{\vec{U}}(t_r)=0\\
\langle\kappa(t_r),C(t_r)\rangle & o.w. \end{cases}\)
\end{center}

It is clear that every extension of \(p\) with only one ordinal added is a one step extension. Next we introduce some notations which will describe a general step extension. The idea is simply to classify extensions according to the order of the measures the new elements of the sequence are chosen from.
\begin{definition}
Let \(p=\langle t_1,t_2,...,t_n,\underset{t_{n+1}}{\underbrace{\langle\kappa,B\rangle}}\rangle \in \Mfor\) \begin{enumerate}
\item  For \(1\leq i\leq n+1\)  define \textit{the tree or order types}  \(T_i(p)=\ ^{ \omega>}[o^{\vec{U}}(t_i)]\), with the order $$\langle x_1,...,x_m\rangle\preceq \langle x_1',...,x_{m'}'\rangle$$ iff there are \( 1 \leq i_1 <...<i_m \leq m' \) such that for every \(1 \leq j \leq m'\):
    \begin{enumerate}
        \item If \( \exists 1\leq r\leq m\) such that \(i_r=j\) then \(x_r=x'_j\).
        \item Otherwise \( \exists \ 1 \leq r \leq n+1\) such that if \( i_{r-1}<j<i_{r}\) then  \(x'_j<x_r\).
    \end{enumerate}
 We think of \(x_r\)'s as placeholders of ordinals from \(B(t_i,x_r)\) and the ordering is induced by definition \ref{Definition of MU order} (3). 
 \item \(T(p)=\prod^{n+1}_{i=1}T_i(p)\) with \(\preceq\) is the product order.
 \item Let \(X_i \in T_i(p) \ \  1\leq i \leq n+1\) , \(|X_i|=l_i\), \(X=\langle X_1,...,X_{n+1}\rangle\in T(p)\).
 \begin{enumerate}

\item Define
$$\vec{\alpha}_i=\langle\alpha_1,...,\alpha_{l_i}\rangle \in  B(p,X_i):=\prod_{j=1}^{l_i} B(t_i , X_i(j))$$
    \(X_i\) is called an \textit{extension-type below \(t_i\)} and \(\langle\alpha_1,...,\alpha_{l_i}\rangle\) is of \textit{type \(X_i\)}.
\item Define
\begin{center}
\(\vec{\alpha}=\langle\vec{\alpha_1},...,\vec{\alpha_{n+1}}\rangle \in B(p,X):=\overset{n+1}{\underset{ i=1}{\prod}}\overset{l_i}{\underset{ j=1}{\prod}} B(t_i , X_i(j))\)
\end{center}
\(X\) is called an \textit{extension-type of \(p\)} and \(\vec{\alpha}\) is of \textit{type \(X\)}. Notice that $X$ is uniquely determined by $\vec{\alpha}$ and the $o^{\vec{U}}()$ function.
\end{enumerate}
 
 \end{enumerate}
 \end{definition}
\begin{flushright}
 \(\blacksquare\)
 \end{flushright}
Notice that by our assumption \(|T(p)|<\min(\nu |0< o^{\vec{U}}(\nu))=\delta_0\). We also use:
\begin{itemize}
\item \(|X_i|=l_i\).
\item \(l_x=\max(i \mid X_i\neq\emptyset)\).
    \item \(x_{i,j}=X_i(j)\) \(\alpha_{i,j}=\vec{\alpha}_i(j)\).
    \item \(x_{i,l_i+1}=o^{\vec{U}}(t_i)\) and \(\alpha_{i,n+1}=\kappa(t_i)\).
    \item \(x_{mc}=x_{l_X,l_{l_X}}\) (i.e. the last element of X).
    \item \(o^{\vec{U}}(\vec{\alpha})=\langle o^{\vec{U}}(\alpha_{i,j})\mid x_{i,j}\in X\rangle\) is the type of \(\vec{\alpha}\).
\end{itemize}
\textbf{A general extension of p of type X} is of the form:
\begin{center}
\(p^{\frown}\langle\vec{\alpha},(C(x_{i,j}))_{x_{i,j}\in X}, (C(t_r))^{n+1}_{r=1}\rangle=p^{\frown}\langle\vec{\alpha},(C(x_{i,j}))_{^{i\leq n+1}_{j\leq l_i+1}}\rangle\)

\end{center}
where 
\begin{center}
\(p^{\frown}\langle\vec{\alpha},(C(x_{i,j}))_{^{i\leq n+1}_{j\leq l_i+1}}\rangle=\langle\vec{s_1},t'_1,...,\vec{s_n},t'_n,\vec{s_{n+1}},\langle\kappa,C\rangle\rangle\)

\end{center}
\begin{enumerate}

    \item \(\vec{\alpha}\in B(p,X)\).
    \item \(t'_s=\begin{cases}
t_s & o^{\vec{U}}(t_s)=0\\
\langle\kappa(t_s),C(t_s)\rangle & o.w. \end{cases}\) \\
For some pre chosen sets \(C(t_s)\in\underset{\xi<o^{\vec{U}}(t_s)}{\bigcap}U(\kappa(t_s),\xi)\)  ,\ \( C(t_s)\subseteq B(t_s)\).
\item \(\vec{s}_i(j)=\begin{cases}
\alpha_{i,j} & x_{i,j}=0\\
\langle\alpha_{i,j},C(x_{i,j})\rangle & o.w. \end{cases}\)\\
For some pre chosen sets \(C(x_{i,j})\in\underset{\xi<x_{i,j}}{\bigcap}U(\alpha_{i,j},\xi)\) , \(C(x_{i,j})\subseteq B(t_i)\cap\alpha_{i,j}\).
\item \(C\in \underset{\xi<o^{\vec{U}}(\kappa)}{\bigcap}U(\kappa,\xi)\) and  \(\min(C)>\max(\vec{s}_{n+1})\).
\end{enumerate}
 By the discussion succeeding definition \ref{Definition of MU order}, we conclude that $$p^{\frown}\langle\vec{\alpha},(C(x_{i,j}))_{^{i\leq n+1}_{j\leq l_i+1}}\rangle \in \Mfor$$ 
 We will more frequently use \(p^{\frown}\langle\vec{\alpha}\rangle\) with the same definition as above except we do not shrink any sets and simply take \(\alpha_{i,j}\cap B(t_i)=C(x_{i,j})\). Define
 \begin{center}
  \(p^{\frown} X=\{p^{\frown}\langle\vec{\alpha}\rangle \mid \vec{\alpha} \in B(p,X)\}\)
 \end{center}
 
 The set \(p^{\frown} X\)'s induce a partition of \(\Mfor\) above \(p\) as stated in the next proposition which is well known and follows directly from definition \ref{Definition of MU order}.
\begin{proposition}\label{partition of extension type}
 \textit{Let \(p\in \Mfor\) be any condition and
 \(p\leq q\in \Mfor\). Then there exists a unique \(\vec{\alpha}\in B(p,X)\) such that \(p^{\frown}\langle\vec{\alpha}\rangle\leq^* q\).}
\begin{flushright}
 \(\blacksquare\)
 \end{flushright}
 
\end{proposition} 
\textbf{Example}:
 Let 
 $$p=\langle\underset{t_1}{\underbrace{\langle\kappa(t_1),B(t_1)\rangle}},\underset{t_2}{\underbrace{\kappa(t_2)}},\underset{t_3}{\underbrace{\langle\kappa(t_3),B(t_3)\rangle}},\underset{t_4}{\underbrace{\langle\kappa(t_4),B(t_4)\rangle}},\underset{t_5}{\underbrace{\langle\kappa,B\rangle}}\rangle$$
 $$o^{\vec{U}}(t_1)=1,\ o^{\vec{U}}(t_2)=0, \ o^{\vec{U}}(t_3)=2, \ o^{\vec{U}}(t_4)=1 , \ o^{\vec{U}}(\kappa)=3$$
 
 Extend $p$ to $q$ as follows:
 $$q=p^{\frown}\langle\underset{\vec{\alpha_1}}{\underbrace{\langle\alpha_{1,1},\alpha_{1,2}}\rangle},\underset{\vec{\alpha_2}}{\underbrace{\langle\rangle}},\underset{\vec{\alpha_3}}{\underbrace{\langle\alpha_{3,1},\alpha_{3,2},\alpha_{3,3}\rangle}},\underset{\vec{\alpha_4}}{\underbrace{\langle\alpha_{4,1}\rangle}},\underset{\vec{\alpha_5}}{\underbrace{\langle\alpha_{5,1},\alpha_{5,2},\alpha_{5,3}\rangle}}\rangle$$
 \begin{center}
 \(o^{\vec{U}}(\alpha_{i,j})=\begin{cases}
 0 & \langle i,j\rangle=\langle 1,1\rangle,\langle 1,2\rangle,\\
 &\langle 3,2\rangle,\langle 4,1\rangle,\langle 5,1\rangle\\
 \\
 1 &  \langle i,j\rangle=\langle 3,1\rangle,\langle 3,3\rangle,\\
 &\langle 5,2\rangle\\
 \\
 2 &\langle i,j\rangle=\langle 5,3\rangle
 \end{cases}\)
 \end{center}
 Then the extension-type of \(\langle\vec{\alpha_1},\vec{\alpha_2},\vec{\alpha_3},\vec{\alpha_4},\vec{\alpha_5}\rangle\) is 
 \begin{center}
 \(X=\langle\underset{X_1}{\underbrace{\langle0,0\rangle}},\underset{X_2}{\underbrace{\langle\rangle}},\underset{X_3}{\underbrace{\langle 1,0,1\rangle}},\underset{X_4}{\underbrace{\langle 0\rangle}},\underset{X_5}{\underbrace{\langle 0,1,2\rangle}}\rangle\)
 \end{center}
 This can be illustrated as follows:\\
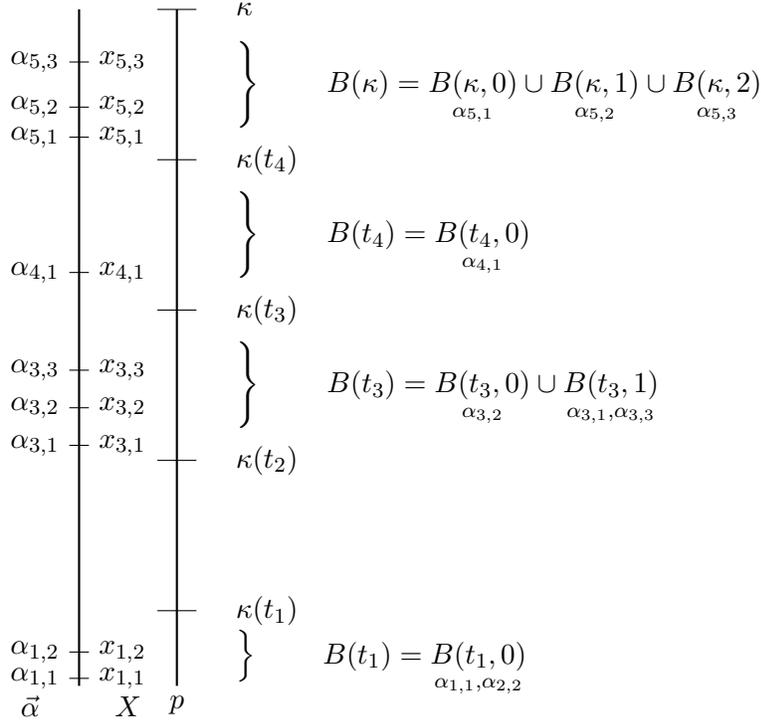
\begin{tikzpicture}[xscale=1.3]
\draw [thick] (1,0) -- (1,9);
\draw (0.9,0.1) -- (1.1,0.1);
\node[align=left, left] at (0.9,0.1)
{\(\alpha_{1,1}\)};
\node[align=left, right] at (1.1,0.1)
{\(x_{1,1}\)};
\draw (0.9,0.45) -- (1.1,0.45);
\node[align=left, left] at (0.9,0.45)
{\(\alpha_{1,2}\)};
\node[align=left, right] at (1.1,0.45)
{\(x_{1,2}\)};

\draw (0.9,3.2) -- (1.1,3.2);
\node[align=left, left] at (0.9,3.2)
{\(\alpha_{3,1}\)};
\node[align=left, right] at (1.1,3.2)
{\(x_{3,1}\)};
\draw (0.9,3.7) -- (1.1,3.7);
\node[align=left, left] at (0.9,3.7)
{\(\alpha_{3,2}\)};
\node[align=left, right] at (1.1,3.7)
{\(x_{3,2}\)};
\draw (0.9,4.2) -- (1.1,4.2);
\node[align=left, left] at (0.9,4.2)
{\(\alpha_{3,3}\)};
\node[align=left, right] at (1.1,4.2)
{\(x_{3,3}\)};

\draw (0.9,5.5) -- (1.1,5.5);
\node[align=left, left] at (0.9,5.5)
{\(\alpha_{4,1}\)};
\node[align=left, right] at (1.1,5.5)
{\(x_{4,1}\)};

\draw (0.9,7.3) -- (1.1,7.3);
\node[align=left, left] at (0.9,7.3)
{\(\alpha_{5,1}\)};
\node[align=left, right] at (1.1,7.3)
{\(x_{5,1}\)};
\draw (0.9,7.7) -- (1.1,7.7);
\node[align=left, left] at (0.9,7.7)
{\(\alpha_{5,2}\)};
\node[align=left, right] at (1.1,7.7)
{\(x_{5,2}\)};
\draw (0.9,8.3) -- (1.1,8.3);
\node[align=left, left] at (0.9,8.3)
{\(\alpha_{5,3}\)};
\node[align=left, right] at (1.1,8.3)
{\(x_{5,3}\)};

\draw [thick] (2,0) -- (2,9);
\draw (1.8,1) -- (2.2,1);
\node[align=left, right] at (2.5,0.3)
{\(\Big\} \ \ \ \ \ \ \ B(t_1)=\underset{\alpha_{1,1},\alpha_{2,2}}{B(t_1,0)}\)};
\node[align=left, right]  at (2.5,1)
{\(\kappa(t_1)\)};
\draw (1.8,3) -- (2.2,3);
\node[align=left, right]  at (2.5,3)
{\(\kappa(t_2)\)};
\draw(1.8,5) -- (2.2,5);
\node[align=left, right] at (2.5,4)
{\(\Bigg\} \ \ \ \ \ \ \ B(t_3)=\underset{\alpha_{3,2}}{B(t_3,0)}\cup \underset{\alpha_{3,1},\alpha_{3,3}}{B(t_3,1)}\)};
\node[align=left, right]  at (2.5,5)
{\(\kappa(t_3)\)};
\draw (1.8,7) -- (2.2,7);
\node[align=left, right] at (2.5,6)
{\(\Bigg\} \ \ \ \ \ \ \ B(t_4)=\underset{\alpha_{4,1}}{B(t_4,0)}\)};
\node[align=left, right]  at (2.5,7)
{\(\kappa(t_4)\)};
\draw (1.8,9) -- (2.2,9);
\node[align=left, right] at (2.5,8)
{\(\Bigg\} \ \ \ \ \ \ \ B(\kappa)=\underset{\alpha_{5,1}}{B(\kappa,0)}\cup\underset{\alpha_{5,2}}{ B(\kappa,1)}\cup \underset{\alpha_{5,3}}{B(\kappa,2)}\)};
\node[align=left, right]  at (2.5,9)
{\(\kappa\)};
\node[align=left, below] at (0.5,0)
{\(\vec{\alpha}\)};
\node[align=left, below] at (1.5,0)
{\(X\)};
\node[align=left, below] at (2,0)
{\(p\)};

\end{tikzpicture}
 \par
As presented in proposition \ref{partition of extension type}, a choice from the set \(p^{\frown}X\) is essentially a choice from $B(p,X)$, which is of the form \(\overset{n}{\underset{i=1}{\prod }}A_{i}\) where for every $1\leq i\leq n$, $A_{i} \in U_{i}$ and $U_i$'s are normal measures on a non decreasing sequence of measurable  \(\kappa_{1} \leq \kappa_{2} \leq ... \leq \kappa_{n} \). For the rest of this section, we prove some combinatorical properties of such sets.
\begin{lemma}\label{stabprod}
\textit{ Let \(\kappa_{1} \leq \kappa_{2} \leq ... \leq \kappa_{n} \) be any collection of measurable cardinals with normal measures \(U_{1},...,U_{n} \) respectively. Assume \(
F:\overset{n}{\underset{i=1}{\prod }}A_{i}\longrightarrow \nu\) where \( \nu <\kappa_1\) and \(A_{i} \in U_{i}\). Then there exists \( H_{i}\subseteq A_{i} \) \ \(H_{i} \in U_{i} \) such that \( \overset{n}{\underset{i=1}{\prod}}H_{i}\) is homogeneous for \(F\) i.e. $F\restriction\overset{n}{\underset{i=1}{\prod}}H_{i}$ is constant.}

\end{lemma}
\pr By induction on \(n\). The case \(n=1\) is a well known property of normal measures.
Assume that the lemma holds for \(n-1\) , and fix \(\vec{\eta}=\langle\eta_{1},...,\eta_{n-1}\rangle \in \overset{n-1}{\underset{i=1}{\prod }}A_{i}\).  Define
$$F_{\vec{\eta}}:A_{n}\setminus(\eta_{n-1}+1) \longrightarrow \nu , \ \ 
F_{\vec{\eta}}(\xi)=F(\eta_{1},...,\eta_{n-1},\xi)$$
 By the case n=1 there exists an homogeneous \(A_n\supseteq H(\vec{\eta}) \in U_{n} \) with color \(C(\vec{\eta}) < \nu\). Define  
\begin{center}
\( \underset{\vec{\eta} \in \prod_{i=1}^{n-1}A_{i}}{\Delta}H(\vec{\eta})=:H_{n}\) 
\end{center}
By the induction hypotheses applied to \(C:\overset{n-1}{\underset{i=1}{\prod }}A_{i}\rightarrow\nu\), there is an homogeneous set of the form  \(\overset{n-1}{\underset{i=1}{\prod}}H_{i}\) where \( A_{i} \supseteq H_{i} \in U_{i} \). To see that \(\overset{n}{\underset{i=1}{\prod}}H_{i}\) is homogeneous for \(F\), let \(\vec{\eta'},  \vec{\eta}\in \overset{n}{\underset{i=1}{\prod }}H_{i}\), and denote by $\eta_n,\eta'_n$, the ordinals $\max(\vec{\eta}),\max(\vec{\eta}')$ respectively. It follows that:
$$F(\vec{\eta})=F_{\vec{\eta}\setminus\langle\eta_n\rangle}(\eta_n) \underset{\eta_{n} \in H(\vec{\eta}\setminus\langle\eta_{n}\rangle)}{\underset{\uparrow}{=}} C(\vec{\eta}\setminus\langle\eta_{n}\rangle)  \underset{\vec{\eta}\setminus\langle\eta_{n}\rangle,\vec{\eta'}\setminus\langle\eta_{n}'\rangle \in \overset{n-1}{\underset{i=1}{\prod}}H_{i}}{\underset{\uparrow}{=}}=C(\vec{\eta'}\setminus\langle\eta_{n}'\rangle)=...=F(\vec{\eta'})$$
\begin{flushright}
 \(\blacksquare\)
 \end{flushright}
\begin{lemma}\label{ImportantCoor}
 \textit{Let \(\kappa_{1} \leq \kappa_{2} \leq ... \leq \kappa_{n} \) be a non descending  finite sequence of measurable cardinals with normal measures \(U_{1},...,U_{n} \) respectively. Assume \(
F:\overset{n}{\underset{i=1}{\prod }}A_{i}\longrightarrow B\) where \(B\) is any set, and \(A_{i} \in U_{i}\). Then there exists \(H_{i}\subseteq A_{i}\)   \(H_{i} \in U_{i} \) and set of coordinates \(I \subseteq  \{1,...,n\}\)  such that for every $\vec{\eta},\vec{\xi}\in \overset{n}{\underset{i=1}{\prod }}H_{i}$,
$$F(\vec{\eta})=F(\vec{\xi} )\leftrightarrow \vec{\eta}\restriction I=\vec{\xi}\restriction I$$
In other words, the function
\(F\restriction\overset{n}{\underset{i=1}{\prod}}H_{i}\) is well defined  modulo the equivalence relation:}
\begin{center}
\( \langle\alpha_1,...,\alpha_n\rangle\sim\langle\alpha_1',...,\alpha_n'\rangle \) \ \textit{iff} \( \forall i \in I \ \alpha_i=\alpha_i'\)
\end{center} \textit{and the induced function, \(\bar{F}\), is injective. We call the set $I$, a set of important coordinates.} 

\end{lemma}
\pr
By induction on \(n\), if \(n=1\) then it is immediate since for any \ \(f:A \rightarrow B\)  such that \(A \in U\) where \(U\) is a normal measure on a measurable cardinal \(\kappa\), \(B\) is any set, then there exists \(A \supseteq A' \in U\) for which \(F\restriction A'\) is either constant or injective.
Assume that the lemma holds for \(n-1 \ ,n>1\) and let \(F:\overset{n}{\underset{i=1}{\prod }}A_{i} \longrightarrow B\) be a function satisfying the conditions of the lemma. Define for every \(x_1 \in A_1\), \(F_{x_1}:\overset{n}{\underset{i=2}{\prod }}A_{i}\setminus (x_1+1)\longrightarrow B\)
\begin{center}
\(F_{x_1}(x_2,...,x_n)=F(x_1,x_2,...,x_n)\)
\end{center}
By the induction hypothesis , for every \(x_1 \in A_1\) there are \( A_{i} \supseteq A_{i}(x_1) \in U_{i} \) and  set of important coordinates \(I(x_1) \subseteq  \{2,...,n\}\).
Therer is $A_1\supseteq A'_1\in U_1$ such that function \(I:A_1 \rightarrow P(\{2,...,n\})\) is constant on \(A'_1\) with value \(I'\). For every \(i=2,...,n\) define \(A'_i= \underset{x_1 \in A_1}{\Delta}A_i(x_1)\). So far, \(\overset{n}{\underset{i=1}{\prod }}A'_{i}\) has the property:
$$(1) \ \ \forall
\langle x_1,x_2,...,x_n\rangle,\langle x_1,x'_2,...,x'_n\rangle \in \overset{n}{\underset{i=1}{\prod }}A'_{i} \text{ with the same first coordinate}$$
\begin{center}
\( F(x_1,x_2,...,x_n)=F(x_1,x'_2,...,x'_n)\) iff \(\forall i \in I'.  \ x_i=x'_i\)
\end{center}
In particular, \(\bar{F}\) is a well defined function modulo \(I' \cup \{1\}\). Next we determine if \(1\) is important. For every \(\langle\alpha,\alpha'\rangle \in A'_1 \times A'_1 \), define \(t_{\langle\alpha,\alpha'\rangle}:\overset{n}{\underset{i=2}{\prod }}A'_{i} \setminus (\alpha'+1)\rightarrow 2\)
\begin{center}

\(t_{\langle\alpha,\alpha'\rangle}(x_2,...,x_n)=1 \longleftrightarrow F(\alpha,x_2,...,x_n)=F(\alpha',x_2,...,x_n)\)
\end{center}
By lemma \ref{stabprod}, for \(i=2,...,n\) there are \( A'_{i} \supseteq A_{i}(\alpha,\alpha') \in U_{i} \) such that \(\overset{n}{\underset{i=2}{\prod }}A_{i}(\alpha,\alpha')\) is homogeneous for \(t_{\langle\alpha,\alpha'\rangle}\) with color \(C(\alpha,\alpha')\).
Taking the diagonal intersection over \(A'_1 \times A'_1 \) of the sets \(A_{i}(\alpha,\alpha')\), at each coordinate \(i=2,...,n\), we obtain \(H_{i} \in U_{i}\) such that \(\overset{n}{\underset{i=2}{\prod }}H_{i}\) is homogeneous for every \(t_{\langle\alpha,\alpha'\rangle} \). Finally, the function  \(C:A'_1 \times A'_1 \rightarrow 2\) yield an homogeneous \( A'_{1} \supseteq H_{1} \in U_{1} \)  with color \(C'\). 

\underline{Case 1:} \(C'=1\). 
Let us show that the important coordinates are \(I'\). 

Let \(\langle x_1,...,x_n\rangle,\langle x'_1,...,x'_n\rangle \in \underset{i=1}{\overset{n}{\prod}}H_{i}\), then,
$$ F(\scriptstyle x_1,...,x_n\textstyle)=F(\scriptstyle x'_1,...,x'_n\textstyle)\underset{F(x'_1,x'_2...,x'_n)=F(x_1,x'_2,...,x'_n)}{\underset{\uparrow}{\longleftrightarrow}}  F(\scriptstyle x_1,x_2,...,x_n\textstyle)=F(\scriptstyle x_1,x'_2,...,x'_n\textstyle)\underset{(1)}{\underset{\uparrow}{\leftrightarrow}}  \forall i. \in I'  \ x_i=x'_i$$

\underline{Case 2:} \(C'=0\).
Then we have a second property:
$$(2) \ \ \forall x_1,x_1'\in H_1\text{ and }\langle x_2,...,x_n\rangle\in\prod^n_{i=2}H_i. \  \ 
x_1=x_1'\leftrightarrow F(x_1,x_2...,x_n)=F(x_1',x_2,...,x_n)$$

We would like to claim that in this case the important coordinates are \(I=I' \cup \{1\}\) but we still have to shrink \(H_i\)'s, to eliminate all remaining counter examples for \(\bar{F}\) not being injective i.e. $\langle x_1,...,x_n\rangle, \langle x'_1,...,x'_n\rangle\in  \underset{i=1}{\overset{n}{\prod}}H_{i}$ such that
$$\langle x_1,...,x_n\rangle\neq \langle x'_1,...,x'_n\rangle mod \ I \ \text{ and } \  F(x_1,...,x_n)=F(x'_1,...,x'_n)$$

Take Any counter example and set
 \begin{center}
  \(\{x_1,...,x_n\}\cup\{x_1',...,x_n'\}=\{y_1,...,y_k\}\) (increasing enumeration)
 \end{center}
 To reconstruct \(\{x_1,...,x_n\},\{x_1',...,x_n'\}\) from \(\{y_1,...,y_k\}\) is suffices to know for example how \(\{x_1,...,x_n\}\) are arranged between \(\{x_1',...,x_n'\}\). There are finitely many ways \footnote{In general, the number of possibilities to arrange two counter examples into one increasing sequence depends on \(I\). Nevertheless, there is an upper bound: Think of  \(x_i\)'s as balls we would like to divide into \(n+1\) cells. The cells are represented by the intervals \((x_{i-1}', x_i']\) plus the cell for elements above \(x_n'\). There are \(\binom{2n}{n}\) such divisions. For any such division, we decide either the cell is \((x_{i-1}', x_i']\) or \((x_{i-1}', x_i')\). Hence, there are at most \(\binom{2n}{n}\cdot 2^n\) such arrangements. } for Such an arrangement. Therefore, if we succeed with eliminating examples of a fixed arrangement, then by $\sigma$-completeness of the measures we will be able to eliminate all counter example.
 
Fix such an arrangement,
 the increasing sequence \(\langle y_1,...,y_k\rangle \) is in the product of some \(k\) large sets \(\overset{k}{\underset{i=1}{\prod }}H_{n_i}\). We have to be careful since the sequence of measurables induced by \(n_1,...,n_k\) is not necessarily non descending. To fix this we can cut the sets \(H_i\)  such that in the sequence \(\langle\kappa_i\mid i=1,...,n\rangle\), wherever \(\kappa_i<\kappa_{i+1}\) then \(\min(H_{i+1})>\kappa_i=\sup(H_i)\). Therefore, assume that \(\langle \kappa_{n_i} \mid i=1,...,k\rangle\) is non descending. Define \(G: \overset{k}{\underset{i=1}{\prod }}H_{n_i} \rightarrow 2\)
\begin{center}
\(G(y_1,...,y_k)=1 \Leftrightarrow F(x_1,...,x_n)=F(x'_1,...,x'_n)\)
\end{center}
By lemma \ref{stabprod} there must be  \(U_i\ni H'_{i}\subseteq H_{i}\) homogeneous for \(G\) with value \(D\). If \(D=0\) we have eliminated from \(H_i\)'s all counter examples of that fixed ordering. Toward a contradiction, assume that \(D=1\), then every \( y_1,...,y_k \) yield a counter example \(\langle x_1,...,x_n\rangle,\langle x_1',...,x_n'\rangle\). By propety $(1)$,
\(x_1=x_1'\), hence assume without loss of generality that \(x_1<x_1'\),  fix \(x<w<y_2<...<y_n\), where \(x,w\in H'_1\) and \(y_i \in H'_{n_i}\) for \(i=2,...,k\).
Since $D=1$, it follows that \(G(x,y_2,...,y_k)=G(w,y_2,...,y_k)=1\), thus,
$$F(x,x_2,...,x_n)=F(x_1',x'_2,...,x'_n)=F(w,x_2,...,x_n)$$
which is a contradiction to property (2).
\begin{flushright}
 \(\blacksquare\)
 \end{flushright}

\newpage

\section{The Main Result Up to \(\kappa\)}

As stated in corollary \ref{Cor of changing cofinality}, Magidor forcing adds a closed unbounded sequence of  length \(\omega^{o^{\vec{U}}(\kappa)}\) to \(\kappa\).
It is possible to obtain a family of forcings that adds a sequence of any limit length to some measurable cardinal, using a variation of Magidor forcing as we defined it\footnote{Magidor's original formulation of \(\Mfor\) in \cite{ChangeCofinality} gives such a family}. Namely, let \(\vec{U}\) be a coherent sequence and \(\lambda_0<\min(\nu\mid o^{\vec{U}}(\nu)>0)\)  a limit ordinal
\begin{center}
(not necessarily C.N.F) \(\lambda_0=\omega^{\gamma_1}+...+\omega^{\gamma_n}\)  \ \ \ ,\(\gamma_n>0\)
\end{center}
Let  \(\langle\kappa_1,...\kappa_n\rangle\) be an increasing sequence such that \(o^{\vec{U}}(\kappa_i)=\gamma_i\). Define the forcing \(\Mforc\) as follows: \\
The root condition will be
\begin{center}
 \(0_{\Mforc}=\langle\langle\kappa_1,B_1\rangle,...,\langle\kappa_n,B_n\rangle\rangle\)
\end{center}
where \(B_1,...,B_n\) are as in the discussion following definition \ref{Definition of MU order}. The conditions of this forcing are any finite sequence that extends \(0_{\Mforc}\)  in the sense of definition \ref{Definition of MU order}. Since each \(\langle \kappa_i ,B_i\rangle\) acts autonomously, this forcing is essentially the same as \(\Mfor\). In fact, \(\Mfor\) is just \(\mathbb{M}_{\langle\kappa\rangle}[\vec{U}]\). The notation we used for \(\Mfor\) can be extended to  \(\Mforc\) since the conditions are also of the form \(\langle t_1,...,t_r,\langle\kappa,B\rangle\rangle\). Let 
\begin{center}
 \(\langle\langle\nu_1,C_1\rangle,...,\langle\nu_m,C_m\rangle\rangle\in \Mforc\)
\end{center}
then \(\mathbb{M}_{\langle\nu_1,...,\nu_m\rangle}[\vec{U}]\) is an open subset of \(\Mforc\) (i.e. \(\leq\)-upwards closed). Moreover, if \(G\subseteq \Mforc\) is any generic set with \(\langle\langle\nu_1,C_1\rangle,...,\langle\nu_m,C_m\rangle\rangle\in G\) then,
\begin{center}
    \((G)_{\langle\nu_1,...,\nu_m\rangle}:=G\cap\mathbb{M}_{\langle\nu_1,...,\nu_m\rangle}[\vec{U}]=\{p\in G\mid p\geq \langle\langle\nu_1,C_1\rangle,...,\langle\nu_m,C_m\rangle\rangle\}\)
    \end{center}
    is generic for \(\mathbb{M}_{\langle\nu_1,...,\nu_m\rangle}[\vec{U}]\). The filter \((G)_{\vec{\nu}}\) is essentially the same generic as \(G\) since it yield the same Magidor sequence and in particular \(V[(G)_{\vec{\nu}}]=V[G]\).
    
From now on the set \(B\) in \(\langle t_1,...,t_r,\langle\kappa,B\rangle\rangle\) will be suppressed and replaced by \(t_{r+1}=\langle \kappa,B\rangle\). An alternative way to describe \(\Mforc\) is through the following product
\begin{center}
\(\Mforc\simeq\Mfor_{\langle\kappa_1\rangle}\times(\Mfor_{\langle\kappa_2\rangle})_{>\kappa_1}\times...\times(\Mfor_{\langle\kappa_n\rangle})_{>{\kappa_{n-1}}}\)\\
\((\mathbb{M}_{\langle\nu_1,...,\nu_m\rangle}[\vec{U}])_{>\alpha}=\{\langle t_1,...,t_{r+1}\rangle\in \mathbb{M}_{\langle\nu_1,...,\nu_m\rangle}[\vec{U}]\mid \kappa(t_1)>\alpha\}\)
\end{center}
This isomorphism is induced by the embeddings
\begin{center}

 \(i_r:((\Mfor_{\langle\kappa_r\rangle})_{>\kappa_{r-1}}\rightarrow\Mforc\) \ \ , \(r=1,...,n\) \\
\(i_r(\langle s_1,...,s_{k+1}\rangle)=\langle\scriptstyle\langle\kappa_1,B_{1}\rangle,...,\langle\kappa_{r-1},B_{r-1}\rangle ,s_1,...,s_k,\underset{s_{k+1}}{\underbrace{\scriptstyle\langle\kappa_r,B(s_{k+1})\rangle}},...,\langle\kappa_{n},B_{n}\rangle \textstyle \rangle\)
\end{center} 
From this embeddings, it is clear that the generic  sequence produced by \((\Mfor_{\langle\kappa_r\rangle})_{>\kappa_{r-1}}\) is just \(C_G\cap(\kappa_{r-1},\kappa_r)\).\\
\\
The formula to compute coordinates still holds: 

Let \(p=\langle t_1,...,t_m,t_{m+1}\rangle\in \Mforc\). For each \(1\leq i\leq m\),  the coordinate of \(\kappa(t_i)\) in any Magidor sequence extending \(p\) is \(C_G(\gamma)=\kappa(t_i) \), where
\begin{center}
 \( \gamma=\underset{j\leq i}{\sum}\omega^{o^{\vec{U}}(t_j)}=:\gamma(t_i,p)<\lambda_0 \)
 
\end{center}
\begin{lemma}\label{Properties of CG}
\textit{Let \(G\) be generic for \(\Mforc\) and the sequence derived}
\begin{center}
 
\(C_G=\bigcup \{ \{\kappa(t_1),...,\kappa(t_l)\}\mid \langle t_1,...,t_l,t_{l+1}\rangle\in G\}\)
\end{center}
\begin{enumerate}
    \item \(\otp(C_G)=\lambda_0\).
    \item \textit{If \(\kappa_{i}<C_G(\gamma)<\kappa_{i+1}\) where \(\gamma\) is limit, then there exists \(\vec{\nu}=\langle\nu_1,...,\nu_m\rangle\) such that
    \((G)_{\vec{\nu}^{\frown}\langle\kappa_{i+1},...,\kappa_n\rangle}\)
    is generic for  \(\mathbb{M}_{\vec{\nu}^{\frown}\langle\kappa_{i+1},...,\kappa_n\rangle}[\vec{U}]\),  \(C_G=C_{(G)_{\vec{\nu}^{\frown}\langle\kappa_{i+1},...,\kappa_n\rangle}}\) and the sequences obtained by the split}
    \begin{center}
    \(\mathbb{M}_{\vec{\nu}}[\vec{U}]\times(\mathbb{M}_{\langle\kappa_{i+1},...,\kappa_n\rangle}[\vec{U}])_{>\nu_m}\simeq \mathbb{M}_{\vec{\nu}^{\frown}\langle\kappa_{i+1},...,\kappa_n\rangle}[\vec{U}]\)
    \end{center}
    \textit{are \(C_G\cap C_G(\gamma), C_G\setminus C_G(\gamma)\).
    More accurately, if}
    \begin{center}
\(\gamma=\underset{\xi}{\underbrace{\omega^{\gamma_1}+...+\omega^{\gamma_i}}}+\omega^{\gamma'_{i+1}}+...+\omega^{\gamma'_m}\) \ \ \textit{(C.N.F)}
\end{center} 
\textit{then} 
\begin{center}
\(\vec{\nu}=\langle\nu_1,...,\nu_m\rangle=\langle\kappa_1,...,\kappa_i,C_G(\xi+\omega^{\gamma'_{i+1}}),...,C_G(\gamma)\rangle\)
\end{center}
\end{enumerate}

\end{lemma}
\pr For (1), the same reasoning as in lemmas \ref{Order type of generic club}-\ref{Cor of changing cofinality} should work. For (2),  notice that by proposition \ref{Reconstruction of generic}, \(0_{\mathbb{M}_{\vec{\nu}^{\frown}\langle\kappa_{i+1},...,\kappa_n\rangle}}\in G\).
Thus \((G)_{\vec{\nu}^{\frown}\langle\kappa_{i+1},...,\kappa_n\rangle}\) is generic for \(\mathbb{M}_{\vec{\nu}^{\frown}\langle\kappa_{i+1},...,\kappa_n\rangle}[\vec{U}]\). The embeddings
\begin{center}
\(i_1:\mathbb{M}_{\langle\nu_1,...,\nu_m\rangle}[\vec{U}]\rightarrow\mathbb{M}_{\vec{\nu}^{\frown}\langle\kappa_{i+1},...,\kappa_n\rangle}[\vec{U}]\)\\
\(i_1(\langle t_1,...,t_{r+1}\rangle)=\langle t_1,...,t_{r+1},\langle\kappa_{i+1},B_{i+1}\rangle,...,\langle\kappa_n,B_n\rangle\rangle\)\\
\end{center}
and
\begin{center}
\(i_2:(\mathbb{M}_{\langle\kappa_{i+1},...,\kappa_n\rangle}[\vec{U}])_{>\nu_m}\rightarrow\mathbb{M}_{\vec{\nu}^{\frown}\langle\kappa_{i+1},...,\kappa_n\rangle}[\vec{U}]\)\\
\(i_2(\langle s_1,...,s_{k+1}\rangle)=\langle\langle\kappa_1,B_{1}\rangle,...,\langle\kappa_i,B_i\rangle ,s_1,...,s_{k+1}\rangle\)
\end{center} 
induces the isomorphism of \(\mathbb{M}_{\vec{\nu}^{\frown}\langle\kappa_{i+1},...,\kappa_n\rangle}[\vec{U}]\) with the product. Therefore,  \(i_1^{-1}(G)\), \( i_2^{-1}(G)\) are generic for \(\mathbb{M}_{\langle\nu_1,...,\nu_m\rangle}[\vec{U}]\),\ \((\mathbb{M}_{\langle\kappa_{i+1},...,\kappa_n\rangle}[\vec{U}])_{>\nu_m}\) respectively. By the definition of \(i_1,i_2\) this generics obviously yield the sequences \(C_G\cap C_G(\gamma)\) and \(C_G\setminus C_G(\gamma)\).
\begin{flushright}
 \(\blacksquare\)
\end{flushright}
In general we will identify \(G\) with \((G)_{\vec{\nu}}\) when using lemma \ref{Properties of CG}. \par
Notice that, the information used in order to compute \(\gamma(t_i,p)\) is just \(o^{\vec{U}}(t_i)\) which is provided by the suitable extension type. Let \(X\) be an extension type of \(p\), one can compute the coordinates of any extension \(\vec{\alpha}\) of type \(X\). In particular, for any \(\alpha_{i,r}\) substituting \(x_{i,r}\in X\) the coordinate of  \(\alpha_{i,r}\) is
 $$\gamma=\gamma(t_{i-1},p)+\omega^{x_{i,1}}+...+\omega^{x_{i,r}}=:\gamma(x_{i,r},p^{\frown}X)$$
 
 In this situation we say that \textit{\(X\) unveils the \(\gamma\)-th  coordinate}. If \(x_{i,r}=x_{mc}\), we say that \textit{X unveils \(\gamma\) as maximal coordinate}. 
\begin{proposition}\label{coordinates in the generic}
 
 \textit{Let \(p=\langle t_1,...,t_n,t_{n+1}\rangle\in \Mforc\) and \(\gamma\) such that for some \(0\leq i\leq n\),  \(\gamma(t_i,p)<\gamma<\gamma(t_{i+1},p)\). Then there exists an extension-type \(X\) unveiling \(\gamma\) as maximal coordinate. Moreover, if}
$$\gamma(t_i,p)+\underset{j\leq m}{\sum}\omega^{\gamma_j}=\gamma \  (C.N.F)$$
\textit{then the extension type is 
\(X=\langle X_i\rangle\) where \( X_i=\langle\gamma_1,...,\gamma_m\rangle\).}
\begin{flushright}
 \(\blacksquare\)
 \end{flushright}
 
\end{proposition}
 \textbf{Example:} Assume \(\lambda_0=\omega_1+\omega^2\cdot 2+\omega\), let \(\kappa_1<\kappa_2<\kappa_3<\kappa_4=\kappa\) be such that $$o^{\vec{U}}(\kappa_1)=\omega_1
, o^{\vec{U}}(\kappa_2)=o^{\vec{U}}(\kappa_3)=2\text{ and } o^{\vec{U}}(\kappa)=1$$ Let
$$p=\langle\underset{t_1 }{\underbrace{\scriptstyle\langle \nu_1,B(\nu_1)\rangle}},\underset{t_2 }{\underbrace{ \scriptstyle\nu_2}},\underset{t_3 }{\underbrace{\scriptstyle\langle \kappa_1,B(k_1)\rangle}},\underset{t_4}{\underbrace{\scriptstyle\langle \nu_4,B(\nu_3)\rangle}},\underset{t_5}{\underbrace{\scriptstyle\langle \kappa_2,B(\kappa_2)\rangle}},\underset{t_6}{\underbrace{\scriptstyle\langle \kappa_3,B(\kappa_3)\rangle}},\underset{t_7}{\underbrace{\scriptstyle\langle \kappa,B\rangle}}\rangle$$
 $$o^{\vec{U}}(t_1)=\omega ,\ o^{\vec{U}}(t_2)=0,\  o^{\vec{U}}(t_4)=1$$
 
Let \(G\) be any generic with \(p\in G\). Calculating \(\gamma(t_i,p) \) for \(i=1,...,7\):
\begin{enumerate}
\item \(\gamma(t_1,p)=\omega^{o^{\vec{U}}(t_1)}=\omega^{\omega}\ \Rightarrow \ \ C_G(\omega^{\omega})=\nu_1\).
\item \(\gamma(t_2,p)=\omega^{\omega}+\omega^{{o^{\vec{U}}(t_2)}}=\omega^{\omega}+1 \ \Rightarrow \ \ C_G(\omega^{\omega}+1)=\nu_2\).
\item \(\gamma(t_3,p)=\omega^{\omega}+1+\omega^{\omega_1}=\omega^{\omega_1}=\omega_1\).
\item \(\gamma(t_4,p)=\omega_1+\omega \ \Rightarrow \ \ C_G(\omega_1+\omega)=\nu_3\).
\item \(\gamma(t_5,p)=\omega_1+\omega+\omega^2=\omega_1+\omega^2\).
\end{enumerate}
To demonstrate proposition \ref{coordinates in the generic} let \(\gamma=\omega^{\omega}+\omega^5\cdot3+5\), then,
 $$\gamma(t_2,p)=\omega^{\omega}+1<\gamma<\omega_1=\gamma(t_3,p)$$
 $$(\omega^{\omega}+1)+\omega^5\cdot3+5=\gamma$$
The extension type unveiling \(\gamma\) as maximal coordinate is then 
\begin{center}
\(X=\langle\langle\rangle,\langle\rangle, X_3\rangle, \ X_3=\langle5,5,5,0,0,0,0,0\rangle\)
\end{center}
i.e. every extension \(\vec{\alpha}=\langle \alpha_{3,1},...\alpha_{3,8}\rangle\in B(p,X)\) will satisfy that 
\begin{center}
 
\(\gamma(\alpha_{mc},p^{\frown}\vec{\alpha})=\gamma(\alpha_{3,8},p^{\frown}\alpha)=\gamma(x_{3,8},p^{\frown}X)=\gamma\)
\end{center}

Let us state the main theorem of this paper.
\begin{theorem}\label{Main theorem}
\textit{Let \(\vec{U}\) be a coherent sequence in \(V\),  \(\langle\kappa_1,...\kappa_n\rangle\) be a sequence such that \(o^{\vec{U}}(\kappa_i)<\min(\nu\mid 0<o^{\vec{U}}(\nu))=:\delta_0\), let \(G\) be \(\Mforc\)-generic and let \(A\in V[G]\) be a set of ordinals. Then there exists \(C'\subseteq C_G\) such that \(V[A]=V[C']\).}

\end{theorem} 
We will prove Theorem \ref{Main theorem} by induction on \(\otp(C_G)\). For \(\otp(C_G)=\omega\) it is just the Prikry forcing which is know by \cite{PrikryCase}. Let \(\otp(C_G)=\lambda_0\) be a limit ordinal,
\begin{center}
\(\lambda_0=\omega^{\gamma_n}+...+\omega^{\gamma_1}\) \ (C.N.F)
\end{center}

If \(\sup(A)<\kappa\), then by lemma 5.3 in \cite{ChangeCofinality}, \(A\in V[C_G\cap \sup(A)]\). By lemma \ref{Properties of CG}, \(V[C_G\cap \sup(A)]\) is a generic extension of some \(\mathbb{M}_{\langle \nu_1,...,\nu_m\rangle}[\vec{U}]\) with order type smaller the \(\lambda_0\), thus by induction we are done. In fact, if there exists \(\alpha<\kappa\) such that \(A\in V[C_G\cap\alpha]\) then the induction hypothesis works. Let us assume that \(A\notin V[C_G\cap\alpha]\) for every \(\alpha<\kappa\)  this kind of set will be called \textit{recent set}. Since \(\kappa_1,...,\kappa_n\) will be fixed through the rest of this chapter we shall abuse notation and denote \(\Mforc=\Mfor\). 

\subsection{The Main Results for Sets of Cardinality Less Than $\kappa$}
First let us show that for \(A\) with small enough cardinality the theorem holds regardless of the induction. 
\begin{lemma}\label{singleExtension}
\textit{ Let \( \underaccent{\sim}{x}\) be a \(\Mfor\)-name and \(p \in \) \(\Mfor\) such that \(p \Vdash \underaccent{\sim}{x} \) is an ordinal. Then there exists \(p\leq^*p^*\in \Mfor\) and an extension-type \(X\in T(p)\) such that }
\begin{center}

   \((*)\)\ \ \ \ \ \ \ \ \ \ \ \ \ \ \(\forall p^{*\frown}\langle\vec{\alpha}\rangle\in p^{*\frown}X\)\ \  \(p^{*\frown}\langle\vec{\alpha}\rangle || \  \underaccent{\sim}{x}\)
    
\end{center}

\end{lemma}
\pr
 Let \(p=\langle t_1,...,t_n,t_{n+1}\rangle\in \Mfor\).
 \begin{claim}
  \textit{There exists \(p\leq^* p'\)  such that for some extension type \(X\)}
 \begin{center}
 \(\forall \vec{\alpha}\in B(p',X)\) \(\exists C(x_{i,j})\) s.t. \(p'^{\frown}\langle\vec{\alpha},(C(x_{i,j}))_{i,j}\rangle \ || \ \underaccent{\sim}{x}\)
 \end{center}
 \end{claim}
 \textit{Proof of Claim:} Define sets \( B_X(t_i,j)\), for any fixed \(X\in T(p)\) as follows: 
 Recall the notation \( l_X , \ x_{mc}\) and let
 \(\vec{\alpha}\in B(
 p,X\setminus \langle x_{mc} \rangle)\).
 Define 
 \begin{center}
 \(B^{(0)}_{X}(\vec{\alpha})= \{ \theta \in B(t_{l_X},x_{mc} ) \mid \exists (C(x_{i,j}))_{i,j}  \ \ p^{\frown}\langle\vec{\alpha},\theta,(C(x_{i,j}))_{i,j}\rangle || \  \underaccent{\sim}{x} \} \)
 \end{center}
and \( B^{(1)}_{X}(\vec{\alpha})=B(t_{l_X},x_{mc}) \setminus B^{(0)}_{X}(\vec{\alpha})\). One and only one of \(B^{(0)}_{X}(\vec{\alpha})\ ,B^{(1)}_{X}(\vec{\alpha}) \) is in \(U(\kappa(t_{l_X}),x_{mc})\).
 Set \(B_X(\vec{\alpha})\) and \( F_X(\vec{\alpha}) \in \{0,1\} \) such that
\begin{center}
\( B_X(\vec{\alpha})=B^{(F_X(\vec{\alpha}))}_{X}(\vec{\alpha}) \in U(\kappa(t_{l_X}),x_{mc})\)
 \end{center}
 Define 
 \begin{center}
 \(B'_X(t_{l_X},x_{mc})= \underset{\vec{\alpha} \in B(p,X\setminus\langle x_{mc}\rangle)}{\Delta}B_X(\vec{\alpha})\)
 \end{center}
 Consider the function \(F:B(p,X\setminus\langle x_{mc}\rangle) \rightarrow \{ 0,1\} \). Applying lemma \ref{stabprod} to \(F\), we get an homogeneous \(\underset{x_{i,j}\in X\setminus\langle x_{mc}\rangle}{\prod}B'_X(t_i,x_{i,j})\) where
 \begin{center}

 \( B'_X(t_i,x_{ij})  \subseteq B(t_i,x_{ij}), \ B'_X(t_i,x_{ij})\in U(t_i,x_{i,j}), \ x_{ij}\in X\setminus\langle x_{mc}\rangle \)
 \end{center}
 For \(\xi\notin X_i\), Set
 \begin{center}
 \( B'_X(t_i,\xi)=B(t_i,\xi)\)
 \end{center}
Since  \( |T(p)|<\kappa(t_1)\), for each \(1\leq i\leq n+1\) and \(\xi<o^{\vec{U}}(t_i)\)
 \begin{center}
 
\(B'(t_i,\xi):=\underset{X\in T(p)}{\bigcap}B'_X(t_i,\xi)\in U(\kappa(t_i),\xi)\) 
\end{center}
Finally, let  \(p'=\langle t'_1,...,t'_n,t_{n+1}'\rangle\) where
 \begin{center}
  \(t_i'=\begin{cases}
t_i & o^{\vec{U}}(t_i)=0\\
\langle\kappa(t_i),B'(t_i)\rangle & otherwise

\end{cases}\)
 \end{center}
 It follows that \(p\leq^* p'\in \Mfor\).
 
 Let \(H\) be \(\Mfor\)-generic, \( p' \in H\). By the assumption on \(p\), there exists \(\delta<\kappa\) such that \(V[H] \models (\underaccent{\sim}{x})_{H}=\delta\) . Hence , there is \(p' \leq q \in M[\vec{U}] \) such that 
 \(q \Vdash \underaccent{\sim}{x}=\overset{\vee}{\delta}\). By proposition \ref{partition of extension type}, there is a unique \(  p'^{\frown}\langle\vec{\alpha},\theta\rangle\in p'^{\frown}X \)   for some extension type X, such that \(p'^{\frown}\langle\vec{\alpha},\theta\rangle\leq^*q\). \(X,p'\) are as wanted: \\By the definition of \(p'\) it follows that \(\vec{\alpha} \in B(p',X\setminus\langle x_{mc}\rangle)\) and \(\theta \in B_X(\vec{\alpha})\). Since \(q \Vdash \underaccent{\sim}{x}=\overset{\vee}{\delta}\), we have that \(F_X(\vec{\alpha})=0\). Fix \(\langle\vec{\alpha'},\theta'\rangle\) of type X. \(\vec{\alpha'}\) and \(\vec{\alpha}\) belong to the same homogeneous set, thus \(F(\vec{\alpha'})=F(\vec{\alpha})=0\) and
 \begin{center}
 \(\theta' \in B^{(0)}_{X}(\vec{\alpha'})\Rightarrow\exists (C(x_{i,j}))_{i,j}\) \(s.t.\) \( p'^{\frown}\langle\vec{\alpha'},\theta',(C(x_{i,j}))_{i,j}\rangle || \  \underaccent{\sim}{x}\) 
   
 \end{center}
 \begin{flushright}
 \(\blacksquare_{of\ claim}\)
 \end{flushright}
 For every \(\vec{\alpha}\in B(p',X)\), fix some \((C_{i,j}(\vec{\alpha}))_{^{i\leq n+1}_{j\leq l_i+1}}\) such that 
 \begin{center}
 \(p'^{\frown}\langle \vec{\alpha},(C_{i,j}(\vec{\alpha}))_{^{i\leq n+1}_{j\leq l_i+1}}\rangle || \  \underset{\sim}{x}\)
 \end{center}
 It suffices to show that we can find \(p'\leq^* p^*\) such that for every \(\vec{\alpha}\in B(p^*,X)\)
 \begin{center}
 
\(B(t^*_i)\cap(\alpha_s,\alpha_{i,j})\subseteq C_{i,j}(\vec{\alpha}) \) , \ \ \(1\leq i\leq n+1 , \ 1\leq j\leq l_i+1\)

 \end{center}
 Where \(\alpha_s\) is the predecessor of \(\alpha_{i,j}\) in \(\vec{\alpha}\) .  
 In order to do that, define \(p'\leq^* p_{i,j}\)  \( i\leq n+1 , \ j\leq l_i+1\) then \(p^*\geq^*p_{i,j}\) will be as wanted. Define \(p_{i,j}\) as follows:\\
 Fix \(\vec{\beta}\in B(p',\langle x_{1,1},...,x_{i,j}\rangle)\), by lemma \ref{stabprod}, the function 
 \begin{center}
 \(C_{i,j}(\vec{\beta},*):B(p',X\setminus\langle x_{1,1},...,x_{i,j}\rangle)\rightarrow P(\beta_{i,j})\)
 \end{center}
  has homogeneous sets \(B^*(\vec{\beta},x_{r,s}) \subseteq B(p',x_{r,s})\) for \(x_{r,s}\in X\setminus\langle x_{1,1},...,x_{i,j}\rangle\). Denote the constant value by \(C^*_{i,j}(\vec{\beta})\). Define 
 \begin{center}
\(B^*(t_r,x_{r,s})=\underset{\vec{\beta}\in B(p',\langle x_{1,1},...,x_{i,j}\rangle)}{\Delta}B^*(\vec{\beta},x_{r,s})\), \ \(x_{r,s}\in X\setminus\langle x_{1,1},...,x_{i,j}\rangle\)
 \end{center}
 Next, fix \(\alpha\in B(t'_i,x_{i,j})\) and let
\begin{center}
\(C^*_{i,j}(\alpha)=\underset{\vec{\alpha'}\in B(p',\langle x_{1,1},...,x_{i,j-1}\rangle)}{\Delta}C^*_{i,j}(\vec{\alpha'},\alpha)\)
\end{center}
Thus \(C^*_{i,j}(\alpha)\subseteq\alpha\).  Moreover, \(\kappa(t_i)\) is in particular an ineffable cardinal and therefore there are \(B^*(t_i,x_{i,j})\subseteq B(t'_i,x_{i,j})\) and \(C^*_{i,j}\) such that 
\begin{center}
\(\forall \alpha\in B^*(t_i,x_{i,j}) \ \ \ C^*_{i,j}\cap\alpha=C^*_{i,j}(\alpha)\)
\end{center}
 By coherency,  \(C^*_{i,j}\in \bigcap U(t_i,\xi)\).
 Finally, define  \(p_{i,j}=\langle t^{(i,j)}_1,...,t^{(i,j)}_n,t^{(i,j)}_{n+1} \rangle\)
 \begin{center}
 \(B(t^{(i,j)}_i)=B^*(t_i)\cap(\underset{j}{\bigcap}C^*_{i,j})\) \ \ \ \ \ \(1\leq i\leq n+1\)

 \end{center}
 To see that \(p^*\) is as wanted, let \(\vec{\alpha}\in B(p^*,X)\) and fix any \(i,j\). Then \(\vec{\alpha}\in B(p_{i,j},X)\) and \(\alpha_{i,j}\in B^*(t_i,x_{i,j})\). Thus
 \begin{center}
 \(B(t^*_i)\cap(\alpha_s,\alpha_{i,j})\subseteq C^*_{i,j}\cap\alpha_{i,j}\setminus\alpha_s=C^*_{i,j}(\alpha_{i,j})\setminus\alpha_s\subseteq C^*_{i,j}(\alpha_{1,1},...,\alpha_{i,j})=C_{i,j}(\vec{\alpha})\)
 \end{center}
 \begin{flushright}
 \(\blacksquare\)
 \end{flushright}
\begin{lemma}\label{Thm for small set}
 \textit{Let \(G\subseteq\Mfor\) be $V$-generic filter and \(A\in V[G]\) be any set of ordinals, such that  \(|A|<\delta_0\). Then there is \(C' \subseteq C_G\) such that \(V[A]=V[C']\).}
 
\end{lemma}
 \pr Let  \(A=\langle a_{\xi} \mid \xi < \delta\rangle \in V[G]\) , where \(\delta<\delta_0\) and \(\underaccent{\sim}{A}=\langle\underaccent{\sim}{a}_{\xi} \mid \xi < \delta\rangle\)  be a $\Mfor$-name for \(\langle a_{\xi} \mid \xi < \delta\rangle
 \). Let \(q\in G\) such that \(q\Vdash \underaccent{\sim}{A}\subseteq On\). We proceed by a density argument, fix \(q\leq p\in \) \(\Mfor\). By lemma \ref{singleExtension}, for each \( \xi < \delta\) there exists \(X(\xi)\) and \(p\leq^* p^*_\xi\) satisfying \((*)\).
 By \(\delta^+-_{\leq^*}\)closure above \(p\) we have \(p^*\in\Mfor\) such that \(\forall\xi<\delta \ p^*_\xi\leq p^*\).
 For each \(\xi\), define \(F_{\xi}:B(p^*,X(\xi))\longrightarrow \kappa\)
 \begin{center}
 \(F_{\xi}(\vec{\alpha})=\gamma\) for the unique \(\gamma\) such that \(p^{*\frown}\langle\vec{\alpha}\rangle \Vdash   \underaccent{\sim}{a}_{\xi}=\overset{\vee}{\gamma}\).
 \end{center}
 Using lemma \ref{ImportantCoor}, we obtain for every \( \xi<\delta\) a set of important coordinates
 \begin{center}
\(I_{\xi}\subseteq\{\langle i,j\rangle\mid1\leq i\leq n+1\ , \ 1\leq j\leq l_i \}\)
 \end{center}
\underline{Example}: Assume \(o^{\vec{U}}(k)=3\)  ,\(C_G=\langle C_G(\alpha) \mid \alpha<\omega^3\rangle\).
\begin{center}
\(a_0=C_G(80) , a_1=C_G(\omega+2)+C_G(3) ,a_2=C_G(\omega^2\cdot 2+\omega+1)\)
\end{center}
and
\begin{center}
\(p=\langle \nu_0,\langle\nu_{\omega},B(\nu_{\omega},0)\rangle,\langle\kappa,\underset{B(\kappa)}{\underbrace{B(\kappa,0)\cup B(\kappa,1)\cup B(\kappa,2)}} \rangle\rangle\)
\end{center}
 We use as index the coordinate in the final sequence to improve clarity. To determine \(a_0\), unveil the first 80 elements of the Magidor sequence i.e. any element of the form
\begin{center}\(p_0=\langle\nu_0,\nu_1,...,\nu_{80},\langle\nu_{\omega},B(\nu_{\omega},0)\setminus \nu_{80}+1\rangle,\langle\kappa,B(\kappa)\rangle \rangle\)
\end{center}
will decide the value of \(a_0\). Thus the extension type X(0) is
\begin{center}
 \(X(0)=\langle\langle\underset{80 \ times}{\underbrace{0,...,0}}\rangle,\langle\rangle\rangle\)
\end{center}
The important coordinates to decide the value of \(a_0\) is only the 80th coordinate. It is easily seen to be one to one modulo the irrelevant coordinates $1,...,79$. For \(a_1\), the form is
\begin{center}\(p_1=\langle\nu_0,\nu_1,\nu_2,\nu_3,\langle\nu_{\omega},B(\nu_{\omega},0)\setminus \nu_3+1\rangle,\nu_{\omega+1},\nu_{\omega+2},\langle\kappa,B(\kappa)\setminus (\nu_{\omega+2}+1)\rangle\rangle\)
\end{center}
The extension type is
\begin{center}
\( X(1)=\langle\langle0,0,0\rangle,\langle0,0\rangle\rangle\)
\end{center}
The important coordinates are the 3rd and the 5th.
For \(a_2\) we have
$$\scriptstyle p_2=\langle\nu_0,\langle\nu_{\omega},B(\nu_{\omega},0)\rangle,\langle\nu_{\omega^2},B(\nu_{\omega^2})\rangle,\langle\nu_{\omega^2\cdot2},B(\nu_{\omega^2\cdot2})\rangle,\langle\nu_{\omega^2\cdot2+\omega},B(\nu_{\omega^2\cdot2+\omega})\rangle,\langle\kappa,B(\kappa) \setminus \nu_{\omega^2\cdot2+\omega}\rangle\rangle$$
\begin{center}
 \(X(2)=\langle\langle\rangle,\langle2,2,1\rangle\rangle\)
\end{center}
 Back to the proof, since \(G\) was generic, there is \(\langle t_1,...,t_n,t_{n+1}\rangle=p^{\star} \in G\) with such functions. 
 Find \(D_\xi \subseteq C_{G}\) such that
  \begin{center}
  \(D_\xi\in B(p^{\star},X_\xi)\)
  \end{center}
   By proposition \ref{Reconstruction of generic} and since \(p^{\star} \in G\), \(D_\xi\) exists. 
  Since \(V[G]\models (\underaccent{\sim}{a}_{\xi})_G=a_\xi\) we conclude that,
$$p^{\star\frown} \langle D_{\xi}\rangle \Vdash\underaccent{\sim}{a}_{\xi}=\overset{\vee}{a_{\xi}}$$
hence,
\(F_\xi(D_\xi)=a_\xi\).
 Set
 \(C_\xi=D_\xi\restriction I_\xi\)
 and \(C'=\underset{\xi<\delta}{\bigcup}C_{\xi}\). Let us show that \(V[\langle a_{\xi} | \xi < \delta\rangle]=V[C']\): \\
 \\
 In \(V[C']\), fix some enumeration of \(C'=\{C'_i\mid i<\otp(C')\}\). For each $\xi<\delta$, $C_\xi$ can be extracted from \(C'\) and \( \Inx(C_\xi,C')\in V\) (See the notation section for the definition of \Inx(A,B)). Since $\delta<\delta_0$, $\l\Inx(C_\xi,C')\mid\xi<\delta\r\in V$, which implies that $\l C_\xi\mid\xi<\delta\r\in V[C']$. Still in $V[C']$, for every \(\xi<\delta\) find
 \begin{center}
 \(D_\xi' \in B(p^{\star},X_\xi)\)
 such that
 \(D'_\xi\restriction I_\xi=C_\xi\)

 \end{center}
 Such \(D'_\xi\) exists as \(D_\xi\) witnesses (the sequence \(\langle D_\xi \mid \xi<\delta\rangle\) may not be in \(V[C']\)). Since \(D'_\xi\sim_{I_\xi} D_{\xi} \), and by the property of $I_\xi$,
 \begin{center}
 \(F_{\xi}(D'_\xi)=F_{\xi}(D_{\xi})=a_{\xi}\)
 
 \end{center}
 hence \(\langle a_{\xi} \mid \xi<\delta\rangle=\langle F_{\xi}(D'_\xi) \mid \xi<\delta\rangle \in V[C']\).\\
 In the other direction, Given \(\langle a_{\xi} \mid \xi<\delta\rangle\), for each \( \xi<\delta \) pick  \(D'_{\xi} \in F_{\xi}^{-1}(a_{\xi}) \) (Note that \(F_{\xi}^{-1}(a_{\xi}) \neq \emptyset\) follows from the fact that \(D_{\xi} \in dom(F_{\xi})\) and \(F_{\xi}(D_{\xi})=a_{\xi}\)). Since \(F_{\xi}\) is 1-1 modulo \(I_{\xi}\) and \(F_{\xi}(D_{\xi})=F_{\xi}(D'_{\xi})\) we have
 \begin{center}
 \(D_{\xi}\sim_{I_\xi} D'_{\xi}\) and \(C_{\xi}=D_\xi\restriction I_\xi=D'_\xi\restriction I_{\xi}\)
 \end{center}
 Hence 
 \begin{center}
 \(\langle C_{\xi} \mid\xi <\delta\rangle=\langle D'_\xi\restriction I_{\xi} \mid\xi < \delta\rangle \in V[\langle a_{\xi} \mid \xi < \delta\rangle]\)
 and \(C'\in V[\langle a_{\xi} \mid \xi < \delta\rangle]\).
   \end{center}
 \begin{flushright}
 \(\blacksquare\)
 \end{flushright}
 \subsection{The Main Result for Subsets of $\kappa$}
We shall proceed by induction on  \(\sup(A)\) for a recent set \(A\). As we have seen in the discussion following Theorem \ref{Main theorem}, if \(A\subseteq\kappa\) is recent then \(\sup(A)=\kappa\). For such 
\(A\), the next lemma gives a sufficient conditions.
\begin{lemma}\label{sufficient condition}
 \textit{Let \( A \in V[G] ,\ \sup(A)=\kappa \).  Assume that \(\exists C^*\subseteq C_G\) such that}
\begin{enumerate}
 \item\(C^*\in V[A]\) \(and\)
 \(\forall\alpha<\kappa \  A\cap\alpha\in V[C^*]\).
 \item \( cf^{V[A]}(\kappa)<\delta_0\).
 \end{enumerate}
 \textit{Then \( \exists C' \subseteq C_G \) such that \( V[A]=V[C']\).}
 
\end{lemma}
 \pr Let \(cf^{V[A]}(\kappa)=\eta\) and \( \langle\gamma_{\xi} \mid \xi < \eta\rangle \in V[A]\) be a cofinal sequence in \( \kappa\). Work in \(V[A]\) , pick an enumerations of \(P(\gamma_{\xi})=\langle X_{\xi,i} \mid i<2^{\gamma_{\xi}}\rangle\in V[C^*]\). Since \(A\cap\gamma_{\xi} \in V[C^*]\), there exists \(i_{\xi}<2^{\gamma_{\xi}}\) such that \(A\cap\gamma_{\xi}=X_{\xi,i_{\xi}}\). The sequences
 \begin{center}
 \(C^*, \ \ \langle i_{\xi} \mid\xi<\eta\rangle, \ \ \langle\gamma_{\xi} \mid \xi < \eta\rangle\) 
 
 \end{center}
 can be coded in \(V[A]\) to a sequence \(\langle x_{\alpha} \mid \alpha < \eta\rangle\). By lemma \ref{Thm for small set}, \( \exists C' \subseteq C_G \) such that \( V[\langle x_{\alpha} \mid \alpha < \eta\rangle]=V[C']\). Let us argue that \(V[A]=V[\langle x_{\alpha} \mid \alpha < \delta\rangle]\), clearly
 \(V[A] \supseteq V[\langle x_{\alpha} \mid \alpha < \eta\rangle]\). For the other direction, note that \(A=\underset{\xi<\eta}{\bigcup}X_{\xi,i_{\xi}}\in V[\langle x_{\alpha} \mid \alpha < \eta\rangle]\).
 \begin{flushright}
 \(\blacksquare\)
 \end{flushright}
Let us consider two kind of subsets of $\kappa$:
\begin{enumerate}
    \item \(\exists \alpha^*<\kappa\) such that \(\forall \beta<\kappa \ \ A\cap\beta\in V[A\cap\alpha^*]\) and we say that \(A\cap\alpha\) stabilizes.
    An example of such $A$ is a generic Prikry sequence $\{C_G(n)\mid n<\omega\}$, simply take \(\alpha^*=0\).
    \item For all \(\alpha<\kappa\) there exists \(\beta<\kappa\) such that \(V[A\cap\alpha]\subsetneq V[A\cap\beta]\) as example we can take Magidor forcing with \(o^{\vec{U}}(\kappa)=2\) and \(A\) the entire Magidor sequence \(C_G\).
    
\end{enumerate}
First we consider \(A\)'s such that \(A\cap\alpha\) does not stabilize.
\begin{lemma}
\textit{ Assume that \(A\cap\alpha\) does not stabilize, then there exists \(C'\subseteq C_G\) such that \(V[A]=V[C']\).}
\end{lemma}
\pr 
Work in \(V[A]\), define the sequence \(\langle\alpha_\xi \mid \xi<\theta\rangle\):
\begin{center}
\( \alpha_0=\min(\alpha \mid V[A\cap\alpha]\supsetneq V)\)
\end{center}
Assume that \(\langle \alpha_\xi \mid \xi<\lambda\rangle\) has been defined and for every \(\xi ,\ \alpha_\xi<\kappa\). 
If \(\lambda=\xi+1\) then set
\begin{center}
\(\alpha_\lambda=\min(\alpha \mid V[A\cap\alpha]\supsetneq V[A\cap\alpha_\xi])\)
\end{center}

If \(\alpha_\lambda=\kappa\) , then \(\alpha_{\lambda}\) satisfies that
\begin{center}
\(\forall \alpha<\kappa \ \ A\cap\alpha\in V[A\cap\alpha_{\lambda^*}]\)

\end{center}
Thus \(A\cap\alpha\) stabilizes which contradicts our assumption.

If \(\lambda\) is limit, define 
\begin{center}
\(\alpha_\lambda=\sup(\alpha_\xi \mid \xi<\lambda)\)
\end{center}
 if \(\alpha_\lambda=\kappa\) define \(\theta=\lambda\) and stop.  The sequence \(\langle \alpha_\xi \mid \xi<\theta\rangle\in V[A]\) is a continues, increasing unbounded sequence in \(\kappa\). Therefore, \( cf^{V[A]}(\kappa)=cf^{V[A]}(\theta)\). Let us argue that \(\theta<\delta_0\). Work in \(V[G]\), for every \(\xi<\theta \) pick \(C_\xi\subseteq C_G\) such that \(V[A\cap\alpha_\xi]=V[C_\xi]\). This is a 1-1 function from \(\theta\) to \(P(C_G)\). The cardinal \(\delta_0\) is still a strong limit cardinal (since there are no new bounded subsets below this cardinal and it is measurable in \(V\)). Moreover, \(\lambda_0:=\otp(C_G)<\delta_0\), thus 
\begin{center}
\(\theta\leq|P(C_G)|=|P(\lambda_0)|<\delta_0\) 
\end{center}
The only thing left to prove, is that we can find \(C^*\) as in Lemma \ref{sufficient condition}. Work in \(V[A]\), for every \( \xi<\theta\),  \(C_\xi\in V[A]\) (The sequence \(\langle C_\xi \mid \xi<\theta\rangle\) may not be in \(V[A]\)). \(C_\xi\) witnesses that
\begin{center}
 \( \exists d_\xi\subseteq \kappa\)  \((|d_\xi|<2^{\lambda_0}\) and \(V[A\cap\alpha_\xi]=V[d_\xi]\))
\end{center}
So \(d=\bigcup\{d_{\xi} | \xi<\theta\}\in V[A]\) and \( |d|\leq2^{\lambda_0}\). Finally, by lemma \ref{Thm for small set}, there exists \(C^*\subseteq C_G\) such that \(V[C^*]=V[d]\subseteq V[A]\). Note that for every $\xi<\theta$, $\Inx(d_\xi,d)\in V$ and also since $\theta<\delta_0$, the sequence $\l \Inx(d_\xi,x)\mid\xi<\theta\r\in V$. It follows that for every $\xi<\theta$, $d_\xi\in V[C^*]$, and in turn $A\cap\alpha_\xi\in V[C^*]$. Since $\alpha_\xi$ is unbounded in $\kappa$, for all \(\alpha<\kappa \ A\cap\alpha\in V[C^*]\). Apply \ref{sufficient condition}, to conclude the lemma. 
\begin{flushright}
 \(\blacksquare\)
 \end{flushright}
 For the rest of this section, we assume that the sequence \(A\cap\alpha\) stabilizes on \(\alpha^*\). Let \(C^*\) be such that \(V[A\cap\alpha^*]=V[C^*]\) and \(\kappa^*=\sup(C^*)\) is limit in \(C_G\). Notice that, \(\kappa^*<\kappa\), since if $\kappa^*=\kappa$, then $\kappa$ is singular in $V[C^*]$, but on the other hand \(A\cap\alpha^*\in V[C_G\cap\alpha^*]\) which implies $\kappa$ is regular in $V[A\cap\alpha^*]=V[C^*]$.
 
In order to apply lemma \ref{sufficient condition}, we only need to argue that for $A$ which is recent, \(\kappa\) changes cofinality in \(V[A]\). To do this, consider the initial segment \(C_G\cap\kappa^*\) and assume that  \(\kappa_{j-1}\leq\kappa^*<\kappa_j\).
Denote by
$$\mathbb{M}_{\leq\kappa^*}:=\mathbb{M}_{\langle\nu_1,...,\nu_i,\kappa^*\rangle}[\vec{U}], \  \mathbb{M}_{>\kappa^*}[\vec{U}]:=(\mathbb{M}_{\langle\kappa_j,...,\kappa\rangle}[\vec{U}])_{>\kappa^*}$$
By lemma \ref{Properties of CG} we can split \(\Mfor\) to
\begin{center}
\(\mathbb{M}_{\leq\kappa^*}[\vec{U}]\times(\mathbb{M}_{\langle\kappa_j,...,\kappa\rangle}[\vec{U}])_{>\kappa^*}\)

\end{center}
such that \(C_G\) is generic for \(\mathbb{M}_{\leq\kappa^*}[\vec{U}]\times\mathbb{M}_{>\kappa^*}[\vec{U}]\) and \(C_G\cap\kappa^*\) is generic for \(\mathbb{M}_{\leq\kappa^*}[\vec{U}]\). By \cite[Thm. 15.43]{Jech}, there is a forcing \(\mathbb{P}\subseteq RO(\Mfor_{\leq\kappa^*})\) \footnote{ $RO(\mathbb{Q})$ is the complete boolean algebra of regular open subsets of $\mathbb{Q}$}, such that \(V[C^*]=V[G^*]\) for some generic \(G^*\) of \(\mathbb{P}\). Also there is a projection of $\pi:\mathbb{M}_{\leq\kappa^*}[\vec{U}]\rightarrow\mathbb{P}$. Recall that if \(\pi:\mathbb{M}_{\leq\kappa^*}[\vec{U}]\rightarrow\mathbb{P}\) is the projection, then the quotient forcing is define: 
\begin{center}
\(\mathbb{M}_{\leq\kappa^*}[\vec{U}]/G^*=\pi^{-1''}[G^*]\)
\end{center}
In \(V[G^*]\) define \(\mathbb{Q}=\mathbb{M}_{\leq\kappa^*}[\vec{U}]/G^*\subseteq\mathbb{M}_{\leq\kappa^*}[\vec{U}]\). It is well known that \(G\restriction\kappa^*\) is \(V[C^*]\)-generic filter for \(\mathbb{Q}\) and clearly \(V[C^*][C_G\cap\kappa^*]=V[C_G\cap\kappa^*]\). In section $4$, we give a more concrete description of \(\pi\) and \(\mathbb{Q}\), however, in this section we will only need the existence of such a forcing and the fact that the projection if on the part below \(\kappa^*\) which implies that \(\mathbb{Q}\) is of small cardinality. 

Forcing \(\mathbb{M}_{>\kappa^*}[\vec{U}]\) above $V[G\restriction\kappa^*]$ is essentially forcing a Magidor forcing adding a sequence to $\kappa$ above $\kappa^*$. To see this, note that all the measures in \(\vec{U}\) above \(\kappa^*\) generates measures in $V[G\restriction\kappa^*]$. In conclusion, we have managed to find a forcing \(\mathbb{Q}\times\mathbb{M}_{>\kappa^*}[\vec{U}]\in V[C^*]\) such that \(V[G]\) is one of its generic extensions and \(\forall\alpha<\kappa \ A\cap\alpha\in V[C^*]\).\\
\\
Work in \(V[C^*]\), let \(\underaccent{\sim}{A}\) be a \(\mathbb{Q}\times\mathbb{M}_{>\kappa^*}[\vec{U}]\)-name for $A$. 
Since $A$ stabilizes, and by the definition of \(C^*\), we can find \(\langle q,p\rangle\in G\) such that $$\langle q,p\rangle\Vdash \forall\alpha<\kappa \  \underaccent{\sim}{A}\cap\alpha\text{
is old (where old means in } V[C^*])$$ Formally, the next argument is a density argument above \(\langle q,p\rangle\).
Nevertheless, in order to simplify notation, assume that \(\langle q,p\rangle=0_{Q\times \Mfor_{>\kappa^*}}\).
 Lemmas \ref{Finding aprox.}-\ref{ineffability} prove that a certain property holds densely often in \(\Mfor_{>\kappa^*}\). In order to Make these lemmas more clear, we consider an ongoing example.\\
 \underline{Example}: Let \(\lambda_0=\otp(C_G)=\omega^2\),
 \begin{center} 
  \(A=\{C_G(n)\mid n\leq\omega \text{ is even}\}\cup\{C_G(\omega\cdot n)+C_G(n)\mid 0<n<\omega\}\)
 \end{center}
 Therefore
 \begin{center}
 
\(C^*=\{C_G(2n)\mid n<\omega\}, \  \kappa^*=C_G(\omega)\)
 \end{center}
 The forcing \(\mathbb{Q}\) can be thought of as adding the missing coordinates to \(C_G\restriction\omega\) i.e. the odd coordinates. For the sake of the example, let
 \begin{center}
  \(p=\langle\underset{t_1}{\underbrace{\langle\nu_{\omega\cdot 2}, B_{\omega\cdot 2}\rangle}},\underset{t_2}{\underbrace{\nu_{\omega\cdot2+1}}},\underset{t_3}{\underbrace{\langle\kappa,B(\kappa)\rangle}}\rangle\in \Mfor_{>\kappa^*} \)
 \end{center}
\begin{lemma}\label{Finding aprox.}
 \textit{For every \(p\in \Mfor_{>\kappa^*}\) there exists \(p\leq^*p^*\) such that for every  extension type \(X\)  of \(p^*\) and \(q\in\mathbb{Q}\) (Recall  \(\max(\vec{\alpha})\alpha_{mc}\), if there is \(p^{*\frown}\vec{\alpha}\in p^{*\frown}X \) and \(p^{**}\geq^*p^{*\frown}\vec{\alpha}\)} such that \( \langle q,p^{**}\rangle|| \underaccent{\sim}{A}\cap\alpha_{mc})\Rightarrow\), then
\begin{center}
 \((*)\)  \ \ \ \ \ \ \ \ \ \((\forall p^{*\frown}\vec{\alpha}\in p^{*\frown} X \  \langle q,p^{*\frown}\vec{\alpha}\rangle|| \underaccent{\sim}{A}\cap\alpha_{mc}=:a(q,\vec{\alpha}))\)  \textit{(a propery of \(q,X\))}
 \end{center} 

\end{lemma} 
\textbf{Example:} Let
$$q=\langle\nu_1,\nu_3,\langle\kappa^*,B(\kappa^*)\rangle\rangle, \  X=\langle\underset{X_1}{\underbrace{\langle0,0\rangle}},\underset{X_2}{\underbrace{\langle \ \rangle}},\underset{X_3}{\underbrace{\langle 1,0\rangle}}\rangle-\text{ extension of }p$$
and let $$\vec{\alpha}=\langle\langle\alpha_{\omega+1},\alpha_{\omega+2}\rangle,\langle\ \rangle ,\langle\alpha_{\omega\cdot3},\alpha_{\omega\cdot3+1}\rangle\rangle\in B(p,X)$$
 If \(H\) is any generic with \(\langle q,p^{\frown}\langle\vec{\alpha}\rangle\rangle\in H\) then
all the elements in \(q\) and \(p^{\frown}\langle\vec{\alpha}\rangle\) have there coordinates in \(C_H\) as specified above, thus, $$(\underaccent{\sim}{A})_H\cap\alpha_{mc}=(\underaccent{\sim}{A})_H\cap\alpha_{\omega\cdot3+1}=\{C_H(n)\mid n\leq\omega\text{ is even}\}\cup\{C_H(\omega\cdot n)+C_H(n)\mid n<\omega\}\cap C_H(\omega\cdot3+1)$$
If \(\alpha_{\omega\cdot3}+\nu_3\geq \alpha_{\omega\cdot3+1}\) then
\begin{center}
\(a(q,\vec{\alpha})=(\underaccent{\sim}{A})_H\cap\alpha_{mc}=C_H\restriction_{even}\cup\{C_H(\omega),C_H(\omega)+\nu_1,\nu_{\omega\cdot2}+C_H(2)\}\)
\end{center}
If \(\alpha_{\omega\cdot3}+\nu_3< \alpha_{\omega\cdot3+1}\) then
\begin{center}
\(a(q,\vec{\alpha})=(\underaccent{\sim}{A})_H\cap\alpha_{mc}=C_H\restriction_{even}\cup\{C_H(\omega),C_H(\omega)+\nu_1,\nu_{\omega\cdot2}+C_H(2),\alpha_{\omega\cdot3}+\nu_3\}\)
\end{center}
Anyway, we have that \(a(q,\vec{\alpha})\in V[C^*]\) and therefore \(\langle q,p^{\frown}\vec{\alpha}\rangle|| \underaccent{\sim}{A}\cap\alpha_{mc}\) for every extension \(\vec{\alpha}\) of type X. Namely, \(q,X\) satisfy (*).
\\
\\
\textit{Proof of \ref{Finding aprox.}}: Let \(p=\langle t_1,...,t_n,t_{n+1}\rangle\). For every
\begin{center}
     \(X=\langle X_1,...,X_{n+1}\rangle \)- extension of \(p\) \ \ ,
     \(q \in\mathbb{Q}\)\ \ ,
   \(\vec{\alpha}\in B(p,X\setminus \langle x_{mc}\rangle)\) 
\end{center}
Recall that \(l_X=\min(i\mid X_i\neq\emptyset)\) and define \(B^X_{(0)}(q,\vec{\alpha})\) to be the set
\begin{center}
\(\{\theta \in B(t_{l_X},x_{mc}) \mid \exists a \exists (C(x_{i,j}))_{x_{i,j}} \ \ 
\langle q,p^{\frown}\langle\vec{\alpha},\theta,C(x_{i,j})\rangle\Vdash \underaccent{\sim}{A}\cap \theta=a\}\)\\

\end{center}
Also let \(B^X_{(1)}(q,\vec{\alpha})=B(t_{l_X},x_{mc})\setminus B^X_{(0)}(q,\vec{\alpha})\). One and only one of  \(B^X_{(1)}(q,\vec{\alpha}),B^X_{(0)}(q,\vec{\alpha})\) is in \(U(t_{l_X},x_{mc})\).
Define \(B^X(q,\vec{\alpha})\) and \(F^X_q(\vec{\alpha})\in\{0,1\}\) such that
\begin{center}
\(B^X(q,\vec{\alpha})=B^X_{(F^X_q(\vec{\alpha}))}(q,\vec{\alpha})\in U(t_{l_X},x_{mc})\)
\end{center}
Since \(|\mathbb{Q}|\leq 2^{\kappa^*}<\kappa(t_{l_X})\) we have
\(B^X(\vec{\alpha})=\underset{q}{\bigcap}B^X(q,\vec{\alpha})\in U(t_{l_X},x_{mc})\).
Define
\begin{center}

\(B^X(t_{l_X},x_{mc})=\underset{\vec{\alpha}}{\Delta}B^X(\vec{\alpha})\in U(t_{l_X},x_{mc})\)
\end{center}
Use lemma \ref{stabprod} to find \( B^X(t_i,x_{i,j})\subseteq B(t_i,x_{i,j})\), \(B^X(t_i,x_{i,j})\in U(t_i,x_{i,j})\) homogeneous for every \(F^X_q\). As before, if \(\lambda \notin X_i\) set  \(B^X(t_i,\lambda)=B(t_i,\lambda)\). Let 
$$p^*=p^{\frown}\langle (B^*(t_i))^{n+1}_{i=1}\rangle, \ B^{*}(t_i,\lambda)=\underset{X}{\bigcap}B^X(t_i,\lambda)$$
So far what we established the following property:
 if \(q,\vec{\alpha},(C(x_{i,j}))_{i,j},a\) are such that $$\langle q,p^{*\frown}\langle\vec{\alpha},(C(x_{i,j}))_{i,j}\rangle\rangle\Vdash \underaccent{\sim}{A}\cap\alpha_{mc}=a$$
since \(\alpha_{mc}\in B^X(q,\vec{\alpha}\setminus\langle\alpha_{mc}\rangle)\) we conclude that \(F^X_q(\vec{\alpha}\setminus\langle\alpha_{mc}\rangle)=0\). Let \(\vec{\alpha}'\) be another extension of type X, then \(\vec{\alpha}'\setminus\langle\alpha'_{mc}\rangle\) and \(\vec{\alpha}\setminus\langle\alpha_{mc}\rangle\)  belong to the same homogeneous set, thus 
$$F^X_q(\vec{\alpha}'\setminus\langle\alpha'_{mc}\rangle)=F^X_q(\vec{\alpha}\setminus\langle\alpha_{mc}\rangle)=0$$
By the definition of \(F^X_q(\vec{\alpha}'\setminus\langle\alpha'_{mc}\rangle)\) it follows that \(\alpha'_{mc}\in B^X_{(0)}(q,\vec{\alpha}'\setminus\langle\alpha'_{mc}\rangle)\) as wanted. For every \(\vec{\alpha}\in B(p',X)\) and \(q\in \mathbb{Q}\) fix some \((C_{i,j}(q,\vec{\alpha}))_{^{i\leq n+1}_{j\leq l_i+1}}\) such that 
 $$\langle q,p^{*\smallfrown}\langle \vec{\alpha},(C_{i,j}(q,\vec{\alpha}))_{^{i\leq n+1}_{j\leq l_i+1}}\rangle\rangle || \underaccent{\sim}{A}\cap\alpha_{mc}$$
 Let us argue that we can extend \(p^*\) to \(p^{**}\) such that for all \(1\leq i\leq n+1\) , \(1\leq j\leq l_i+1\) and \(\vec{\alpha}\in B(p^*,X)\), 
$$B(t^{**}_i)\cap(\alpha_s,\alpha_{i,j})\subseteq C_{i,j}(\vec{\alpha})$$ 
 Where \(\alpha_s\) is the predecessor of \(\alpha_{i,j}\) in \(\vec{\alpha}\) .
 In order to do that, fix \( i, \ j\) and stabilize \(C_{i,j}(\vec{\alpha})\) as follows:\\
 Fix \(\vec{\beta}\in B(p^*,\langle x_{1,1},...,x_{i,j}\rangle)\) By lemma \ref{stabprod} , the function 
 $$C_{i,j}(q,\vec{\beta},*):B(p^*,X\setminus\langle x_{1,1},...,x_{i,j}\rangle)\rightarrow P(\beta_{i,j})$$
  has homogeneous sets \(B'(\vec{\beta},x_{r,s},q)\subseteq B(t^*_r,x_{r,s})\) for \(x_{r,s}\in X\setminus\langle x_{1,1},...,x_{i,j}\rangle\). Denote the constant value by \(C^*_{i,j}(q,\vec{\beta})\). Define 
 $$B'(t^*_r,x_{r,s})=\underset{q\in\mathbb{Q}}{\underset{\vec{\beta}\in B(p^*,\langle x_{1,1},...,x_{i,j}\rangle)}{\Delta}}B'(\vec{\beta},x_{r,s},q), \ \ \text{ for }x_{r,s}\in X\setminus\langle x_{1,1},...,x_{i,j}\rangle$$
 Next, fix \(\alpha\in B(t^*_i,x_{i,j})\) and let
$$C^*_{i,j}(\alpha)=\underset{q\in\mathbb{Q}}{\underset{\vec{\alpha'}\in B(p^*,\langle x_{1,1},...,x_{i,j-1}\rangle)}{\Delta}}C^*_{i,j}(q,\vec{\alpha'},\alpha)$$
Thus \(C^*_{i,j}(\alpha)\subseteq\alpha\). Since \(\kappa(t_i)\) is ineffable, there is \(B'(t^*_i,x_{i,j})\subseteq B(t^*_i,x_{i,j})\) and \(C^*_{i,j}\) such that  for every \( \alpha\in B'(t^*_i,x_{i,j})\), \( C^*_{i,j}\cap\alpha=C^*_{i,j}(\alpha)\).
 By coherency, \(C^*_{i,j}\in \bigcap U(t_i,\xi)\).
 Finally, define  \(p^{**}=\langle t^{**}_1,...,t^{**}_n,t^{**}_{n+1} \rangle\), where
 $$B(t^{**}_i)=B'(t^*_i)\cap(\underset{j}{\cap}C^*_{i,j}) \ \ \ \ \ 1\leq i\leq n+1$$
 To see that \(p^{**}\) is as wanted, let \(\vec{\alpha}\in B(p^{**},X)\) and fix any \(i,j\). Then \(\vec{\alpha}\in B(p^{**},X)\) and \(\alpha_{i,j}\in B(t^{**}_i,x_{i,j})\), hence for any \(i,j\)
 \begin{center}
 \(B(t^{**}_i)\cap(\alpha_s,\alpha_{i,j})\subseteq C^*_{i,j}\cap\alpha_{i,j}\setminus\alpha_s=C^*_{i,j}(\alpha_{i,j})\setminus\alpha_s\subseteq C^*_{i,j}(\alpha_{1,1},...,\alpha_{i,j})=C_{i,j}(\alpha)\)
 \end{center}
  \begin{flushright}
 \(\blacksquare\)
 \end{flushright}
\begin{lemma}\label{ineffability}
\textit{Let \(p^*\) be as in lemma \ref{Finding aprox.} There exist \(p^{*} \leq p^{**}\) such that for every extension type X of \(p^{**}\) and \(q\in \mathbb{Q}\) that satisfies (*) there exists sets \(A(q,\vec{\alpha})\subseteq\kappa\) , for \( \vec{\alpha}\in B(p^{**},X\setminus\langle x_{mc}\rangle) \), such that for all \(\alpha\in B(p^{**},x_{mc})\)}
$$ A(q,\vec{\alpha})\cap\alpha=a(q,\vec{\alpha},\alpha)$$
\end{lemma}
\textbf{Example:} Recall that we have obtained the sets
\begin{center}
\(a(q,\vec{\alpha})=C_H\restriction_{even}\cup\{C_H(\omega),C_H(\omega)+\nu_1,\nu_{\omega\cdot2}+C_H(2)\}\cup b(q,\vec{\alpha})\)\end{center}
\begin{center}
\(b(q,\vec{\alpha})=\begin{cases}
 \emptyset & \alpha_{\omega\cdot3}+\nu_3\geq\alpha_{mc}\\
\{\alpha_{\omega\cdot3}+\nu_3\} & \alpha_{\omega\cdot3}+\nu_3<\alpha_{mc} \end{cases}\)

\end{center}
The element \(\alpha_{mc}\) is chosen from the set \(B(t_3,x_{mc})=B(t_3,0)\), by shrinking this set, we can directly extend \(p\) to \(p^*\) such that for every \(\vec{\alpha}\in B(p^*,X)\) , \(\alpha_{\omega\cdot3}+\nu_3<\alpha_{mc}\). Therefore, \begin{center}
\(A(q,\vec{\alpha})= C_H\restriction_{even}\cup\{C_H(\omega),C_H(\omega)+\nu_1,\nu_{\omega\cdot2}+C_H(2),\alpha_{\omega\cdot3}+\nu_3\}\)

\end{center}
\textit{Proof of \ref{ineffability}}: 
Fix \(q,X\) satisfying (*) and \(\vec{\alpha} \in B(p^{*},X\setminus\langle x_{mc}\rangle)\), since \(\kappa(t_i)\) is ineffable we can shrink the set \(B(t^*_{l_X},x_{mc})\) to \(B'(q,\vec{\alpha})\) to find sets \(A(q)\subseteq t_i\) such that 
\begin{center}
\(\forall \alpha\in B'(q,\vec{\alpha}) \ \ A(q,\vec{\alpha})\cap\alpha=a(q,\vec{\alpha},\alpha)\)
\end{center}define \(B_q(t^*_i,x_{mc})=\underset{\vec{\alpha}\in B(p^{*},X\setminus\langle x_{mc}\rangle)}{\Delta}B^{**}(q,\vec{\alpha})\). Intersect over all \(X,q\) and find \(p^*\leq p^{**}\) as before.
 \begin{flushright}
 \(\blacksquare\)
 \end{flushright}
 Thus there exists \(p_* \in G_{>\kappa^*}\) with the properties described in Lemma's \ref{Finding aprox.}-\ref{ineffability}. Next we would like to claim that for some sufficiently large family of \(q\in \mathbb{Q}\) and extension-type \(X\) we have \(q,X\) satisfy (*). 
\begin{lemma}\label{antichain for inaff}

\textit{ Let \(p_*\in G_{>\kappa^*}\) be as above and let X be any extension-type of \(p_*\). Then there exists a maximal antichain \(Z_{X}\subseteq \mathbb{Q}\) and extension-types \(X\preceq X_q\) for \(q\in Z_X\), unveiling the same maximal coordinate as \(X\)  such that for every \(q\in Z_{X}\), \ \(q,X_q\) satisfy (*).}
\end{lemma}
\textbf{Example:} The anti chain \(Z_X\) can be chosen as follows:
For any possible \(\nu_1,\nu_3\) choose a condition \(\langle\nu_1,\nu_3,\langle\kappa^*,B^*\rangle\rangle\in\mathbb{Q}\). This set definitely form a maximal anti chain, and by the same method of the previous examples taking \(X_q=X\) works. In general, if the maximal coordinate of X is some \(\omega\cdot(2n+1)\), \(Z_X\) will be the anti chain consisting of representative conditions for the \(2n+1\) first coordinates.\\
\\
\pr The existence of \(Z_{X}\) will follow from Zorn's Lemma and the method proving existence of \(X_q\) for some \(q\). Fix any \(\vec{\alpha}\in B(p_*,X)\), there exists a generic \(H\subseteq \mathbb{Q}\times\mathbb{M}_{>\kappa^*}[\vec{U}]\) with \(\langle 1_{\mathbb{Q}} , p^{\frown}_* \vec{\alpha}\rangle\in H=H_{\leq\kappa^*}\times H_{>\kappa^*}\). Consider the decomposition of \(\Mfor_{>\kappa^*}\) above \(p^{\frown}_* \vec{\alpha}\) induced by \(\alpha_{mc}\) and let \(p^{\frown}_*\vec{\alpha}=\langle p_1,p_2\rangle\),  i.e. \(\langle p_1,p_2\rangle\in (\Mfor_{>\kappa^*})_{\leq \alpha_{mc}}\times(\Mfor_{>\kappa^*})_{>\alpha_{mc}} \). \(H\) is generic for the forcing \(\mathbb{Q}\times(\Mfor_{>\kappa^*})_{\leq \alpha_{mc}}\times(\Mfor_{>\kappa^*})_{>\alpha_{mc}}\). Define  \( H_1=H_{\leq\kappa^*}\times(H_{>\kappa^*})_{\leq\alpha_{mc}}\) and \(H_2=H_{>\alpha_{mc}}\). Denote by \((\underaccent{\sim}{A})_{H_1} \in V[H_1]\) to be the name obtained by filtering only the part of $H_1$. It is a name of \(A\) in the forcing \(\Mfor_{>\alpha_{mc}}\). Above \(p_2\) we have sufficient closure to determine \( (\underaccent{\sim}{A})_{H_1}\cap\alpha_{mc}\)  
\begin{center}
\(\exists p^*_2\geq^*p_2 \ s.t.\ p^{*}_2\Vdash_{\Mfor_{>\alpha_{mc}}}(\underaccent{\sim}{A})_{H_1}\cap\alpha_{mc}=a\)
\end{center}
for some \(a\in V[C^*]\). Hence there exists \(\langle1_{\mathbb{Q}_{\leq\kappa^*}} , p_1\rangle\leq\langle q,p^{*}_1\rangle\) such that 
\begin{center}
\(\langle q,p^{*}_1\rangle\Vdash_{\mathbb{Q}\times\mathbb{M}_{\leq\alpha_{mc}}[\vec{U}]}\overset{\vee}{p^{**}_2}\Vdash_{\Mfor_{>\alpha_{mc}}}\underaccent{\sim}{A}\cap\alpha_{mc}=a\)
\end{center}
It is clear that \(\langle q,p^{*}_1,p^{*}_2\rangle ||_{\mathbb{Q}\times\mathbb{M}_{>\kappa^*}[\vec{U}]}\ \underaccent{\sim}{A}\cap\alpha_{mc}\). Finally, \(X_q\) is simply the extension type of \(p^{*}_1\). Since \(p^{*}_1\in\mathbb{M}_{\leq\alpha_{mc}}[\vec{U}]\), \(X_q\) unveils the same maximal coordinate as \(X\). By lemma \ref{Finding aprox.}, \( X_q,q\) satisfies \((*)\). 
 \begin{flushright}
 \(\blacksquare\)
 \end{flushright}
\begin{lemma}\label{change cofinality for stab}

\textit{ \( \kappa\) changes cofinality in \(V[A]\).}
\end{lemma}
\pr Let \(p_*=\langle t^*_1,...,t_n^*,t_{n+1}^*\rangle\in G_{>\kappa^*}\) be as before, \(\lambda_0=\otp(C_G)\) and \(\langle C_G(\xi) \mid\xi < \lambda_0\rangle\) be the Magidor sequence corresponding to \(G\). Work in V[A], define a sequence \(\langle\nu_i \mid \gamma(t_n^*,p_*)\leq i<\lambda_0\rangle \subset \kappa\):
\begin{center}

\(\nu_{\gamma(t_n^*,p_*)}=C_G(\gamma(t_n^*,p_*))+1=\kappa(t^*_n)+1\)

\end{center}
Assume that \(\langle\nu_{\xi'} \mid \xi'<\xi<\lambda_0\rangle\) is defined such that it is increasing and \(\nu_{\xi'}<\kappa\).
If \(\xi\) is limit define 
\begin{center}
\(\nu_{\xi}=\sup(\nu_{\xi'})+1\).
\end{center}
If \(\sup(\nu_{\xi'})=\kappa\) we are done, since \( \kappa \) changes cofinality to \(cf(\xi)<\lambda_0\) (which cannot hold for regular \(\lambda_0\)). Therefore, \(\nu_{\xi}<\kappa\).
If \(\xi=\xi'+1\), by proposition \ref{coordinates in the generic}, there exist an extension type \(X_\xi\) of \(p_*\) unveiling \(\xi\) as maximal coordinate. By lemma \ref{antichain for inaff} we can find \(Z_\xi\) and \(X_\xi\preceq X_q\) unveiling \(\xi\) as maximal coordinate such that \(q,X_q\)  satisfies (*). By lemma \ref{ineffability} there exists 
\begin{center}
\(A(q,\vec{\alpha})\)'s for \(q \in Z_\xi \ \ \vec{\alpha}\in B(p^*,X_q\setminus\langle x_{mc}\rangle)\).
\end{center}
Since \(A\notin V[C^*]\), \(A\neq A(q,\vec{\alpha})\). Thus define \(\eta(q,\vec{\alpha})=\min(A(q,\vec{\alpha})\Delta A)+1\)
\begin{center}

 \(\beta_{\xi}=\sup(\eta(q,\vec{\alpha}) \mid \vec{\alpha} \in [\nu_{\xi'}]^{<\omega}\cap B(p^*,X_q\setminus\langle x_{mc}\rangle), \  q\in 
 Z_\xi)\)
 \end{center}
It follows that \(\beta_{\xi}\leq\kappa\). Assume \(\beta_{\xi}=\kappa\), then \(\kappa\) changes cofinality but it might be to some other cardinal larger than \(\delta_0\), this is not enough in order to apply \ref{sufficient condition}\footnote{Actually, after proving Theorem \ref{Main theorem}, we can conclude that this phenomena cannot hold.}. Toward a contradiction, fix an unbounded and increasing sequence \(\langle\eta(q_i,\vec{\alpha_i}) \mid i<\theta<\kappa\rangle\) for some $q_i\in Z_\xi$ and $\vec{\alpha}_i\in [\nu_{\xi'}]^{<\omega}$. Notice that since \(\eta(q_i,\vec{\alpha}_i)<\eta(q_{i+1},\vec{\alpha}_{i+1})\) it must be that \(A(q_i,\vec{\alpha}_i)\neq A(q_{i+1},\vec{\alpha}_{i+1})\) and
\begin{center}
\(A(q_i,\vec{\alpha}_i)\cap\eta(q_i,\vec{\alpha}_i)=A\cap\eta(q_i,\vec{\alpha}_i)=A(q_{i+1},\vec{\alpha}_{i+1})\cap\eta(q_i,\vec{\alpha}_i)\) 

\end{center}
Define \(\eta_i=\min(A(q_i,\vec{\alpha_i})\Delta A(q_{i+1},\vec{\alpha}_{i+1}))\geq \eta(q_{i},\vec{\alpha}_{i})\). It follows that \(\langle \eta_i \mid i<\theta\rangle\) is a short cofinal sequence in \(\kappa\). This definition is independent of $A$ an only involve \(\langle\langle q_i,\vec{\alpha_i} \rangle\mid i<\theta<\kappa\rangle\), which can be coded as a bounded sequence of \(\kappa\). By the induction hypothesis there is \(C''\subseteq C\), bounded in \(\kappa\) such that $$V[C'']=V[\langle\langle q_i,\vec{\alpha}_i \rangle\mid i<\theta<\kappa\rangle]$$ Define \(C'=C^*\cup C''\), the model \(V[C']\) should keep \(\kappa\) measurable, since $C'$ is bounded, but also include the sequence \(\langle \eta_i \mid i<\theta\rangle\), contradiction.

Therefore, \(\beta_\xi<\kappa\), set \(\nu_\xi=\beta_\xi+1\). This concludes the construction of the sequence \(\nu_\xi\) . To see that the sequence is unbounded in \(\kappa\), let us show that \(C_G(\xi)<\nu_{\xi}\):

Clearly \(C_G(\gamma(t^*_n,p_*))<\nu_{\gamma(t^*_n,p_*)}\). Inductively sssume that \(C_G(i)<\nu_i, \ \gamma(t^*_n,p_*)\leq i<\xi\). If \(\xi\) is limit, since Magidor generic sequences are closed,
\begin{center}
\(C_G(\xi)=\sup(C_G(i) \mid i<\xi)\leq \sup(\nu_i \mid \gamma(t^*_n,p_*)\leq i<\xi)<\nu_{\xi}\)

\end{center}
If \(\xi=\xi'+1\) is successor, let \(\{q_\xi\}= Z_\xi\cap G_{\leq\kappa^*}\) 
\begin{center}

\(p_\xi=p^{\frown}_*\langle C_G(i_1),...,C_G(i_n),C_G(\xi)\rangle\in p^{\frown}_*X_\xi\cap G_{>\kappa^*}\)
\end{center}
By induction \(C_G(i_r)<\nu_{\xi'} \), therefore,  \(\eta(q_\xi,\langle C_G(i_1),...,C_G(i_n)\rangle)<\nu_{\xi}\) . Finally, \(\langle q_\xi, p_\xi\rangle\in G\),
\(\langle q_\xi, p_\xi\rangle\Vdash \underaccent{\sim}{A}\cap C_G(\xi)=A(q_\xi,\langle C_G(i_1),...,C_G(i_n)\rangle)\cap C_G(\xi)\), thus
$$A\cap C_G(\xi)=A(q_\xi,\langle C_G(i_1),...,C_G(i_n)\rangle)\cap C_G(\xi), \text{ hence } C_G(\xi)\leq\eta(q_\xi,\langle C_G(i_1),...,C_G(i_n)\rangle)<\nu_{\xi}.$$\begin{flushright}
 \(\blacksquare\)
 \end{flushright}

\newpage
\section{The Main Result Above \(\kappa\)}

In order to push the induction to sets above \(\kappa\) we will need a projection of \(\Mfor\) onto some forcing that adds a subsequence of \(C_G\). The majority of this chapter is the definition of this projection and some of its properties. The inductive argument will continue at lemma \ref{Thm for subsets of kappaplus}. 

Let G be generic and \(C_G\) the corresponding Magidor sequence. Let \(C^*\subseteq C_G\) be a subsequence and \(I=\Inx(C^*,C_G)\). Then \(I\) is a subset of \(\otp(C_G):=\lambda_0\), hence \(I\in V\).
Assume that \(\kappa^*=\sup(C^*)\) is a limit point in \(C_G\) and that \(C^*\) is closed i.e. containing all of its limit points below \(\kappa^*\)\footnote{ Note that for any $C\subseteq C_G$, we can take $Cl(C)=\{\alpha<\sup(C)\mid \alpha=\sup(C\cap\alpha)\}\cup C$, and $V[Cl(C)]=V[C]$, since $\Inx(C,Cl(C))\in V$.}. As we will see in the next lemma, one can find a forcing of the form \(\mathbb{M}_{\langle\nu_1,...,\nu_m\rangle}[\vec{U}]\), such that \(G\subseteq \mathbb{M}_{\langle\nu_1,...,\nu_m\rangle}[\vec{U}]\) is $V$-generic, which will be easier to project.
\begin{proposition}\label{making C a club}
 \textit{ Let \(G\) be \(\Mforc\)-generic and \(C^*\subseteq C_G\) such that \(C^*\) is closed and \(\kappa^*=\sup(C^*)\) is a limit point of \(C_G\). Then there exists \(\langle\nu_1,...,\nu_m\rangle\) such that G is generic for \(\mathbb{M}_{\langle\nu_1,...,\nu_m\rangle}[\vec{U}]\) and  for all \(1\leq i\leq m\), 
\(C^*\cap(\nu_{i-1},\nu_i)\) is either empty or a club in \(\nu_i\). (as usual denote \(\nu_0=0\))} 
\end{proposition}
\textbf{Example:} Assume that \(\lambda_0=\omega_1+\omega^2\cdot 2+\omega\), \(C^*\) is
\begin{center}
 \(C_G\restriction(\omega_1+1)\cup\{C_G(\omega_1+\omega+2),C_G(\omega_1+\omega+3)\}\cup\{C_G(\omega_1+\alpha)\mid \omega^2\cdot2 <\alpha <\lambda_0\}\)
\end{center}
Let \(\kappa_1<\kappa_2<\kappa_3<\kappa_4=\kappa\) be such that \(o^{\vec{U}}(\kappa_1)=\omega_1\)
, \(o^{\vec{U}}(\kappa_2)=o^{\vec{U}}(\kappa_3)=2\) and \(o^{\vec{U}}(\kappa)=1\).
We have
\begin{enumerate}
\item\((0,\kappa_1)\cap C^*=C_G\restriction\omega_1\).
\item\((\kappa_1,\kappa_2)\cap C^*=\{C_G(\omega_1+\omega+2),C_G(\omega_1+\omega+3)\}\).
\item\((\kappa_2,\kappa_3)\cap C^*=\emptyset\).
\item\((\kappa_3,\kappa_4)\cap C^*=\{C_G(\omega_1+\alpha)\mid \omega^2\cdot2 <\alpha <\lambda_0\}\).
\end{enumerate}
Then (1),(3),(4) are either empty or a club, but (2) is not. To fix this, we simply add \(\{C_G(\omega_1+\omega+2),C_G(\omega_1+\omega+3)\}\) to \(\kappa_1<\kappa_2<\kappa_3<\kappa_4\).
\vskip 0.2 cm
\textit{Proof of \ref{making C a club}}: By induction on \(m\), let us define a sequence
\begin{center}
 \(\vec{\nu}_m=\langle\nu_{1,m},...,\nu_{n_m,m}\rangle\)
 
\end{center}
such that for every \(m\),  \(G\) is generic for \(\mathbb{M}_{\vec{\nu}_m}[\vec{U}]\).
Define
 \(\vec{\nu}_0=\langle\kappa_1,...,\kappa_n\rangle\).
Assume that \(\vec{\nu}_{m}\) is defined with \(G\) generic, if for every \(1\leq i\leq n_m+1\) we have \(C^*\cap(\nu_{i-1,m},\nu_{i,m})\) is either empty or  unbounded (and therefore a club), stabilize the sequence at \(m\). Otherwise, let \(i\) be maximal such that \(C^*\cap(\nu_{i-1,m},\nu_{i,m})\) is nonempty and bounded. Thus,
\begin{center}
 \(\nu_{i-1,m}<\sup(C^*\cap(\nu_{i-1,m},\nu_{i,m}))<\nu_{i,m}\)
 
\end{center}
Since \(C^*\) is closed,   \(C_G(\gamma)=\sup(C^*\cap(\nu_{i-1,m},\nu_{i,m}))\in C^*\) for some \(\gamma\). As in lemma \ref{Properties of CG} we can find 
\begin{center}
 
\(\vec{\nu}_{m+1}=\langle \nu_{1,m},...,\nu_{i,m},\xi_{1},...,\xi_k,\nu_{i+1,m},...,\nu_{n_m,m}\rangle\subseteq C_G\)
\end{center}
such that \(C_G(\gamma)=\xi_k\) is unveiled and the forcing \(\mathbb{M}_{\vec{\nu}_{m+1}}[\vec{U}]\subseteq\mathbb{M}_{\vec{\nu}_{m}}[\vec{U}]\) is a subforcing of \(\mathbb{M}_{\vec{\nu}_{m}}[\vec{U}]\) with \(G\) one of its generic sets. Note that the maximal ordinal in the sequence \(\vec{\nu}_{m+1}\) such that \(C^*\cap(\nu_{j-1,m+1},\nu_{j,m+1})\) is nonempty and bounded is strictly less than \( \nu_{i,m}\). Therefore this iteration stabilizes at some \(N<\omega\). Consider the forcing \(\mathbb{M}_{\vec{\nu}_N}[\vec{U}]\), by the construction of the \(\vec{\nu}_r\)'s, for every \(1\leq i\leq n_N+1\)  \(C^*\cap(\nu_{i-1,N},\nu_{i,N})\) is either empty or unbounded (Since \(\vec{\nu}_{N+1}=\vec{\nu}_N\)).
\begin{flushright}
 \(\blacksquare\)
 \end{flushright}
\vskip 0.2 cm

Assume that \(\Mforc\) and \(C^*\) satisfy the property of \ref{making C a club}. 
Let us define a projection of 
$$\Mforc=\prod^n_{i=1}(\mathbb{M}_{\kappa_i})_{>\kappa_{i-1}}$$
 onto some forcing \(\prod^n_{i=1}\mathbb{P}_i\). We can define such a projection, by projecting each factor
$$\pi_i:(\mathbb{M}_{\kappa_i})_{>\kappa_{i-1}}\rightarrow\mathbb{P}_i \ \ \ (1\leq i\leq n)$$
and derive \(\pi:\Mforc\rightarrow\prod^n_{i=1}\mathbb{P}_i\).  First, if \(C^*\cap(\kappa_{i-1},\kappa_i)\) is empty, the projection is going to be to the trivial forcing. Otherwise, \(C^*\cap(\kappa_{i-1},\kappa_i)\) is a club at $\kappa_i$. In order to simplify notation, we will assume that \((\mathbb{M}_{\kappa_i})_{>\kappa_{i-1}}=\Mfor_{\langle\kappa\rangle}=\Mfor\) and \(C^*=C^*\cap(\kappa_{i-1},\kappa_i)\) is a club in \(\kappa\).  It seems natural that the projection will keep only the coordinates in \(I=\Inx(C^*,C_G)\), namely:
\begin{definition}\label{definition of projection}
Let \(p=\langle t_1,...,t_{n+1}\rangle\in\Mfor\), define \textit{the projection to the $I$ coordinates} by, $$\pi_I(p)=\langle t'_i \mid \gamma(t_i,p)\in I\rangle^{\frown}\langle t_{n+1}\rangle\text{, \  where  }\   
t'_i= \begin{cases} 
 \kappa(t_i) & \gamma(t_i,p)\in \succ(I) \\
t_i & \gamma(t_i,p) \in \Lim(I) 
\end{cases}$$
\end{definition}
Let us define a forcing notion \(\mathbb{P}_i=\MforI\) (the range of the projection \(\pi_I\)) that will add the subsequence \(C^*\), such that the forcing \(\Mfor\) (more precisely, a dense subset of \(\Mfor\)) projects onto \(\MforI\) via the projection \(\pi_I\) we have just defined in \ref{definition of projection}.
\subsection{The forcing \textbf{\(\MforI\)}}

Considering \(C^*\) as a function with domain \(I\), we would like to have a function similar to \(\gamma(t_i,p)\) that tells us which coordinate is currently unveiled. Given \(p=\langle t_1,...,t_n,t_{n+1}\rangle\), define recursively:
\begin{enumerate}
    \item \(I(t_0,p):=0\).
    \item If \(\{j \in I\setminus I(t_{i-1},p)+1 \mid o_L(j)=o^{\vec{U}}(t_i)\}=\emptyset\), then $I(t_i,p)=N/A$ undefined.
    \item If \(\{j \in I\setminus I(t_{i-1},p)+1 \mid o_L(j)=o^{\vec{U}}(t_i)\}\neq\emptyset\), define
    $$I(t_i,p):=\min(j \in I\setminus I(t_{i-1},p)+1 \mid o_L(j)=o^{\vec{U}}(t_i))$$
\end{enumerate}
If for every $0\leq i\leq n$, $I(t_i,p)\neq N/A$, we say that $I$ is defined on $p$.

\textbf{Example:} Consider Magidor forcing adding a sequence of length \(\omega^2\) i.e. $o^{\vec{U}}(\kappa)=2$ and \(C_G=\{C_G(\alpha)\mid \alpha<\omega^2\}\). Assume
 \(C^*=\{C_G(0)\}\cup\{C_G(\alpha)\mid \omega\leq\alpha<\omega^2\}\), hence, \(I=\{0\}\cup(\omega^2\setminus\omega)\). The \(\omega\)-th element of \(C_G\) is no longer limit in \(C^*\). Let
 \begin{center}
  \(p=\langle\underset{t_1}{\underbrace{\langle\kappa(t_1),B(t_1)\rangle}},\underset{t_2}{\underbrace{\langle\kappa,B(t_2)\rangle}}\rangle\)
 \end{center}
 Where \(o^{\vec{U}}(t_1)=1\). Computing \(I(t_1,p)\),
 \begin{center}
  \(I(t_1,p)=\omega=\gamma(t_1,p)\)
 \end{center}
Therefore \(\pi_I(p)=\langle\kappa(t_1),t_2\rangle\).\\
\\
\begin{definition}\label{Definition of condition of MIU}
 The \textit{Magidor forcing adding a sequence prescribed by I}, denoted by \(\MforI\), consist of conditions of the form \(p=\langle t_1,...,t_{n+1}\rangle\) such that:
\begin{enumerate}
\item $I$ is defined on $p$.
\item \(\kappa(t_1)<...<\kappa(t_n)<\kappa(t_{n+1})=\kappa\).
\item For \(i=1,...,n+1\),
\begin{enumerate}

    \item If \(I(t_i,p) \in \succ(I)\), \begin{enumerate}
    \item \(t_i=\kappa(t_i)\).
   \item \(I(t_{i-1},p)\) is the predecessor of \(I(t_i,p)\) in \(I\).
   \item If \(I(t_{i-1},p)+\overset{m}{\underset{i=1}{\sum}}\omega^{\gamma_i}=I(t_i,p)\)  (C.N.F) , then \footnote{Recall that \(Y(\gamma)=\{\alpha<\kappa\mid o^{\vec{U}}(\alpha)=\gamma\}\)},
   $$\Big(Y(\gamma_1)\times ... \times Y(\gamma_{m-1})\Big)\cap[(\kappa(t_{i-1}),\kappa(t_i))]^{<\omega}\neq\emptyset$$
     \end{enumerate}
    \item If \( I(t_i,p) \in \Lim(I) \), then,
    \begin{enumerate}
    \item \(t_i=\langle\kappa(t_i),B(t_i)\rangle\). \ ,  \(B(t_i)\in\underset{\xi<o^{\vec{U}}(t_i)}{\bigcap} U(t_i,\xi) \)
    \item \(I(t_{i-1},p)+\omega^{o^{\vec{U}}(t_i)}=I(t_i,p)\).
    \item \(\min(B(t_i))>\kappa(t_{i-1})\).
    \end{enumerate}
\end{enumerate}
\end{enumerate}
\end{definition}
\begin{definition}\label{Definition of order of MIU}
 Let \(p=\langle t_1,...,t_n,t_{n+1}\rangle,q=\langle s_1,...,s_m,s_{m+1}\rangle\in \MforI\). Define the \textit{order of $\MforI$}, \(\langle t_1,...,t_n,t_{n+1}\rangle\leq_I\langle s_1,...,s_m,s_{m+1}\rangle\) iff
 \(\exists 1\leq i_1<...<i_n\leq m<i_{n+1}=m+1\) such that:
\begin{enumerate}
\item \(\kappa(t_r)=\kappa(s_{i_r})\) and \(B(s_{i_r})\subseteq B(t_r)\).
\item If \(i_k<j<i_{k+1}\), then
\begin{enumerate}
\item \(\kappa(s_j) \in B(t_{k+1})\). 
\item If \(I(s_j,q)\in \succ(I)\), then, $$\Big(B(t_{k+1},\gamma_1)\times...\times B(t_{k+1},\gamma_{k-1})\Big)\cap[(\kappa(s_{j-1}),\kappa(s_j))]^{<\omega}\neq\emptyset$$
where 
\(I(s_{i-1},q)+\overset{k}{\underset{i=1}{\sum}}\omega^{\gamma_i}=I(s_i,q)\) (C.N.F).

\item If \(I(s_j,q)\in \Lim(I)\) then \( B(s_j)\subseteq B(t_{k+1})\cap \kappa(s_j)\).
\end{enumerate}
\end{enumerate}
\end{definition}
\begin{definition}\label{Definition of direct order of MIU}
 Let \(p=\langle t_1,...,t_n,t_{n+1}\rangle , q=\langle s_1,...,s_m,s_{m+1}\rangle \in \MforI\),  \(q\) is a \textit{direct extension of \(p\)}, denoted \(p\leq^*_Iq\), iff:
\begin{enumerate}
\item \(p\leq_Iq\).
\item \(n=m\).
\end{enumerate}
\end{definition}
  \textbf{Remarks:} \begin{enumerate}
  
  \item In definition \ref{Definition of condition of MIU} (b.i), although it seems superfluous to take all the measures corresponding to \(t_i\) as well as those which do not take an active part in the development of \(C^*\), the necessity is apparent when examining definition \ref{Definition of order of MIU} (2.b)- the \(\gamma_i\)'s may not be the measures taking active part in \(C^*\). In lemma \ref{pi is a projection} this condition will be crucial when completing \(C^*\) to \(C_G\). 
\item As we have seen in earlier chapters, the function \(\gamma(t_i,p)\) returns the same value when extending \(p\). \(I(t_i,p)\) have the same property, let 
\(p=\langle t_1,...,t_n,t_{n+1}\rangle , \ q=\langle s_1,...,s_m,s_{m+1}\rangle\in \MforI \), such that \(p\leq_I q\), by \ref{Definition of condition of MIU} (2.b.ii),
\( I(t_r,p)=I(s_{i_r},q)\).
\item In definition \ref{Definition of direct order of MIU}, since \(n=m\) we only have to check (1) of definition \ref{Definition of order of MIU}.
\item Let \(p=\langle t_1,...,t_{n+1}\rangle\in \MforI\) be any condition. Assume we would like to unveil a new index  \(j\in I\) between \(I(t_i,p)\) and \(I(t_{i+1},p)\). It is possible if for example \(j\) is the successor of \(I(t_i,p)\) in \(I\):\\ Assume \(I(t_i,p)+\overset{m}{\underset{l=1}{\sum}}\omega^{\gamma_l}=j\) (C.N.F), then \(\gamma_l<o^{\vec{U}}(t_{i+1})\). Extend \(p\) by choosing \(\alpha \in B(t_{i+1},\gamma_m)\) above some sequence 
$$\langle\beta_1,...,\beta_k\rangle\in B(t_{i+1},\gamma_1)\times...\times B(t_{i+1},\gamma_{m-1})$$
Then
 $$I(\alpha,p^{\frown}\langle\alpha\rangle)=\min(r\in I\setminus I(t_i,p)\mid o_L(r)=o_L(j))=j$$
Another possible index is any \(j\in \Lim(I)\) such that \(I(t_{i},p)+\omega^{o_L(j)}=j\). For such \(j\), extend \(p\) by picking  \(\alpha \in B(t_{i+1},o_L(j))\) above some sequence \(\langle\beta_1,...,\beta_k\rangle\),
 to obtain
\begin{center}
\(p\leq_I \langle t_1,...,t_i ,\langle\alpha,\underset{\xi<o_L(j)}{\bigcap}B
(t_{i+1},\xi)\cap\alpha\rangle,\langle\kappa(t_{i+1}),B(t_{i+1})\setminus(\alpha+1)\rangle,...,t_{n+1}\rangle\)
\end{center} 
A routine verification of definition \ref{Definition of condition of MIU} asserts that in both cases the extension of \(p\) is in \(\MforI\).
\end{enumerate} 
 The forcing \(\MforI\) has lots of the properties of \(\Mfor\), however, they are irrelevant for the proof. Therefore, we will state only few of them.
 \begin{lemma}\label{MUI is kappaplus cc}
 \textit{\(\MforI\) satisfy \(\kappa^+-c.c\)}.
 
 \end{lemma}
 \pr Let \(\{\langle t_{\alpha,1},...,t_{\alpha,n_{\alpha}}\rangle=p_{\alpha}\mid\alpha<\kappa^+\}\subseteq \MforI\). Find \(n<\omega\) and
\(E\subseteq\kappa^+ ,\ |E|=\kappa^+\) and \(\langle\kappa_1,...,\kappa_n\rangle \) such that \(\forall \alpha \in E\),
\begin{center}

 \(n_\alpha=n,\text{ and } \langle\kappa(t_{\alpha,1}),...,\kappa(t_{\alpha,n_{\alpha}})\rangle=\langle\kappa_1,...,\kappa_n\rangle\)
\end{center}
Fix any \(\alpha,\beta \in E\). Define \(p^*=\langle t_1,...,t_n,t_{n+1}\rangle\) where
\begin{center}
\(B^*(t_i)=B(t_{i,\alpha})\cap B(t_{i,\beta})\in \underset{\xi<o^{\vec{U}}(\kappa_i)}{\bigcap} U(\kappa_i,\xi)\)
\end{center}
\begin{center}
\(t_i=
\begin{cases} 
\langle \kappa_i,B^*(t_i)\rangle & I(t_i,p)\in \Lim(I) \\
\kappa_i & otherwise
\end{cases}\)
\end{center}
Since \(p_{\alpha},p_{\beta} \in \MforI\), it is clear that \(p^* \in \MforI\) and also \(p_\alpha,p_\beta\leq^*_Ip^*\).
\begin{flushright}
 \(\blacksquare\)
 \end{flushright}
\begin{lemma}\label{Properties of MUI}
\textit{Let \(G_I\subseteq \MforI\) be generic , define }
\begin{center}

\(C_I=\bigcup\{\{\kappa(t_i) | i=1,...,n\} \mid \langle t_1,...,t_n,t_{n+1}\rangle\in G_I \} \)
\end{center}
\textit{Then}
\begin{enumerate}
\item\textit{ \(\otp(C_I)=\otp(I)\) (thus we may also think of \(C_I\) as a function with domain \(I\)).
\item \(G_I\) consist of all conditions \(p=\langle t_1,...,t_n,t_{n+1}\rangle\in \MforI\) such that} 
\begin{enumerate}
\item  \(C_I(I(t_i,p))=\kappa(t_i)\). 
\item \(C_I\cap(\kappa(t_{i-1}),\kappa(t_{i})) \subseteq B(t_i)\) \ \ \ \(1\leq i\leq n+1\).
\item \textit{\(\forall i\in \succ(I)\cap(I(t_r,p),I(t_{r+1},p))\) with predecessor \(j\in I\) such that \(j+\overset{k}{\underset{l=1}{\sum}}\omega^{\gamma_l}=i\) (C.N.F) we have }
\begin{center}
\([(C_I(j),C_I(i))]^{<\omega}\cap B(t_{r+1},\gamma_1)\times...\times B(t_{r+1},\gamma_{k-1})\neq\emptyset\)
\end{center}
\end{enumerate}    
\end{enumerate}

\end{lemma}
\pr For (1) , let us consider the system of ordered sets of ordinals \( (\kappa(p),i_{p,q})_{p,q}\) where
\begin{center}

\(\kappa(p)=\{\kappa(t_1),...,\kappa(t_n)\}\) for \( p=\langle t_1,...,t_{n+1}\rangle\in G_I\) 
\end{center}
 \(i_{p,q}:\kappa(p)\rightarrow\kappa(q)\) are defined for \(p=\langle t_1,...,t_{n+1}\rangle\leq_I \langle s_1,...,s_{m+1}\rangle=q\) as the inclusion:
\begin{center}
\(i_{p,q}(\kappa(t_r))=\kappa(t_r)=\kappa(s_{i_r})\)  (\(i_r\) are as in the definition of \(\leq_I\) ) 
\end{center}
Since \(G_I\) is a filter, \( (\kappa(p),i_{p,q})_{p,q}\) form a directed system with a direct ordered limit $$\underrightarrow{\Lim} \ \kappa(p)=\underset{p\in G_I}{\bigcup }\kappa(p)=C_I\text{ and inclusions } 
i_p:\kappa(p)\rightarrow C_I$$
We already defined for  \(p,q\in G_I\) such that \(p\leq_Iq\), commuting functions
$$I(*,p):\kappa(p)\rightarrow I, \ \ (*,p)=I(*,q)\circ i_{p,q}$$
Thus \((I(*,p))_{p\in G}\) form a compatible system of functions, and by the universal property of directed limits, we obtain 
$$I(*):C_I\rightarrow I,  \ I(*)\circ i_p= I(*,p)$$
Let us show that $I$ is an isomorphism of ordered set: Since  \(I(*,p)\) are injective \(I(*)\) is also injective. Assume \(\kappa_1<\kappa_2 \in C_I\), find \(p\in G_I\) such that \(\kappa_1,\kappa_2\in \kappa(p)\). Therefore, \(I(\kappa_i,p)=I(\kappa_i)\) preserve the order of \(\kappa_1,\kappa_2\).
Fix \(i\in I\), it suffices to show that there exists some condition \(p\in G_I\) such that \(i\in Im(I(*,p))\). To do this, let us show that the set of all conditions \(p\in\MforI\) with \(i\in Im(I(*,p))\) is a dense subset of \(\MforI\). Let \(p=\langle t_1,...,t_{n+1}\rangle\in \MforI\) be any condition , if \( i\in Im(I(*,p))\) then we are done. Otherwise, there exists \(0\leq k\leq n\) such that,
$$I(t_k,p)<i<I(t_{k+1},p)$$
therefore \(I(t_{k+1},p)\in \Lim(I)\). By induction on $i$, let us argue that it is possible to extend \(p\) to a condition \(p'\), such that \(i \in Im(I(*,p'))\). If 
$$\overset{k}{\underset{l=1}{\sum}}\omega^{\gamma_l}=i=\min(I)\ \text{ (C.N.F)}$$ 
it follows that
\(i<I(t_1,p)\). 
By definition \ref{Definition of condition of MIU} (2.b.ii), \(I(t_1,p)=\omega^{o^{\vec{U}}(t_1)}\). To extend \(p\) just pick any \(\alpha\) above some sequence
$$\langle\beta_1,...,\beta_k\rangle\in B(t_1,\gamma_1)\times...\times B(t_1,\gamma_{k-1})$$
and
$$p\leq_I\langle\alpha,\langle\kappa(t_1),B(t_1)\setminus(\alpha+1)\rangle,t_2,...,t_{n+1}\rangle\in \MforI$$
If \(i\in \succ(I)\) with predecessor \(j\in I\). By the induction hypothesis, we can assume that for some \(k\), \(j=I(t_k,p)\in Im(I(*,p))\). Thus by the remarks following definition \ref{Definition of direct order of MIU} we can extend \(p\) by some \(\alpha\) such that \(i\in Im(I(*,p))\).
Finally if \(i\in \Lim(I)\), then
$$i=\alpha+\omega^{o_L(i)},\text{ where }\alpha:=\overset{m}{\underset{i=1}{\sum}}\omega^{\gamma_i} \  (C.N.F)$$
therefore 
 \(\forall\beta\in(\alpha,i) ,\ \beta+\omega^{o_L(i)}=i\).
Take any \(i'\in I\cap(\alpha,i)\). Just as before, it can be assumed that
\(i'=I(t_k,p)\), thus \(I(t_k,p)+\omega^{o_L(i)}=i\). By the same remark, we can extend \(p\) to some \(p'\in\MforI\) with \(j\in Im(I(*,p'))\).\\
\\
For (2), let 
\(p=\langle t_1,...,t_{n+1}\rangle\in G_I\).
(a) is satisfied by the argument in (1).
Fix \(\alpha \in C_I\cap(\kappa(t_i),\kappa(t_{i+1}))\), there exists \(p\leq_I p'=\langle s_1,...,s_m\rangle\in G_I\) such that \(\alpha \in \kappa(p')\)  thus \(\alpha\in B(t_{i+1})\) by definition. Moreover, if \(I(\alpha,p')\in \succ(I)\) with predecessor \(j\in I\) , then by definition \ref{Definition of condition of MIU} (2.a.ii), there is \(s_k\) such that  \(j=I(s_k,p')\) and by definition \ref{Definition of order of MIU} (2.b) 
\begin{center}
 \([(\kappa(s_{k-1}),\kappa(s_k))]^{<\omega}\cap B(t_{i+1},\gamma_1)\times...\times B(t_{i+1},\gamma_{k-1})\neq\emptyset\)
\end{center}
From (a),
\begin{center}

\(\kappa(s_k)=C_I(j)\) and \(\kappa(s_{k+1})=C_I(i)\)
\end{center}
In the other direction, if \(p=\langle t_1,...,t_{n+1}\rangle\in \MforI\) satisfies (a)-(c). By (a), there exists some \(p''\in G_I\) with \(\kappa(p)\subseteq\kappa(p'')\). Set \(E\) to be
\begin{center}
 \(\{\langle w_1,...,w_{l+1}\rangle\in (\MforI)_{\geq_Ip''}\mid \kappa(w_j)\in B(t_i)\cup\{\kappa(t_i)\} \rightarrow \ B(w_j)\subseteq B(t_i)\}\)
\end{center}
\(E\) is dense in \(\MforI\) above \(p''\). Find \(p''\leq_Ip'=\langle s_1,...,s_{m+1}\rangle \in G_I\cap D\).  Checking definition \ref{Definition of order of MIU}, Let us show that \(p\leq_Ip'\):
For (1), since \(\kappa(p)\subseteq\kappa(p')\) there is a natural injection \(1\leq i_1<...<i_n\leq m\) which satisfy \(\kappa(t_r)=\kappa(s_{i_r})\). Since \(p'\in E\), \(B(s_{i_r})\subseteq B(t_r)\). (2a), follows from condition (b), (2b) follows from condition (c). Since \(p'\in E\), if \(i_r<j<i_{r+1}\) then \(\kappa(s_j)\in B(t_{r+1})\), thus, (2c) holds.
\begin{flushright}
 \(\blacksquare\)
 \end{flushright}

Given a generic set \(G_I\) for \(\MforI\) , we have \(V[C_I]=V[G_I]\). Once we will show that \(\pi_I\) is a projection, then for every \(G\subseteq\Mfor\) generic, 
\begin{center}
 \(\pi_{I*}(G):=\{p\in \MforI\mid \exists q\in\pi''_IG, \ p\leq_Iq\}\)
\end{center}
will be generic for \(\MforI\) and by the definition of \(\pi_I\) \ref{definition of projection}, we conclude that the corresponding sequence to \(\pi_I(G)\) is \(C^*\), this is stated formally on corollary \ref{Cor to generate subsequence}. Let us turn to the proof that \(\pi_I\) is a projection. 
\begin{definition}\label{Definition of D}
Let \(D\) be the set of all 
$$p=\langle t_1,...,t_n,t_{n+1}\rangle\in \Mfor \ ,\ \pi_I(p)=\langle t'_{i_1},...,t'_{i_m},t_{n+1}\rangle$$ 
such that:
\begin{enumerate}
\item If \(\gamma(t_{i_j},p)\in\Lim(I)\) then \(\gamma(t_{i_{j-1}},p)=\gamma(t_{i_j-1},p)\).
\item If \(\gamma(t_{i_j},p) \in \succ(I)\) then \(\gamma(t_{i_{j-1}},p)\) is the predecessor of \(\gamma(t_{i_j},p)\) in \(I\).
\end{enumerate}

\end{definition}
Condition \((1)\) is to be compared with definition \ref{Definition of condition of MIU} (2.b.ii) and condition (2) with (2.a.ii). The following example justifies the necessity of D.\\
\\
\textbf{Example:} Assume that
\begin{center}
 
\(\lambda_0=\omega^2\) and \(I=\{2n\mid n\leq \omega\}\cup\{\omega+2,\omega+3\}\cup\{\omega\cdot n\mid n<\omega\}\)
\end{center}
let \(p\) be the condition
\begin{center}
 \(\langle\underset{t_1}{\underbrace{\langle \nu_{\omega},B_\omega\rangle}},\underset{t_2}{\underbrace{ \nu_{\omega+1}}},\underset{t_3}{\underbrace{\langle \nu_{\omega\cdot2},B_{\omega\cdot2}\rangle}},\underset{t_4}{\underbrace{\langle \kappa,B\rangle}}\rangle\)
 \end{center}
 \begin{center}
 \(\pi_I(p)=\langle\underset{t_1\mapsto t'_{i_1}}{\underbrace{\langle \nu_{\omega},B_\omega\rangle}},\underset{t_3\mapsto t'
 _{i_2}}{\underbrace{\nu_{\omega\cdot2}}},\underset{t_4}{\underbrace{\langle \kappa,B\rangle\rangle}}\)
\end{center}
The \(\omega+2,\omega+3\)-th coordinates cannot be added. On one hand, they should be chosen below \(\nu_{\omega\cdot2}\), on the other hand, there is no large set associated to $\nu_{\omega\cdot2}$. In $D$, this situation is impossible due to condition $(2)$ of definition \ref{Definition of D}, which $p$ fails to satisfy:
$$\omega\cdot2\in\succ(I)\text{ but }
 \omega+3\in I\text{ is the predecessor and } \gamma(t_{i_2})=\omega\cdot2$$
 Notice that we can extend \(p\) to
\begin{center}
 \(\langle\langle \nu_{\omega},B_\omega\rangle, \nu_{\omega+1},\nu_{\omega+2},\nu_{\omega+3},\langle \nu_{\omega\cdot2},B_{\omega\cdot2}\rangle,\langle \kappa,B\rangle\rangle\)
\end{center}
to fix this problem. 

Next consider
\begin{center}
 \(I=\{2n\mid n\leq \omega\}\cup\{\omega+2,\omega+3\}\cup\{\omega\cdot n\mid n<\omega, \ n\neq 2\}\)
\end{center}
and let \(p\) be the condition
\begin{center}
 \(\langle\underset{t_1}{\underbrace{\langle \nu_{\omega},B_\omega\rangle}},\underset{t_2}{\underbrace{ \langle\nu_{\omega\cdot2},B_{\omega\cdot2}\rangle}},\underset{t_3}{\underbrace{\langle \nu_{\omega\cdot3},B_{\omega\cdot3}\rangle}},\underset{t_4}{\underbrace{\langle \kappa,B\rangle}}\rangle\)
  \end{center}
 \begin{center}
 \(\pi_I(p)=\langle\underset{t_1\mapsto t'_{i_1}}{\underbrace{\langle \nu_{\omega},B_\omega\rangle}},\underset{t_3\mapsto t'
 _{i_2}}{\underbrace{\langle \nu_{\omega\cdot3},B_{\omega\cdot3}\rangle}},\underset{t_4}{\underbrace{\langle \kappa,B\rangle\rangle}}\)

\end{center}
Once again the coordinates \(\omega+2,\omega+3\) cannot be added since \\ \(\min(B_{\omega\cdot3})>\nu_{\omega\cdot2}\). This problem points out condition (1) of definition \ref{Definition of D}, which $p$ fails to satisfy: $$ \gamma(t_{i_1},p)=\omega<\omega\cdot2=\gamma(t_{i_2-1},p)$$
As before, we can extend \(p\) to avoid this problem.
\begin{proposition}\label{D is dense}
 \textit{\(D\) is dense in \(\Mfor\)}.
 
\end{proposition}
\pr Fix \(p=\langle t_1,...,t_{n+1}\rangle\in \Mfor\), define  recursively \(\langle p_k \mid k<\omega\rangle\) as follows:

First, \(p_0=p\). Assume that \(p_k=\langle t^{(k)}_1,...,t^{(k)}_{n_k},t^{(k)}_{n_k+1}\rangle\) is defined. If \(p_k\in D\), define \(p_{k+1}=p_k\). Otherwise, there exists a maximal \(1\leq i_j=:i_j(k)\leq n'+1\) such that \(\gamma(t^{(k)}_{i_j},p_k)\in I \) which fails to satisfy  \((1)\) or fails to satisfy \((2)\) of definition \ref{Definition of D}. Let us split into two cases accordingly:
\begin{enumerate}
    \item Assume \(\neg(1)\), Thus
$$\gamma(t^{(k)}_{i_j},p_k)\in \Lim(I) \text{ and } \gamma(t^{(k)}_{i_{j-1}},p_k)<\gamma(t^{(k)}_{i_j-1},p_k)$$
Since \(\gamma(t^{(k)}_{i_j},p_k)\in \Lim(I)\) there exists \(\gamma\in I\cap (\gamma(t^{(k)}_{i_j-1},p_k),\gamma(t^{(k)}_{i_j},p_k))\).
Use proposition \ref{coordinates in the generic} to find \(p_{k+1}\geq p_k\) with \(\gamma\) added and the only other coordinates added are below \(\gamma\), thus if \( t^{(k)}_{i_j}=t^{(k+1)}_{r}\) then \(\gamma=\gamma(t^{(k+1)}_{{r-1}},p_{k+1})\). Thus, every \(l\geq r\) satisfies   \((1)\) and \(2)\). If \(p_{k+1}\notin D\) then the problem must accrue below \(\gamma(t^{(k)}_{i_j},p_k)\).
\item Assume \(\neg(2)\), thus 
$$\gamma(t^{(k)}_{i_j},p) \in \succ(I)\text{ and } \gamma(t^{(k)}_{{i_j-1}},p)\text{ is not the predecessor of }(\gamma(t^{(k)}_{i_j},p))$$
Let \(\gamma\) be the predecessor in $I$ of \(\gamma(t^{(k)}_{i_j},p)\). By proposition \ref{coordinates in the generic}, there exist \(p_{k+1}\geq p_k\) with \(\gamma\) added and the only other coordinates added are below \(\gamma\). As before, if \( t^{(k)}_{i_j}=t^{(k+1)}_{r}\) then \(\gamma=\gamma(t^{(k+1)}_{{r-1}},p_{k+1})\) and for every \(l\geq r\), \(\gamma(t^{(k+1)}_{l},p_{k+1})\)  satisfies \((1)\) and \((2)\).
\end{enumerate}
The sequence \(\langle p_k \mid k<\omega\rangle\) is defined. It necessarily stabilizes, otherwise the sequence \(\gamma(t^{(k)}_{i_j(k)},p_k)\) form a strictly decreasing infinite sequence of ordinals. Let \(p_{n^*}\) be the stabilized condition, it is an extension of \(p\) in \(D\).
\begin{flushright}
 \(\blacksquare\)
 \end{flushright}
\begin{lemma}\label{pi is a projection}
  \(\pi_I\restriction D:D\rightarrow \MforI\) is a projection, i.e:
\begin{enumerate}
    \item \(\pi_I\) is onto.
    \item \(p_1\leq p_2 \Rightarrow \pi_I(p_1)\leq_I \pi_I(p_2) \) (also \(\leq^*\) is preserved).
    \item \( \forall p\in \Mfor \  \forall q\in\MforI\ (  \pi_I(p)\leq_I q \rightarrow \exists p'\geq p \ \ (q=\pi_I(p'))\).
\end{enumerate}

\end{lemma}
\pr
Let \(p\in D\), such that \(\pi_I(p)=\langle t_{i_1}',...,t_{i_{n'}}',t_{n+1}\rangle\)\\
\\
\underline{Claim:} \(\pi_I(p)\) computes \(I\) correctly i.e. for every \(0\leq j\leq n'\), we have the equality \(\gamma(t_{i_j},p)=I(t_{i,j}',\pi_I(p))\).\\
\\
\textit{Proof of claim:} By induction on \(j\), for \(j=0\), \(\gamma(0,p)=0=I(0,\pi_I(p))\) . For \(j>0\), assume  \(\gamma(t_{i_{j-1}},p)=I(t_{i_{j-1}}',\pi_I(p))\) and \(\gamma(t_{i_j},p)\in\succ(I)\). Since \(p\in D\), \(\gamma(t_{i_{j-1}},p)\) is the predecessor of \(\gamma(t_{i_{j}},p)\) in \(I\).
Use the induction hypothesis to see that
\begin{center}
 \(I(t'_{i_j},\pi_I(p))=\min(\beta\in I\setminus \gamma(t_{i_{j-1}},p)+1\mid o_L(\beta)=o^{\vec{U}}(t_{i_j}))=\gamma(t_{i_j},p)\)
\end{center}
For \(\gamma(t_{i_j},p)\in \Lim(I)\), use condition (1) of  definition \ref{Definition of D}. to conclude that \(\gamma(t_{i_{j-1}},p)+\omega^{o^{\vec{U}}(t_{i_j})}=\gamma(t_{i_{j}},p)\). 
Thus
\begin{center}
 \(\forall r\in I\cap(\gamma(t_{i_{j-1}},p),\gamma(t_{i_j},p))\)  (\(o_L(r)<o^{\vec{U}}(t_{i_j})\))
\end{center}
In Particular,
\begin{center}
 \(I(t'_{i_j},\pi_I(p))=\min(\beta\in I\setminus \gamma(t_{i_{j-1}},p)+1\mid o_L(\beta)=o^{\vec{U}}(t_{i_j}))=\gamma(t_{i_j},p)\)
\end{center}
\begin{flushright}
 \(\blacksquare_{of \ claim}\)
\end{flushright}
Checking definition \ref{Definition of condition of MIU}, show that \(\pi_I(p)\in \MforI\): (1), (2.a.i), (2.b.i), (2.b.iii) are immediate from the definition of \(\pi_I\). Use the claim to verify that (2.a.ii), (2.b.ii) follows from (1),(2) in \(D\) respectively. For (2.a.iii), let \(1\leq j\leq n'\), write
\begin{center}
 \(\gamma(t_{i_{j-1}},p)+\underset{i_{j-1}<l\leq i_j }{\sum}\omega^{o^{\vec{U}}(t_l)}=\gamma(t_{i_j},p)\)
\end{center}
This equation induces a C.N.F equation
\begin{center}
 \(I(t_{i_{j-1}},\pi_I(p))+ \sum_{k=1}^{n_0}\omega^{o^{\vec{U}}(t_{l_k})}=I(t_{i_j},\pi_I(p))\) \ (C.N.F)
\end{center}
Thus 
\begin{center}
 \(\langle\kappa(t_{l_1}),...,\kappa(t_{l_{n_0-1}})\rangle\in Y(o^{\vec{U}}(t_{l_1}))\times ... \times Y(o^{\vec{U}}(t_{l_{n_0-1}}))\bigcap[(\kappa(t_{i_{j-1}}),\kappa(t_{i_j}))]^{<\omega}\)
\end{center}
For $(1)$, let \(q=\langle t'_1,...,t'_{n+1}\rangle\in \MforI\). 
For every \(t'_j\) such that \(I(t'_j,q)\in \succ(I)\), use definition \ref{Definition of condition of MIU} (2.a.iii) to find \(\vec{s}_j=\langle s_{j,1},...,s_{j,m_j}\rangle\) such that
\begin{center}
\( \langle\kappa(s_{j,1}),...,\kappa(s_{j,m_j})\rangle\in  Y(\gamma_1)\times ... \times Y(\gamma_{m-1})\bigcap[(\kappa(t'_{i_{r}-1}),\kappa(t'_{i_r}))]^{<\omega}\)
\end{center}
where 
\(I(t'_{i_{r}-1},q)+\overset{m}{\underset{i=1}{\sum}}\omega^{\gamma_i}=I(t'_{i_r},q)\) (C.N.F).\\
For each \(i=1,...,n\) such that \(o^{\vec{U}}(t'_i)>0\) and \(\kappa(t'_i)\in \succ(I)\) pick some
\(B(t'_i)\in\underset{\xi<o^{\vec{U}(t_i')}}{\bigcap}U(t_i,\xi)\).
Define \(p=\langle t_1,...,t_{n+1}\rangle^{\frown}\langle \vec{s}_r\mid I(t_r,q)\in\succ(I)\rangle\)
\begin{center}
\(t_i= \begin{cases} 
\langle \kappa(t'_i),B(t'_i)\setminus \kappa(s_{i,m_i})+1\rangle & o^{\vec{U}}(t'_i)>0 \\
\kappa(t'_i) & otherwise
\end{cases}\)
\end{center}
Once we prove that \(\gamma(s_{r,j},p)\notin I\) and that \(p\) computes \(I\) correctly  i.e. \(\gamma(t_i,p)=I(t'_i,q)\), it will follow that \(\pi_I(p)=\langle t'_i \mid \gamma(t_i,p)\in I\rangle=q\). By induction on \(i\), for \(i=0\) it is trivial. Let \(0<i\) and assume the statement holds for i. If \(I(t'_{i+1},q)\in\Lim(I)\), then by \ref{Definition of condition of MIU} (b.ii)
\begin{center}
\(I(t'_{i+1},q)=I(t'_{i},q)+\omega^{o^{\vec{U}}(t'_{i+1})}=\gamma(t_i,p)+\omega^{o^{\vec{U}}(t_{i+1})}=\gamma(t_{i+1},p)\)
\end{center}
If \(I(t'_{i+1},q)\in\succ(I)\), then from \ref{Definition of condition of MIU} (a.ii) it follows that  \(I(t'_i,q)\) is the predecessor of \(I(t'_{i+1},q)\). By the choice of \(\vec{s}_{i+1}\), 
\begin{center}
 \(\gamma(t_{i+1},p)=\gamma(t_{i},p)+ \overset{m-1}{\underset{i=1}{\sum}}\omega^{\gamma_1}n_i+\omega^{\gamma_m}(n_{m}-1)+\omega^{o^{\vec{U}(t_{i+1})}}=\) \\ \(=I(t'_i,q)+ \overset{m-1}{\underset{i=1}{\sum}}\omega^{\gamma_1}n_i+\omega^{m_1}(n_{m_1}-1)+\omega^{o^{\vec{U}(t'_{i+1})}}=I(t'_{i+1},q)\)
\end{center} 
Also, for all \(1\leq r\leq m_{i+1}\),  \(\gamma(s_{i+1,r},p)\) is between two successor ordinals in \(I\), hence \(\gamma(s_{i+1,r},p)\notin I\). Finally,  \(p\in D\) follows from \ref{Definition of order of MIU} (a.ii) and \ref{Definition of D} condition (1). If \(\gamma(t_{i},p)\in \Lim(I)\) we did not add \(\vec{s}_{i}\). Thus \(i_{j-1}=i_j-1\).

For $(2)$, assume that \(p,q\in D\), \(p\leq q\). Using the claim, the verification of definition \ref{Definition of order of MIU} it similar to $(1)$.

As for $(3)$, let us prove it for a simpler case to ease the notation. Nevertheless, the general statement if very similar and only require suitable notation. Let \(p=\langle t_1,...,t_{n+1}\rangle\in \Mfor\). Assume that
$$\pi_I(p)=\langle t'_{i_1},...,t'_{i_{n'}}\rangle \leq_I\langle t'_{i_1},...,t'_{i_{j-1}},s_{1},..,s_{m},t'_{i_j},...,t'_{i_n}\rangle=q'\in \MforI$$
For every \(l=1,...,m\) such that \(I(s_l,\pi_I(p))\in \succ(I)\) use definition \ref{Definition of order of MIU} (2b) to find \(\vec{s_l}=\langle s_{l,1},...,s_{l,m_l}\rangle\) such that
\begin{center}
\( \langle\kappa(s_{l,1}),...,\kappa(s_{l,m_l})\rangle\in  B(t_{i_j},\gamma_1)\times ... \times B(t_{i_j},\gamma_{m-1})\bigcap[(\kappa(s_{l-1}),\kappa(s_l))]^{<\omega}\)
\end{center}
where 
\(I(s_{l-1},\pi_I(p))+\overset{m}{\underset{i=1}{\sum}}\omega^{\gamma_i}=I(s_l,\pi_I(p))\) (C.N.F). Define \(p\leq p'\) to be the extension
\(p'=p^{\frown}\langle s'_{1},..,,s'_{m}\rangle^{\frown}\langle \vec{s_l}\mid I(s_l,\pi_I(p))\in\succ(I)\rangle\) where
\begin{center}
\(s'_i= \begin{cases} 
\langle \kappa(s_i),B_i\setminus \kappa(s_{i,m_i})+1\rangle & o^{\vec{U}}(s_i)>0 \\
s_i & otherwise
\end{cases}\)
\end{center}
 As in (1), \(\pi_I(p')=\langle t'_{i_1},...,t'_{i_{j-1}},(s_1')',...,(s_m')',...t_{i_{n'}}\rangle\). Notice that since we only change \(s_l\) such that \(I(s_l,\pi_I(p))\in\succ(I)\), \((s_l')'=s_l\). Thus \(\pi_I(p')=q\) and \(p'\in D\) follows.
 \begin{flushright}
 \(\blacksquare\)
 \end{flushright}
 From the discussion previous to \ref{Definition of D}, we have to following corollary:  
 \begin{corollary}\label{Cor to generate subsequence}
 Let  \(G\subseteq\Mfor\) be a $V$-generic filter, and let $C'\subseteq C_G$ be a closed subset, then there is $G_I\subseteq \MforI$ such that $V[C']=V[G_I]$ and $C_I=C'$, where $I=\Inx(C',C_G)$.
\end{corollary}
 \begin{flushright}
 \(\blacksquare\)
 \end{flushright}
 
\subsection{The Quotient forcing $\Mfor/G_I$}
\begin{definition}\label{Definition of quotient}
Let \(G_I\) be \(\MforI\) generic, the \textit{quotient forcing} is
\begin{center}

\(\Mfor/G_I=\pi^{-1''}_IG_I=\{p\in \Mfor\mid \pi_I(p)\in G_I\}\)
 
\end{center}
\end{definition}

 The forcing \(\Mfor/G_I\) completes \(V[G_I]\) to \(V[G]\) in the sense that if \(G\subseteq\Mfor\) is $V$-generic, and \( \pi^*_I(G)=G_I\), then \(G\subseteq\Mfor/G_I\) is $V[G_I]$-generic. Moreover, if $G\subseteq\Mfor/G_I$ is $V[G_I]$-generic, then $G\subseteq \Mfor$ is $V$-generic, and $\pi^*_I(G)=G_I$.
 \begin{corollary}\label{ Cor. copmpleting CI to CG}
  Let $G_I\subseteq\MforI$ be $V$-generic, then there is $G\subseteq \Mfor$ such that $C_G\restriction I=C_I$.
 \end{corollary}
 The following proposition is straightforward: 
\begin{proposition}\label{forcing to be in quotient}
  \textit{Let \(x,p\in\Mfor\) and \(q\in\MforI\), then}
\begin{enumerate}
    
\item\(\pi_I(p)\leq_I q\Rightarrow q\Vdash_{\MforI}\overset{\vee}{p}\in \Mfor/\underaccent{\sim}{G}_I\).
\item \(q\Vdash_{\MforI}\overset{\vee}{p}\in \Mfor/\underaccent{\sim}{G}_I \Rightarrow \pi_I(p), q \text{ are  compatible} \).
\item \(x\Vdash_{\Mfor}\overset{\vee}{p}\in \Mfor/\underaccent{\sim}{G}_I \Rightarrow \pi_I(p), \pi_I(x) \ are \ compatible\).
\end{enumerate}
 \begin{flushright}
 \(\blacksquare\)
 \end{flushright}
 
\end{proposition}
\begin{lemma}\label{quotient is kappaplus cc in GI}
 \textit{Let \(G_I\subseteq\MforI\) be $V$-generic. Then the forcing \(\Mfor/G_I\) satisfies \(\kappa^+-c.c.\) in \( V[G_I]\).}
 
\end{lemma}
\pr Fix \(\{p_\alpha \mid \alpha<\kappa^+\}\subseteq \Mfor/G_I\) and let 
$$r\in G_I ,\ r\Vdash_{\MforI} \forall\alpha<\kappa^+ \ \ \underaccent{\sim}{p}_\alpha\in  \Mfor/\underaccent{\sim}{G}_I$$
Let us argue that
\begin{center}
\(E=\{q\in \MforI \mid (q\bot r)\bigvee(q\Vdash_{\MforI}\exists\alpha,\beta<\kappa^+ \  (\underaccent{\sim}{p}_\alpha,\underaccent{\sim}{p}_\beta\ are \ compatible )\}\)
\end{center}
is a dense subset of \(\MforI\). Assume \(r\leq_I r'\), for every \(\alpha<\kappa^+\) pick some \(r'\leq_I q^*_{\alpha}\in \MforI,\ p^*_\alpha \in \Mfor\) such that:
\begin{enumerate}
    \item \(\pi_I(p^*_\alpha)=q^*_\alpha\).
    \item \(q^*_\alpha\Vdash\underaccent{\sim}{p}_\alpha\leq \overset{\vee}{p^*_\alpha}\in \Mfor/\underaccent{\sim}{G}_I\).
\end{enumerate}
There exists such \( q^*_{\alpha}\ , p^*_\alpha\): Find \(r'\leq_I q'_{\alpha}\) and \(p'_\alpha\) such that \(q'_{\alpha}\Vdash \overset{\vee}{p'_\alpha}=\underaccent{\sim}{p}_\alpha\) then by the proposition \ref{forcing to be in quotient} (2), there is \(q^*_{\alpha}\geq_I \pi_I(p'_\alpha),q'_\alpha\). By lemma \ref{pi is a projection} (3) there is \(p^*_\alpha\geq p'_\alpha\) such that \(q^*_\alpha:=\pi_I(p^*_\alpha)\). It follows from proposition \ref{forcing to be in quotient} (1) that
\begin{center}

\(q^*_\alpha\Vdash  \underaccent{\sim}{p}_\alpha\leq\overset{\vee}{p^*_\alpha}\in \Mfor/\underaccent{\sim}{G}_I\)
\end{center} 
Denote
\(p^*_\alpha=\langle t_{1,\alpha},...,t_{n_\alpha,\alpha},t_{n_\alpha+1,\alpha}\rangle\),
\(q^*_\alpha=\langle t_{i_1,\alpha},...,t_{i_{m_\alpha},\alpha},t_{n_\alpha+1,\alpha}\rangle\).
Find \(S\subseteq\kappa^+ , \ n<\omega\) and \(\langle \kappa_1,...,\kappa_n\rangle\) such that  \(|S|=\kappa^+ \) and for any \(\alpha \in S \), \(n_\alpha=n\) and
\begin{center}

\(\langle\kappa(t_{1,\alpha}),...,\kappa(t_{n_\alpha,\alpha})\rangle=\langle\kappa_1,...,\kappa_n\rangle\).
\end{center}
Since \(\pi_I(p^*_\alpha)=q^*_\alpha\) it follows that
\begin{center}

\(\langle\kappa(t_{i_1,\alpha}),...,\kappa(t_{i_{m_\alpha},\alpha})\rangle=\langle\kappa_{i_1},...,\kappa_{i_m}\rangle\)
\end{center}
for some \(m<\omega\) and \(1\leq i_1<...<i_m\leq n\). \\
Fix any \(\alpha,\beta \in S\) and let
\(p^*=\langle t_1,....,t_n,t_{n+1}\rangle\) 
 where
\begin{center}
\(t_i= \begin{cases} 
\langle \kappa_i,B(t_{i,\alpha})\cap B(t_{i,\beta})\rangle & o^{\vec{U}}(t_{i,\alpha})>0 \\
\kappa_i & otherwise
\end{cases}\) 
\end{center}
Denote \(p_{\alpha}^*\cap p^*_{\beta}=p^*\). Set
\begin{center}
 \(q^*=\pi_I(p^*)=\langle t'_{i_1},...,t'_{i_m}\rangle\)
\end{center}
Then \(r'\leq_I q^*_\alpha\cap q^*_\beta=\pi_I(p^*_\alpha)\cap\pi_I(p^*_\beta)=\pi_I(p^*_\alpha\cap p^*_\beta)=\pi_I(p^*)=q^*\).
It follows that \(q^* \in E\) since by proposition \ref{forcing to be in quotient} (1) \(q^*\Vdash_{\MforI}\overset{\vee}{p^*}\in \Mfor/\underaccent{\sim}{G}_I\) and
\begin{center}
\(q^*\Vdash_{\MforI} \underaccent{\sim}{p}_\alpha\leq \overset{\vee}{p^*_\alpha}\leq^*\overset{\vee}{p^*} \ \wedge \ \underaccent{\sim}{p}_\beta\leq \overset{\vee}{p^*_\beta}\leq^*\overset{\vee}{p^*}\)
\end{center}
The rest is routine.
\begin{flushright}
 \(\blacksquare\)
 \end{flushright}
\begin{lemma}\label{quotient is kappaplus cc in VG}
\textit{ Let \(G\subseteq\Mfor\) be $V$-generic. Then the forcing \(\Mfor/G_I\) satisfies \(\kappa^+-c.c.\) in \( V[G]\).}

\end{lemma}
\pr Fix \(\{p_\alpha \mid \alpha<\kappa^+\}\subseteq \Mfor/G_I\) in \(V[G]\) and let 
\begin{center}

\(r\in G ,\ r\Vdash_{\Mfor}  \forall\alpha<\kappa^+ \ \ \underaccent{\sim}{p}_\alpha\in  \Mfor/\underaccent{\sim}{G}_I\)
\end{center}
Similar to lemma \ref{quotient is kappaplus cc in GI} let us argue that
\begin{center}
\(E=\{x\in \Mfor \mid (q\bot r)\bigvee(q\Vdash_{\Mfor}\exists\alpha,\beta<\kappa^+ (\underaccent{\sim}{p}_\alpha,\underaccent{\sim}{p}_\beta)\ are \ compatible )\}\)
\end{center}
is a dense subset of \(\Mfor\). Assume \(r\leq r'\) , for every \(\alpha<\kappa^+\) pick some \(r'\leq x'_{\alpha}\in \Mfor\ ,\ p'_\alpha \in \Mfor\) such that
\(x'_\alpha\Vdash_{\Mfor}\underaccent{\sim}{p}_\alpha= \overset{\vee}{p'_\alpha}\). By proposition \ref{forcing to be in quotient} (3), we can find \(\pi_I(x'_\alpha),\pi_I(p'_\alpha)\leq_I y_\alpha \). 
By lemma \ref{pi is a projection} (3), There is \(x'_{\alpha}\leq x^*_{\alpha}\) ,\(p'_\alpha\leq p^*_\alpha\) such that 
$$\pi_I(x'_\alpha),\pi_I(p''_\alpha)\leq_Iy_\alpha=\pi_I(p^*_\alpha)=\pi_I(x^*_\alpha)$$
Denote by
$$x^*_\alpha=\langle s_{1_\alpha},...,s_{k_\alpha,\alpha},s_{k_\alpha+1,\alpha}\rangle. \ \ p^*_\alpha=\langle t_{1,\alpha},...,t_{n_\alpha,\alpha},t_{n_\alpha+1,\alpha}\rangle$$
and
$$\pi_I(x^*_\alpha)=\langle t'_{i_1,\alpha},...,t'_{i_{k'_\alpha},\alpha}t'_{k_{\alpha}+1}\rangle=\pi_I(p_\alpha)$$
Find \(S\subseteq\kappa^+ \ |S|=\kappa^+\) and \(\langle\kappa_1,...,\kappa_n\rangle,\langle\nu_1,...,\nu_k\rangle\) such that for any \(\alpha \in S \),
$$\langle\kappa(t_{1,\alpha}),...,\kappa(t_{n_\alpha,\alpha})\rangle=\langle\kappa_1,...,\kappa_n\rangle, \
\langle\kappa(s_{1,\alpha}),...,\kappa(s_{k,\alpha})\rangle=\langle\nu_1,...,\nu_k\rangle$$
Fix any \(\alpha,\beta \in S\) and let
\(p^*=p^*_\alpha\cap p^*_\beta\),
\(x^*=x^*_\alpha\cap x^*_\beta\).
Then \(p'_\alpha,p'_\beta\leq^* p^*\) and \(x_\alpha,x_\beta\leq^*_I x^*\) .
Finally claim that \(x^*\in E\):
$$\pi_I(p^*)=\pi_I(p^*_\alpha)\cap\pi_I(p^*_\beta)=\pi_I(x^*_\alpha)\cap\pi_I(x^*_\beta)=\pi_I(x^*)$$
 thus \(x^*\Vdash_{\Mfor} \overset{\vee}{p^*} \in \Mfor/\underaccent{\sim}{G}_I\). Moreover, \(x_\alpha \leq^* x^*\) which implies that \\ \(x^*\Vdash_{\Mfor} \overset{\vee}{p^*}\geq \underaccent{\sim}{p}_\alpha,\underaccent{\sim}{p}_\beta \).
 \begin{flushright}
 \(\blacksquare\)
 \end{flushright}
 \subsection{The Argument for General Sets}
 Let us conclude first the main result for subsets of $\kappa^+$. 
\begin{lemma}\label{Thm for subsets of kappaplus}
 \textit{If \(A\in V[G] \ ,\ A\subseteq\kappa^+\) then there exists \(C^*\subseteq C_G\) such that \(V[A]=V[C^*]\).}
 
\end{lemma}
\pr Work in \(V[G]\), for every \(\alpha<\kappa^+\) find subsequences \(C_\alpha\subseteq C_G\)  such that \(V[C_\alpha]=V[A\cap\alpha]\) using the induction hypothesis.
The function \(\alpha\mapsto C_\alpha\) has range \(P(C_G)\) and domain \(\kappa^+\) which is regular in \(V[G]\).
Therefore there exist \(E\subseteq\kappa^+\) unbounded in \(\kappa^+\) and \(\alpha^*<\kappa^+\)  such that for every \(\alpha\in E\), \(C_\alpha=C_{\alpha^*}\). Set \(C^*=C_{\alpha^*}\), then 
\begin{enumerate}
\item \( C^*\subseteq C_G\).
\item \(C^*\in V[A\cap\alpha^*]\subseteq V[A]\).
\item \(\forall\alpha<\kappa^+. A\cap\alpha\in V[C^*]\).
\end{enumerate}
Assume that \(C^*\) is a club\footnote{\(Cl(C^*)\) clearly  satisfy $(1)-(3)$}. Unlike \(A\)'s that were subsets of \(\kappa\), for which we added another piece of \(C_G\) to \(C^*\) to obtain \(C'\) such that \(V[A]=V[C']\), here we argue that \(V[A]=V[C^*]\).

By (2),  \(C^*\in V[A]\). For the other direction, denote by \(I\) the indexes of \(C^*\) in \(C\) and consider the forcings \(\MforI,\Mfor/G_I\). Toward a contradiction, assume that \(A\notin V[C^*]\), and let \(\underaccent{\sim}{A}\in V[C^*]\) be a \(\Mfor/G_I\)-name for \(A\), where \(\pi''_IG=G_I\). Work in \(V[G_I]\), by lemma \ref{Properties of MUI} (2), \(V[G_I]=V[C^*]\).  For every \(\alpha<\kappa^+\) define
 \begin{center}
 \(X_\alpha=\{B\subseteq\alpha \mid ||\underaccent{\sim}{A}\cap\alpha=B||\neq 0\}\)
 \end{center}
where the truth value is taken in \(RO(\Mfor/G_I)\)\footnote{ $RO(\mathbb{Q})$ denotes the complete boolean algebra of regular open subsets of  \(\mathbb{Q}\)}. By lemma \ref{quotient is kappaplus cc in GI},
$$\forall\alpha<\kappa^+ \ |X_\alpha|\leq\kappa$$
For every \(B\in X_\alpha\) define \(b(B)=||\underaccent{\sim}{A}\cap\alpha||\).
Assume that \(B'\in X_\beta\) and \(\alpha\leq\beta\) then \(B=B'\cap\alpha \in X_\alpha\). Switching to boolean algebra notation (\(p\leq_B q\) means \(p\) extends $q$) \(b(B')\leq_B b(B)\). Note that for such \(B,B'\) if \(b(B')<_B b(B)\), then there is 
$$0<p\leq_B(b(B)\setminus b(B'))\leq_Bb(B)$$
Therefore,
$$p\cap b(B')\leq_B (b(B)\setminus b(B'))\cap b(B')=0$$
Hence \(p\bot b(B')\). Work in V[G],  denote \(A_\alpha=A\cap \alpha\). Recall that
$$\forall \alpha<\kappa^+ \ A_\alpha\in V[C^*]$$
thus \(A_\alpha\in X_\alpha\). Consider the  \(\leq_B\)-non-increasing sequence \(\langle b(A_\alpha) \mid \alpha<\kappa^+\rangle\). If there exists some \(\gamma^*<\kappa^+\) on which the sequence stabilizes, define 
$$A'=\bigcup\{B\subseteq\kappa^+ \ | \ \exists\alpha \ b(A_{\gamma^*})\Vdash \underaccent{\sim}{A}\cap\alpha=B\}\in V[C^*]$$
To see that \(A'=A\) , notice that if \(B,B',\alpha,\alpha'\) are such that $\alpha\leq\alpha'$, and
\begin{center}
 \(b(A_{\gamma^*})\Vdash \underaccent{\sim}{A}\cap\alpha=B\), 
 \(b(A_{\gamma^*})\Vdash \underaccent{\sim}{A}\cap\alpha'=B'\)
\end{center}
then \(B'\cap\alpha=B\) otherwise, the non zero condition \(b(A_{\gamma^*})\) would force contradictory information. Consequently, for every
\(\xi<\kappa^+\) there exists \(\xi<\gamma<\kappa^+\) such that \(b(A_{\gamma^*})\Vdash \underaccent{\sim}{A}\cap\gamma=A\cap\gamma\), hence \(A'\cap\gamma= A\cap\gamma\). This is a contradiction to \(A\notin V[C^*]\).\\ Therefore, the sequence \(\langle b(A_\alpha) \mid \alpha<\kappa^+\rangle\) does not stabilize. By regularity of \(\kappa^+\), there exists a subsequence \(\langle b(A_{i_\alpha}) \mid \alpha<\kappa^+\rangle\) which is strictly decreasing. By the observation we made in the last paragraph, find \(p_\alpha\leq_B b(A_{i_\alpha}) \) such that \(p_\alpha \bot b(A_{i_{\alpha+1}})\). Since \(b(A_{i_\alpha})\) are decreasing, for any \(\beta>\alpha\) \  \(p_\alpha\bot b(A_{i_\beta})\) and in turn \(p_\alpha\bot p_\beta\). This shows that \(\langle p_\alpha \mid \alpha<\kappa^+\rangle\in V[G]\) is an anti chain of size \(\kappa^+\) which contradicts Lemma \ref{quotient is kappaplus cc in VG}.
Thus \(V[A]=V[C^*]\).
\begin{flushright}
 \(\blacksquare\)
 \end{flushright}
 \vskip 0.2 cm
\textbf{End of the proof of Theorem \ref{Main theorem}:} By induction on \(\sup(A)=\lambda>\kappa^+\). It suffices to assume that \(\lambda\) is a cardinal.
\vskip 0.2 cm
\underline{case1:} Assume \(cf^{V[G]}(\lambda)>\kappa\), then the arguments of lemma \ref{Thm for subsets of kappaplus} works.
\vskip 0.2 cm
\underline{case2:} Assume \(cf^{V[G]}(\lambda)\leq\kappa\), since \(\Mfor\) satisfies \(\kappa^+-c.c.\) we must have that   \(\nu:=cf^V(\lambda)\leq\kappa\). Fix \(\langle\gamma_i | \ i<\nu\rangle\in V\) cofinal in \(\lambda\). Work in \(V[A]\), for every \(i<\nu\) find \(d_i\subseteq \kappa\) such that \(V[d_i]=V[A\cap\gamma_i]\). By induction, there exists \(C^*\subseteq C_G\) such that \(V[\langle d_i\mid i<\nu\rangle]=V[C^*]\), therefore 
\begin{enumerate}
\item \(\forall i<\nu \ A\cap\gamma_i\in V[C^*]\).
\item \(C^*\in V[A]\).
\end{enumerate}
Work in \(V[C^*]\), for \(i<\nu\) define
 \(X_i=\{B\subseteq\alpha \mid\ ||\underaccent{\sim}{A}\cap\gamma_i=B||\neq 0\}\).
By lemma \ref{quotient is kappaplus cc in GI}, \(|X_i|\leq\kappa\).  
For every \(i<\nu\) fix an enumeration 
$$X_i=\langle X(i,\xi) \mid\ \xi<\kappa\rangle\in V[C^*]$$
There exists \(\xi_i<\kappa\) such that \(A\cap\gamma_i=X(i,\xi_i)\). Moreover, since \(\nu\leq\kappa\) the sequence \(\langle A\cap\gamma_i \mid i<\nu\rangle=\langle X(i,\xi_i) \mid i<\nu\rangle\) can be coded in \(V[C^*]\) as a sequence of ordinals below \(\kappa\). By induction there exists \(C''\subseteq C_G\) such that \(V[C'']=V[\langle\xi_i \mid i<\nu\rangle]\). It follows that,
$$V[C'',C^*]=(V[C^*])[\langle\xi_i \mid i<\nu\rangle]=V[A]$$
Finally, we can take for example, \(C'=C''\cup C^*\subseteq C_G\) to obtain
\(V[A]=V[C']\)
\begin{flushright}
 \(\blacksquare_{theorem \ \ref{Main theorem}}\)
 \end{flushright}
 
\newpage
\section{Classification of subforcing of Magidor}

Now that we have Classified models of the form $V[A]$, we can conclude the following:
\begin{corollary}
 Let $G\subseteq\Mfor$ be a $V$-generic filter, and let $M$ be a transitive model of $ZFC$ such that $V\subseteq M\subseteq V[G]$. Then there is $C_M\subseteq C_G$ such that $V[C_M]=M$.
\end{corollary}
\pr By \cite[Thm. 15.43]{Jech}, there is $D\in V[G]$ such that $M=V[D]$. By theorem \ref{Main theorem}, there is $C_M\subseteq C_G$ such that $V[C_M]=V[D]=M$.

As we have seen in the previous section, the models $V[C_M]$ are generic extensions for the forcings $\MforI$ which in turn are projection of $\Mfor$, this yield the classification of subforcings. Although the classification can naturally be extended to a the class of forcings $\Mforc$, we present here only the classification of subforcings of $\Mfor$.
\begin{definition}
 Recall definition \ref{Definition of condition of MIU} of $\MforI$. The forcings $$\{\MforI \mid I\in P(\omega^{o^{\vec{U}}(\kappa)}), \ I \text{ is closed}\}$$ is the family of \textit{Magidor-type forcings} with the coherent sequence \(\vec{U}\).
 
\end{definition}

 In practice, Magidor-type forcings are just Magidor forcing with a subsequence of \(\vec{U}\); If \(I\) is any closed subset of indices, we can read the measures of \(\vec{U}\) from which the elements of the final sequence are chosen from, using the sequence \(\langle o_L(i)\mid i\in I\rangle\) (recall that \(o_L(i)=\gamma_n\) where \(i=\omega^{\gamma_1}+...+\omega^{\gamma_n}\) \textit{C.N.F}).\\
 \\
\textbf{Example}: Assume that \(o^{\vec{U}}(\kappa)=2\) and let a\begin{center}
  \(I=\{1,\omega,\omega+1\}\cup(\omega\cdot3\setminus\omega\cdot2)\cup\{\omega\cdot3,\omega\cdot4,...\}\in P(\omega^2)\)
 \end{center}
 Then \(\langle o_L(i)\mid i\in I\rangle=\langle0,1,\underset{\omega}{\underbrace{0,0,0...}},\underset{\omega}{\underbrace{1,1,1...}}\rangle\). Therefore \(\MforI\) is just Prikry foricing with \(U(\kappa_1,0)\) for some measurable \(\kappa_1<\kappa\) followed by Prikry forcing with \(U(\kappa,1)\).\\
 Although in this example the noise at the beginning can be neglected, there are \(I\)'s for which we do not get "pure" Magidor forcing which uses one measure at a time and combine several measure. For example we can obtain the Tree-Prikry forcing- let $I=\langle \omega^n\mid n<\omega\rangle$ then $\langle o_L(i)\mid i\in I\rangle=\langle n\mid n<\omega$, conditions in the forcing are of the form $\langle t_1,...,t_n,\langle\kappa,B\rangle\rangle$ the extensions is from the measures $U(\kappa,m), \ m>n$ which is essentially $P_T(\langle U(\kappa,n)\mid n<\omega\rangle)$ the tree Prikry forcing such that at level $n$ the tree splits on a large set in $U(\kappa,n)$.
 \vskip 0.4 cm

 

 For the definition of complete subforcing see \cite{Shelah}. 
\begin{theorem}\label{classification}
\textit{ Let \(\mathbb{P}\subseteq \Mfor\) be a complete subforcing of \(\Mfor\) then there exists a maximal antichain \(Z\subseteq\mathbb{P}\) and \(\l I_p\mid  p\in Z\r\) such that \(\mathbb{P}_{\geq p}\) (the forcing \(\mathbb{P}\) above \(p\)) is forcing equivalent to the Magidor-type forcing \(\mathbb{M}_{I_p}[\vec{U}]_{\geq q_p}\) for some condition $q_p\in RO(\mathbb{M}_{I_p}[\vec{U}])$.}
 
\end{theorem}
\pr Let \(H\subseteq \mathbb{P}\) be generic, then there exists \(G\subseteq\Mfor\) generic such that \(H=G\cap\mathbb{P}\), in particular \(V\subseteq V[H]\subseteq V[G]\). By Theorem \ref{Main theorem}, there is a closed \(C'\subseteq C_G\) such that \(V[C']=V[H]\), and let \(I=\Inx(C',C_G)\). The assumption \(o^{\vec{U}}(\kappa)\) is crucial to claim that \(I\in V\). By  corollary \ref{Cor to generate subsequence}, there is $G_I\subseteq \MforI$ such that \(C_I=C'\). Let \(\underaccent{\sim}{C'}_0,\lusim{H}_0\) be a \(\mathbb{P}\)-name of \(C',H\).
  Let \(p\in \mathbb{P}\) such that
\begin{center}
 \(p\Vdash  \underaccent{\sim}{C'}\) is generic sequence for \(\MforI\) and \(V[\underaccent{\sim}{H}]=V[\underaccent{\sim}{C'}]\).\footnote{This is indeed a formula in the forcing language since for any set \(A\), \(V[A]=\underset{z\subseteq ord,z\in V}{\bigcup}L[z,A]\) where \(L[z,A]\) is the class of all constructible sets relative to \(z,A\). }
\end{center}
Denote by $I_p:=I$. For the other direction, let \(\lusim{C'}_1,\underaccent{\sim}{H}_1\) be \(\mathbb{M}_{I_p}[\vec{U}]\)-names for \(C',H\) and let \(q_p\in RO(\mathbb{M}_{I_p}[\vec{U}])\) be the truth value:
 $$\underaccent{\sim}{H}_1\subseteq\mathbb{P}\text{ is }V-\text{generic }, \ p\in \underaccent{\sim}{H}_1\text{ and } V[\underaccent{\sim}{H}_1]=V[\underaccent{\sim}{C'}_1]$$
Clearly, \(\mathbb{M}_{I_p}[\vec{U}]_{\geq q_p}\) and \(\mathbb{P}_{\geq p}\) have the same generic extensions and therefore forcing equivalent.
\begin{flushright}
 \(\blacksquare\)
 \end{flushright}

\newpage

\section{Prikry forcings with non-normal ultrafilters.  }

Let $\kappa$ be a measurable cardinal and  let $\mathbb{U}=\l U_a \mid a \in [\kappa]^{<\omega}\r$ be a tree consisting of $\kappa-$complete non-trivial ultrafilter
over $\kappa$.

Recall the definition due to Prikry of the tree Prikry forcing with $\mathbb{U}$.

\begin{definition}
The \textit{Tree Prikry forcing} $P(\mathbb{U})$, consist of all pairs $\l p, T \r$ such that:
\begin{enumerate}
    \item $p$ is a finite increasing sequence of ordinals below $\kappa$.
    \item $T \subseteq [\kappa]^{<\omega}$ is a tree with trunk $p$ such that
   \\for every $q \in T$ with $q\geq_T p$,
   the set of the immediate successors of $q$ in $T$, i.e.
   $Suc_T(q)$ is in $U_q$.
  \end{enumerate}
\end{definition}

The orders $\leq,\leq^*$ are defined in the usual fashion.


For every $a \in [\kappa]^{<\omega}$, let $\pi_a$ be a projection of $U_a$ to a normal ultrafilter. Namely,
let $\pi_a:\kappa\to \kappa$ be a function which represents $\kappa$ in the ultrapower by $U_a$, i.e. $[\pi_a]_{U_a}=\kappa$.
Once $U_a$ is a normal ultrafilter, then let $\pi_a$ be the identity.

By passing to a dense subset of $P(\mathbb{U})$, we can assume that for each $\l p, T \r \in P(\mathbb{U})$, for every $\l \nu_1,..., \nu_n\r \in T$
we have $$\nu_1<\pi_{\l\nu_1 \r}(\nu_2)\leq \nu_2<...\leq \nu_{n-1}<\pi_{\l \nu_1,..., \nu_{n-1}\r}(\nu_n)$$
and for every $\nu\in Suc_T(\l \nu_1,..., \nu_n\r), \pi_{\l \nu_1,..., \nu_{n}\r}(\nu)>\nu_{n}.$

Note that once the measures over a certain level (or certain levels) are the same - say for some $n<\omega$ and $U$, for every $a \in [\kappa]^n$,
$U_a=U$, then a modified diagonal intersection

$$\Delta^*_{\alpha<\kappa}A_\alpha:=\{\nu<\kappa \mid \forall \alpha<\pi_k(\nu) (\nu\in A_\alpha)\}\in U, $$ once
$\{ A_\alpha \mid \alpha<\kappa \}\subseteq U$, can be used to avoid or to simplify the tree structure.

For example, if $\l {\cal V}_n \mid n<\omega\r$ is a sequence of $\kappa-$complete ultrafilters over $\kappa$, then the Prikry forcing with it
$P(\l {\cal V}_n \mid n<\omega\r)$ is defined as follows:

\begin{definition}
The \textit{tree Prikry forcing with an $\omega$-sequence of ultrafilters}, $P(\l {\cal V}_n \mid n<\omega\r)$, is the set of all pairs $\l p, \l A_n \mid |p|<n<\omega\r \r$ such that:
\begin{enumerate}
    \item $p=\l \nu_1,...,\nu_k\r$ is a finite sequence of ordinals below  $\kappa$, such that \\$\nu_j<\pi_i(\nu_i)$, whenever $1\leq j<i\leq k$.
    \item $A_n \in {\cal V}_n$, for every $n,|p|<n<\omega$.
    \item\label{diag-sq} $\pi_{k+1}(\min(A_{k+1}))>\max(p)$, where
    $\pi_n:\kappa\to \kappa$ is a projection of ${\cal V}_n$ to a normal ultrafilter, i.e. $\pi_n$ is a function which represents $\kappa$ in the ultrapower by ${\cal V}_n$,
    $[\pi]_{{\cal V}_n}=\kappa$.

  \end{enumerate}

\end{definition}

A simpler case is once all ${\cal V}_n$ are the same, say all of them are $U$.
Then we will have the Prikry forcing with $U$:

\begin{definition}\label{def-PU}
The \textit{Prikry forcing with general ultrafilter $P(U)$} is the set of all pairs $\l p, A \r$ such that
\begin{enumerate}
    \item $p=\l \nu_1,...,\nu_k\r$ is a finite sequence of ordinals below  $\kappa$, such that \\$\nu_j<\pi(\nu_i)$, whenever $1\leq j<i\leq k$.
    \item $A \in U$.
    \item\label{diag} $\pi(\min(A))>\max(p)$,
    where $\pi$ is a projection of $U$ to a normal ultrafilter.
  \end{enumerate}

\end{definition}

Let $G$ be a generic for $\l P(\mathbb{U}), \leq \r$. Set
$$C= \bigcup\{p \mid \exists T \quad \l p, T \r \in G \}.$$
It is called a Prikry sequence for $\mathbb{U}$.

 For every natural $n\geq 1$ we would like to define a $\kappa-$complete ultrafilter $U_n$ over $[\kappa]^n$ which correspond to the first $n-$levels
 of trees in $P(\mathbb{U})$.
\\If $n=1$, set $U_1=U_{\l\r}$.
\\ Deal with the next step $n=2$.
Here for each $\nu<\kappa$ we have $U_{\nu}$.
\\Consider the ultrapower by $U_{1}$:

$$i_1:=i_{\l\r}: V\to M_{\l\r}.$$

Then the sequence $i_{\l\r}(\l U_{\l \nu \r} \mid \nu<\kappa\r)$ will have the length $i_{\l\r}(\kappa)$.
\\Let $U_{\l [id]_{U_{\l\r}}\r} $ be its $[id]_{U_{\l\r}}$ ultrafilter in $M_{\l\r}$ over $i_{\l\r}(\kappa)$.
Consider its ultrapower
$$i_{U_{\l [id]_{U_{\l\r}}\r}}:M_{\l\r}\to M_{\l [id]_{U_{\l\r}}\r}$$
Set $$i_2=i_{U_{\l [id]_{U_{\l\r}}\r}}\circ i_{\l\r}.$$
Then $$i_2:V \to  M_{\l [id]_{U_{\l\r}}\r}.$$
Note that if all of $U_{\l \nu \r}$'s are the same or just for a set of $\nu$'s in $U_{\l\r}$ they are the same, then this is just an ultrapower
by the product of $U_{\l\r}$ with this ultrafilter.
In general it is an ultrapower by $$U_{\l\r}-Lim \l U_{\l \nu \r} \mid \nu<\kappa\r,$$
where
$$X \in U_{\l\r}-Lim \l U_{\l \nu \r} \mid \nu<\kappa\r \text{ iff } [id]_{U_{\l [id]_{U_{\l\r}}\r}} \in i_2(X).$$
Note that once most of $U_{\l \nu \r}$'s are normal, then $U_{\l [id]_{U_{\l\r}}\r}$ is normal as well, and so,
$[id]_{U_{\l [id]_{U_{\l\r}}\r}} =i_{\l\r}(\kappa)$.

Define an ultrafilter $U_2$ on $[\kappa]^2$ as follows:

$$X \in U_2 \text{ iff } \l [id]_{U_{\l\r}},[id]_{U_{\l [id]_{U_{\l\r}}\r}}\r \in i_2(X).$$

Define also for $k=1,2$, ultrafilters $U_2^k$ over $\kappa$ as follows:

$$X \in U_2^1 \text{ iff }  [id]_{U_{\l\r}}\in i_2(X),$$

$$X \in U_2^2 \text{ iff }  [id]_{U_{\l [id]_{U_{\l\r}}\r}}\in i_2(X).$$

 Clearly, then $U_2^1 =U_1$ and $U_2^2= U_{1}-Lim \l U_{\l \nu \r} \mid \nu<\kappa\r$.
 Also $U_2^1$ is the projection of $U_2$ to the first coordinate and $U_2^2$ to the second.

Let $\l\l\r, T\r \in P(\mathbb{U})$.
It is not hard to see that $T\upr 2 \in U_2$.

Continue  and define in the similar fashion the ultrafilter $U_n$ over $[\kappa]^n$ and its projections to the coordinates $U_n^k$ for every $n>2, 1\leq k\leq n$.
We will have that for any $\l\l\r, T\r \in P(\mathbb{U})$, $T\upr n \in U_n$.
Also, if $1\leq n\leq m<\omega$, then the natural projection of $U_m$ to $[\kappa]^n$ will be $U_n$.


It is easy to see that $C$ is a Prikry sequence for  $\l U_n^n\mid 1\leq n<\omega\r$, in a sense that for every sequence $\l A_n \mid n<\omega\r \in V$,
with $A_n \in U_n^n$, there is $n_0<\omega$ such that for every $n>n_0$, $C(n) \in U_n^n$.
\\However, it does not mean that $C$ is generic for the forcing $P( \l U_n^n\mid 1\leq n<\omega\r)$ defined above (Definition 6.2).
The problem is with projection to normal. All $U_n^n$'s have the same normal $U_1$.

Suppose now that we have an ultrafilter $W$ over $[\kappa]^{\ell}$ which is Rudin-Keisler below some $\mathfrak{V}$ over $[\kappa]^k$ ($W\leq_{RK} \mathfrak{V}$), for some $k,\ell,1\leq \ell,k<\omega$.
This means that there is a  function $F:[\kappa]^k\to [\kappa]^{\ell}$ such that
$$X \in W \text{ iff } F^{-1}{}''X \in \mathfrak{V}.$$
So $F$ projects $\mathfrak{V}$ to $W$.
Let us denote this by $W=F_*\mathfrak{V}$.

The next statement characterizes $\omega-$sequences in $V[C]$.

\begin{theorem}\label{lem1}
Let $\l \alpha_k \mid k<\omega \r \in V[C]$ be an
increasing cofinal in $\kappa$ sequence. Then $\l \alpha_k \mid k<\omega \r$ is a Prikry sequence for a sequence in $V$ of $\kappa-$complete ultrafilters which are Rudin -Keisler below $\l U_n \mid n<\omega\r$.\footnote{ Let $\l {\cal V}_k \mid k<\omega \r$ be such sequence of ultrafilters over $\kappa$. We do not claim that $\l \alpha_k \mid k<\omega \r$ is Prikry generic for the forcing $P(\l {\cal V}_k \mid k<\omega\r)$, but rather that for every sequence $\l A_k \mid k<\omega\r \in V$,
with $A_k \in {\cal V}_k$, there is $k_0<\omega$ such that for every $k>k_0$, $\alpha_k \in {\cal V}_k $.}
\\Moreover,
there exist a non-decreasing  sequence of natural numbers $\l n_k \mid k<\omega\r$
and a sequence of functions $\l F_k \mid k<\omega \r$ in $V$,  $F_k:[\kappa]^{n_k} \to \kappa$, ($k<\omega$),
such that
\begin{enumerate}

\item
 $\alpha_k=F_k(C\upr n_k),$ for every $k<\omega$.

\item Let $\l n_{k_i} \mid i<\omega \r$ be the increasing subsequence of $\l n_k \mid k<\omega\r$ such that:
\begin{enumerate}
  \item $\{ n_{k_i} \mid i<\omega\}=\{ n_k \mid k<\omega\}$.
  \item $k_i=\min\{k \mid n_k=n_{k_i}\}$.
\end{enumerate}
Set $\ell_i=|\{k \mid n_k=n_{k_i}\}|$.
Then
$\l F_k(C\upr n_{k_i})\mid i<\omega,n_k=n_{k_i} \r $ will be a Prikry sequence for  $\l W_i \mid  i<\omega\r$,
i.e.  for every sequence $\l A_i \mid i<\omega\r \in V$,
with $A_i \in W_i$, there is $i_0<\omega$ such that for every $i>i_0$, $\l F_k(C\upr n_{k_i})\mid i<\omega,n_k=n_{k_i} \r \in A_i$,
where each $W_i$ is an ultrafilter over
$[\kappa]^{\ell_i}$ which is the projection of $U_{n_{k_i}}$ by $\l F_{k_i},..., F_{k_i+\ell_i-1} \r$.

\end{enumerate}

\end{theorem}
\pr Work in $V$. Given a condition $\l q,  S \r$, we will construct by
induction, using the Prikry property of the forcing $P(\mathbb{U})$,  a stronger
condition $ \l p, T \r$ which decides $\lusim{\alpha}_k$ once going up to a certain level $n_k$ of $T$.
Let us assume for simplicity that $q$ is the empty sequence.

Build by induction $\l \l\r, T \r\geq^* \l \l\r, S \r$
and a non-decreasing  sequence of natural
numbers $\l n_k \mid k< \omega \r$ such that for every $k<\omega$
\begin{enumerate}
  \item for every $\l \eta_{1},...,\eta_{n_k} \r \in T$  there is $\rho_{\l \eta_{1},...,\eta_{n_k} \r }<\kappa$ such that:
  \begin{enumerate}
    \item
  $\l  \l \eta_1,...,\eta_{n_k} \r, T_{\l \eta_1,...,\eta_{n_k} \r}\Vdash\lusim{\alpha}_k=\rho_{\l \eta_{1},...,\eta_{n_k}} \r$.
    \item $\rho_{\l \eta_{1},...,\eta_{n_k} \r} \geq \pi_{\l \eta_{1},...,\eta_{n_{k-1}} \r}(\eta_{n_k})$.
  \end{enumerate}

\item
there is no $n, n_k\leq n<n_{k+1}$ such that for some $\l
\eta_{1},...,\eta_n \r  \in T$ and ${E}$, the condition $\l\l
 \eta_1,...,\eta_{n} \r, {E} \r $ decides the value of
$\lusim{\alpha}_{{k+1}}$,

\end{enumerate}

Now, using the density argument and making finitely many changes, if necessary, we can assume that such $\l \l\r, T\r$ in
the generic set.

For every $k<\omega$,
define a function $F_k:Lev_{n_k}(T)\to \kappa$ by
setting $$F_k(\eta_1,...,\eta_{n_k})= \nu \text{ if } \l
 \l \eta_1,...,\eta_{n_k} \r, T_{\l \eta_1,...,\eta_{n_k}\r}  \r
\Vdash  \lusim{\alpha}_k=\nu.$$
\begin{flushright}
   $\blacksquare$
\end{flushright}

We restrict now our attention to ultrafilters $U$ which are P-points.
This will allow us to deal with arbitrary sets of ordinals in $V[C]$.
\\Recall the definition.

\begin{definition}\label{def-P-point}
$U$ is called a \textit{P-point} iff every non-constant (mod $U$) function $f:\kappa\to \kappa$
is almost one to one (mod $U$), i.e.
there is $A \in U$ such that for every $\delta<\kappa$,
$$|\{\nu\in A \mid f(\nu)=\delta\}|<\kappa.$$

\end{definition}

Note that, in particular, the projection to the normal ultrafilter $\pi$ is  almost one to one. Namely,

$$|\{\nu<\kappa \mid \pi(\nu)=\alpha\}|<\kappa, $$ for any $\alpha<\kappa$.
\\Denote by $U^{nor}$ the projection of $U$ to the normal ultrafilter.

\begin{lemma}\label{lem2}
Assume that $\mathbb{U}=\l U_a \mid 1\leq a \in [\kappa]^{<\omega}\r$ consists of P-point ultrafilters.
Suppose that $A \in V[C]\setminus V$ is an unbounded subset of $\kappa$. Then
$\kappa$ has cofinality $\omega$ in $V[A]$.

\end{lemma}
\pr
Work in $V$. Let $\lusim{A}$ be a name of $A$ and $\l s, {S} \r \in
P(\mathbb{U})$. Suppose for simplicity that $s$ is the empty sequence. Define by induction a subtree $T$ of ${S} $. For
each $\nu \in Lev_1(S)$ pick some subtree ${S'_\nu}$ of $S_{\l \nu \r}$  and
$a_\nu \subseteq \pi_{\l \r}(\nu)$ such that
$$\l \l \nu \r, {S}'_\nu \r \llvdash  \lusim{A} \cap \pi_{\l \r}(\nu)=a_\nu.$$
Let $S(0)'$ be a subtree of $S$ obtained by replacing $S_{\l \nu \r}$ by ${S'_\nu}$, for every $\nu \in Lev_1(S)$.

 Consider the function $\nu \rightarrow a_\nu$, where $\nu \in Suc_S(\l\r)$. By
normality of $\pi_{{\l \r}*}U_{\l\r}$ it is easy to find $A(0) \subseteq \kappa$ and
$T(0) \subseteq Lev_1(S(0)'), T(0) \in U_{\l \r}$ such that $A(0) \cap \pi_{\l\r}(\nu) =a_\nu$,
for every $\nu \in T(0)$. Set the first level of $T $ to be $T(0) $. Set $S(0)$ to be a subtree of $S(0)'$ obtained by shrinking the first level
to $T(0)$.
\\
Let now $\l \nu_1, \nu_2 \r\in Lev_2(S(0))$. So,  $\pi_{\l \nu_1 \r}(\nu_2) > \nu_1$. Find a subtree ${S}'_{\nu_1,\nu_2}$ of
$ {S(1)}_{\l \nu_1,\nu_2\r }, $
and $a_{\nu_0,\nu_1} \subseteq \pi_{\l \nu_1 \r}(\nu_2)$ such
that
$$\l \l\nu_1,\nu_2 \r, {S}'_{\nu_0,\nu_1} \r \llvdash \lusim{A} \cap \pi_{\l \nu_1 \r}(\nu_2)=a_{\nu_1,\nu_2}.$$
Let $S(1)'$ be a subtree of $S(0)$ obtained be replacing $S_{\l \nu_1,\nu_2 \r}$ by ${S'_{\nu_1,\nu_2}}$, for every $\l \nu_1, \nu_2 \r\in Lev_2(S(0))$.

Again, we 'consider the function $\nu \rightarrow a_\nu$, where  $\nu \in Suc_{S(1)}(\l\nu_1\r)$. By normality of $\pi_{{\l \nu_1\r}*}U_{\l \nu_1\r}$,  it is easy to find $A(\nu_1)
\subseteq \kappa$ and $T(\nu_1) \subseteq Suc_{S'(1)}(\l \nu_1\r), T(\nu_1) \in U_{\l \nu_1\r}$
such that $A(\nu_1) \cap \pi_{\l \nu_1\r}(\nu) =a_{\nu_1,\nu}$, for every $\nu \in
T(\nu_1)$. 

Define the set of the immediate successors of $\nu_1$
to be $T(\nu_1)$, i.e. $Suc_T(\nu_1)= T(\nu_1)$. Let $S(1)$ be a subtree of $S(1)'$ obtained this way.
\\
This defines the second level of $T$. Continue similar to define
further levels of $T$.
\\
We will have the following
property:

(*) for every $\l \eta_1,...,\eta_n \r \in T,$   $$\l \l
\eta_1,...,\eta_n \r, T_{\l \eta_1,...,\eta_n \r} \r \llvdash
\lusim{A}\cap \pi_{\l
\eta_1,...,\eta_{n -1}\r}(\eta_n) =A( \eta_1,...,\eta_{n-1} ) \cap \pi_{\l
\eta_1,...,\eta_{n-1} \r}(\eta_n).$$

A simple density argument implies that there is a condition which
satisfies (*) in the generic set.
Assume for simplicity that already $\l \l\r, T\r$ is such a condition.
Then, $C$ is a branch through $T^*$. Let $\l \kappa_n \mid n<\omega \r =C$.
So, for every $n<\omega$,
$$ A \cap \pi_{\l \kappa_0,...,\kappa_{n-1}\r}(\kappa_n)= A(\kappa_0,...,\kappa_{n-1}) \cap \pi_{\l \kappa_0,...,\kappa_{n-1}\r}(\kappa_n).$$
Let us work now in $V[A]$ and define by induction a sequence $\l\eta_n \mid n<\omega \r$
as follows.
Consider $A(0)$. It is a set in $V$, hence $A(0)\not = A$. So there is $\eta$ such that for every $\nu \in Lev_1(T)$ with $\pi_{\l\r}(\nu)\geq \eta$ we have $A \cap\pi_{\l\r}(\nu) \not = A(0)\cap \pi_{\l\r}(\nu)$. Set $\eta_0$ to be the least such $\eta$.
\\Turn to $\eta_1$. Let $\xi \in Lev_1(T)$ be such that $\pi_{\l\r}(\xi)< \eta_0$.
Consider $A(\xi)$. It is a set in $V$, hence $A(\xi)\not = A$. So there is $\eta$ such that for every
$\nu \in Lev_2(T_{\l \xi\r})$   with $\pi_{\l \xi \r}(\nu)\geq \eta$ we have $A \cap\pi_{\l \xi \r}(\nu) \not = A(\xi)\cap \pi_{\l \xi \r}(\nu)$. Set $\eta(\xi)$ to be the least such $\eta$.
Now define $\eta_1$ to be $\sup(\{\eta(\xi) \mid \pi_1(\xi)<\eta_0\})$.
The crucial point now is that the number of $\xi$'s with $\pi_{\l  \r}(\xi)<\eta_0$ is less than $\kappa$,  since $U_{\l\r}$ is a P-point.
\\If $\eta_1=\kappa$, then the cofinality of $\kappa$ (in $V[A]$) is at most $\eta_0$. So it must be $\omega$ since the Prikry forcing
used does not add new bounded subsets to $\kappa$, and we are done.\\
Let us argue however that this cannot happen and always
$\eta_1<\kappa$.

\begin{claim}\label{clm<}

$\eta_1<\kappa$.

\end{claim}
\pr
Suppose otherwise. Then
$$\sup(\{\eta(\xi) \mid \pi_{\l  \r}(\xi)<\eta_0\})=\kappa.$$
Hence for every $\alpha<\kappa$ there will be $\xi$ with $\pi_{\l  \r}(\xi)<\eta_0$
such that $$A \cap \alpha=A(\xi)\cap \alpha.$$
Then,  for every  $\alpha<\kappa$ there will be $\xi, \xi'$ with $\pi_{\l  \r}(\xi),\pi_{\l  \r}(\xi')<\eta_0$
such that $$A(\xi) \cap \alpha=A(\xi')\cap \alpha.$$
Now, in $V$, set $\rho_{\xi,\xi'}$ to be the least $\rho<\kappa$ such that $$A(\xi) \cap \rho\not=A(\xi')\cap \rho,$$
if it exists and $0$ otherwise, i.e. if $A(\xi)=A(\xi')$. Let
$$Z=\{\rho_{\xi,\xi'}\mid \pi_{\l  \r}(\xi),\pi_{\l  \r}(\xi')<\eta_0\}.$$
Then $|Z|^V<\kappa$, since the number of possible $\xi,\xi'$ is less than $\kappa$.
But $Z$ should be unbounded in $\kappa$ due to the fact that
 for every $\alpha<\kappa$ there will be $\xi$ with $\pi_{\l  \r}(\xi)<\eta_0$
such that $A \cap \alpha=A(\xi)\cap \alpha$ and $A \not= A(\xi)$.
Contradiction.
\begin{flushright}
   $\blacksquare_{of\ the\ claim}$

\end{flushright}

Suppose that $\eta_0,...,\eta_n<\kappa$ are defined. Define $\eta_{n+1}$.
Let $\l \xi_0,...,\xi_n \r$ be in $T$.
Consider $A(\xi_0,...,\xi_n)$. It is a set in $V$, hence $A(\xi_0,...,\xi_n)\not = A$. So there is $\eta$ such that for every
$\nu \in Lev_{n+2}(T_{\l \xi_0,...,\xi_n \r})  $ with $\pi_{\l \xi_0,...,\xi_n \r}(\nu)\geq \eta$ we have $A \cap\pi_{\l \xi_0,...,\xi_n \r}(\nu) \not = A(\xi_0,...\xi_n)\cap \pi_{\l \xi_0,...,\xi_n \r}(\nu)$. Set $\eta(\xi_0,...\xi_n)$ to be the least such $\eta$.
Now define $\eta_{n+1}$ to be $\sup(\{\eta(\xi_0,...\xi_n) \mid \pi_{\l\r}(\xi_0)<\eta_0,...,\pi_{\l \xi_0,...,\xi_{n-1} \r}(\xi_n)<\eta_n\})$.
\\Each relevant ultrafilter is a P-point, and so, the number of relevant  $\xi_0,...\xi_n$ is bounded in $\kappa$.
So,  $\eta_{n+1}<\kappa$, as in the claim above.

This completes the definition of the sequence $\l \eta_n \mid n<\omega \r$.
\\Let us argue that it is cofinal in $\kappa$.
\\Suppose otherwise.
\\Note that the sequence $\l \pi_{{\l \kappa_0,...,\kappa_{n-1}\r}}(\kappa_n) \mid n<\omega\r$ is unbounded in $\kappa$.
\\
Let $k$ be the least such that $\pi_{{\l \kappa_0,...,\kappa_{k-1}\r}}(\kappa_k) > \sup(\{\eta_n \mid n<\omega \})$.
Then $$A \cap \pi_{{\l \kappa_0,...,\kappa_{k-1}\r}}(\kappa_k) =A(\kappa_0,...,\kappa_{k-1})\cap \pi_{\l \kappa_0,...,\kappa_{k-1}\r}(\kappa_k).$$
This is impossible, since $\eta_k<\pi_{\l \kappa_0,...,\kappa_{k-1}\r}(\kappa_k)$.
\begin{flushright}
   $\blacksquare$
\end{flushright}


\begin{theorem}\label{P-point}
Let $\mathbb{U}=\l U_a \mid a \in [\kappa]^{<\omega}\r$ consists of P-point ultrafilters over $\kappa$.
Then for every new set of ordinals $A$ in $V^{P(\mathbb{U})}$, $\kappa$ has cofinality $\omega$ in $V[A]$.
\end{theorem}
\pr
Let $A$ be a new set of ordinals in $V[G]$, where $G\subseteq {P(\mathbb{U})}$ is generic.
By Lemma \ref{lem2}, it is enough to find a new subset of $A$ of size $\kappa$.
\\Suppose that every subset of $A$ of size $\kappa$ is in $V$.
Let us argue that then $A$ is in $V$ as well. Let $\lambda=\sup(A)$.
\\The argument is similar to
 \cite{PrikryCase}(Lemma 0.7).\\
Note that $(\calP_{\kappa^+}(\lambda))^V$ remains stationary in $V[G]$, since $P(\mathbb{U})$ satisfies $\kappa^+-$c.c.
For each $x\in (\calP_{\kappa^+}(\lambda))^V$ pick $\l s_x, S_x \r \in G$
such that
$$\l s_x, S_x \r \llvdash \lusim{A} \cap x = A \cap
x.$$ There are a stationary $E \subseteq (\calP_{\kappa^+}(\lambda))^V$ and $s \in
[\kappa]^{<\omega}$ such that for each $x \in E$ we have
$s=s_x$. Now, in $V$, we consider $$H=\{\l s, T \r \in P(U)
\mid \exists x\in \calP_{\kappa^+}(\lambda) \exists a \subseteq x \quad \l s,
T \r \llvdash \lusim{A}\cap x=a \}.$$ Note that if $\l s, T \r,
\l s, T' \r \in P(U)$ and for some $x\subseteq y$ in $\calP_{\kappa^+}(\lambda)$,
$a\subseteq x, b \subseteq y$ we have
$$\l s,
T \r \llvdash \lusim{A}\cap x=a \text{ and } \l s, T' \r
\llvdash \lusim{A}\cap y=b,$$ then $b\cap x=a$. Just
conditions of this form are compatible, and so they cannot force
contradictory information.
\\
Apply this observation to $H$. Let
$$X =\{a \subseteq \lambda \mid \exists \l s, S \r \in H \quad \exists x \in \calP_{\kappa^+}(\lambda)
\l s, T \r \llvdash \lusim{A}\cap x=a \}.$$ Then necessarily,
$\bigcup X =A$.
\begin{flushright}
   $\blacksquare$
\end{flushright}

We do not know wether $V[A]$ for $A \in V[C]\setminus V$ is equivalent to a single $\omega-$sequence even for $A\subseteq \kappa^+$.
The problematic case is once $U_n$'s have $\kappa^+-$many different ultrafilters below in the Rudin-Keisler order.

\begin{theorem}\label{cor2}
Assume that there is no inner model with $o(\alpha)=\alpha^{++}$.
Let $U$ be $\kappa-$complete ultrafilter  over $\kappa$ and  $V=L[\vec{E}]$, for a coherent sequence of measures $\vec{E}$.
Force with the Prikry forcing with $U$. Suppose that $A$ is a new set of ordinals in a generic extension.
 Then the cofinality of $\kappa$ is $\omega$ in $V[A]$.

\end{theorem}
\pr
Consider $$i_U:V\to M\simeq V^\kappa/U.$$
By Mitchell \cite{Mit}, $i_{U}$ is an iterated ultrapower using measures from $\vec{E}$ and images of $\vec{E}$. In addition we have that  ${}^\kappa M\subseteq M$.
Hence it should be a finite iteration using measures from $\vec{E}$.
Since $\kappa$ is the critical point, no measures below $\kappa$ are involved and the first one applied is a measure on $\kappa$ in $\vec{E}$.
Denote it by $E_0$ and let
$$i_0:V\to M_1$$ be the corresponding embedding.
Let $\kappa_1=i_0(\kappa)$. Rearranging, if necessary, we can assume that the next step was to use a measure $E_1$ over $\kappa_1$ from $i_0(\vec{E})$.
So, it is either the image of one of the measures of $\vec{E}$ or $E_0-Lim \l E^\xi \mid \xi<\kappa \r$, where $\l E^\xi \mid \xi<\kappa \r$
is a sequence of measures over $\kappa$ from $\vec{E}$ which represents in $M_1$ the measure used over $\kappa_1$.
\\Let $$i_1:M_1\to M_2$$ be the corresponding embedding and $\kappa_2=i_1(\kappa_1)$.
\\$\kappa_2$ can be moved further in our iteration, but  only finitely many times.
Suppose for simplicity that it does not move.
\\If nothing else is moved then $U$ is equivalent to  $E_0-Lim \l E^\xi \mid \xi<\kappa \r$
and theorem 6.7 easily provides the desired conclusion.
\\Suppose $i_1\circ i_0$ is not $i_U$. Then some measures from $i_1\circ i_0(\vec{E})$ with critical points in the intervals
$(\kappa,\kappa_1),(\kappa_1,\kappa_2)$ are applied.
Again, only finitely many can be used.
\\Thus suppose for simplicity that only one is used in each interval. The treatment of a general case is more complicated only
due to notation.
\\So suppose that a measure $E_2$ with a critical point $\delta\in (\kappa,\kappa_1)$ is used on the third step of the iteration.
\\Let $$i_2:M_2\to M_3$$ be the corresponding embedding.
Note that the ultrafilter ${\cal V}$ defined by
$$X \in {\cal V} \text{ iff } i_2(\delta)\in i_2\circ i_1\circ i_0(X)$$
is $P-$point. Thus, a function
 $f:\kappa\to \kappa$ which represents $\delta$ in $M_1$, i.e. $\delta=i_0(f)(\kappa)$, will witness this.
\\Similar an ultrafilter used in the interval  $(\kappa_1,\kappa_2)$ will be $P-$point in $M_1$, and so, in $V$,
it will be equivalent to a limit of $P-$points.
\\So such situation is covered by \ref{P-point}.
\begin{flushright}
   $\blacksquare$
\end{flushright}

\newpage
\section{Prikry forcing may add a Cohen subset.}

Our aim here will be to show the following:

\begin{theorem}\label{Cohen subset} Suppose that $V$ satisfies GCH and $\kappa$ is a measurable cardinal. Then in a generic cofinality preserving extension there is a $\kappa-$complete ultrafilter $U$ over $\kappa$ such that the Prikry forcing with $U$ adds a Cohen subset to $\kappa$ over $V$.
In particular, this forcing has a non-trivial subforcing which preserves regularity of $\kappa$.
\end{theorem}

By \cite{PrikryCase} such $F$ cannot by normal and by 6.6 $F$ cannot be a P-point ultrafilter, since in any Cohen extension, \(\kappa\) stays regular.

Note that the above situation is impossible in $L[\mu]$. Just every $\kappa-$complete ultrafilter over the measurable $\kappa$ is Rudin-Kiesler equivalent to $\mu^n$, for some $n,1\leq  n<\omega$(
see \cite[Lemma 19.21]{Jech}). But the Prikry forcing with $\mu^n$ is the same as the Prikry forcing with $\mu$ which is a normal measure.

 We start with a GCH model with a measurable.
Let \(\kappa\) be a measurable and \(U\) a normal measure on \(\kappa\).\\
Denote by
  \(j_U:V\rightarrow N\simeq Ult(V,U)\) the corresponding elementary embedding.

Define an iteration  \(\big\langle P_\alpha,Q_\beta\mid \alpha\leq\kappa,\beta<\kappa\big\rangle\) with Easton support as follows.
Set
 \(P_0=0\).
Assume that \(P_\alpha\) is defined.  Set \(\lusim{Q}_\alpha\) to be  the trivial forcing unless \(\alpha\) is an inaccessible cardinal.
\\If $\alpha$ is an inaccessible cardinal, then let  ${Q}_\alpha=Q_{\alpha 0}*\lusim{Q}_{\alpha 1}$, where
$Q_{\alpha 0}$ is an atomic forcing consisting of three elements $0_{Q_{\alpha 0}},x_\alpha,y_\alpha$, such that $x_\alpha,y_\alpha$ are two incompatible elements which are stronger than $0_{Q_{\alpha 0}}$.
\\
Let $\lusim{Q}_{\alpha 1}$ be trivial once $y_\alpha$ is picked and let it be the Cohen forcing at $\alpha$, i.e. \begin{center}

$Cohen(\alpha,2)=\{f:\alpha\rightarrow2\mid |f|<\alpha\}$
\end{center}
once $x_\alpha$ was chosen.

Let \(G_{\kappa}\subseteq P_{\kappa}\) be a generic.
We extend now the embedding $$j_U:V \to N,$$ in $V[G_\kappa]$,
to $$j^*_U:V[G_\kappa]\rightarrow N[G_\kappa*G_{[\kappa,j_U(\kappa))}],$$
for some $G_{[\kappa,j_U(\kappa))}\subseteq P_{[\kappa,j_U(\kappa))}$ which is $N[G_\kappa]-$generic for $P_{j_U(\kappa)}/G_\kappa.$
This can be done easily, once over $\kappa$ itself in $Q_{\kappa 0}$, we pick $y_\kappa$, which makes the forcing $Q_\kappa$ a trivial one.
\\This shows, in particular, that $\kappa$ is still a measurable in $V[G_\kappa]$, as witnessed by an extension of $U$.

Consider now the second ultrapower $N_2\simeq \Ult(N, j_U(U))$.
\\Denote $j_U$ by $j_1$, $N$ by $N_1$. Let $$j_{12}:N_1 \to N_2$$ denotes the ultrapower embedding of $N_1$ by $j_1(U)$.
Let $j_2=j_{12}\circ j_1$. Then $$j_2:V \to N_2.$$
Let us extend, in $V[G_\kappa]$, the embedding  $$j_{12}:N_1 \to N_2$$ to
$$j_{12}^*:N_1[G_\kappa*G_{[\kappa,j_1(\kappa))}] \to N_2[G_\kappa*G_{[\kappa,j_1(\kappa))}*G_{[j_1(\kappa),j_2(\kappa))}]$$
in a standard fashion, only this time  we pick $x_{j_1(\kappa)}$ at stage $j_1(\kappa)$ of the iteration.
Then a Cohen function should be constructed over $j_1(\kappa)$, which is not at all problematic to find in $V[G_\kappa]$.

Now we will have
$$j_2\subseteq j_2^*:V[G_\kappa]\to N_2[G_\kappa*G_{[\kappa,j_1(\kappa))}*G_{[j_1(\kappa),j_2(\kappa))}]$$
which is the composition of $j_1^*$ with $j_{12}^*$.

 Define  a $\kappa-$complete ultrafilter \(W\) over $\kappa$ as follows:
\begin{center}
\(X\in W\) iff \(X\subseteq\kappa\) and \(j_1(\kappa)\in j^*_2(X).\)
\end{center}

\begin{proposition}
 \textit{$W$ has the following basic properties:}
\begin{enumerate}
    \item \(W\cap V=U\).
    \item $\{\alpha<\kappa \mid x_\alpha  \text{ was picked at the stage } \alpha \text{ of the iteration }\}\in W$.
    \item \textit{if \(C\subseteq\kappa\) is a club, then \(C\in W\). Moreover} $$\{\nu\in C\mid \nu \text{ is an inaccessible}\}\in W.$$
\end{enumerate}

\end{proposition}
\pr
(1) and (2) are standard. Let us show only (3).
Let $C\subseteq \kappa$ be a club. Then, in $N_2$,
  \(j_2(C)\) is a club at \(j_2(\kappa)\). In addition, \(j_2(C)\cap \kappa_1=j_1(C)\). Now, \(j_1(C)\) is a club in \(j_1(\kappa)\). It follows that \(j_1(\kappa)\in j_2(C)\).
 \\In order to show that
 $$\{\nu\in C\mid \nu \text{ is an inaccessible}\}\in W,$$
 just note that $j_1(\kappa)$ is an inaccessible in $N_2$, and so $W$ concentrates on inaccessibles.
\begin{flushright}
 \(\blacksquare\)
 \end{flushright}
Force with \(Prikry(W)\) over \(V[G_\kappa]\).\\ Let $$C=\l \eta_n \mid n<\omega \r$$ be a generic Prikry sequence. \\By (2) in the previous proposition, there is $n^*<\omega$ such that for every $m\geq n^*$, at the stage $\eta_m$ of the forcing $P_\kappa$, $x_{\eta_m}$ was picked,
and, hence, a Cohen function $f_{\eta_m}:\eta_m \to 2$ was added.

Define now $H:\kappa\to 2$ in $V[G_\kappa, C]$ as follows:

$$H= f_{\eta_{n^*}}\cup \bigcup_{n^*\leq m<\omega} f_{\eta_{m+1}}\upr [\eta_m,\eta_{m+1}).$$
\\
\begin{proposition}

 \textit{\(H\) is a Cohen generic function for $\kappa$ over \(V[G_\kappa]\).}\end{proposition}
\pr Work in $V[G_\kappa]$. Let \(D\in V[G_\kappa]\) be a dense open subset of \(Cohen(\kappa)\). Consider a set
$$C=\{\alpha<\kappa\mid \text{ if } \alpha \text{ is an inaccessible, then } D \cap V_\alpha[G_\alpha] \text{ is a dense open subset of } Cohen(\alpha) \text{ in } V[G_\alpha]  \}.$$
\\
\\
\textbf{claim 1} $C$ is a club.
\\
\pr
Suppose otherwise. Then $S=\kappa\setminus C$ is stationary. It consists of inaccessible cardinals by the definition of $C$.\\
 Pick a cardinal $\chi$ large enough and consider an elementary submodel $X$ of
$\l H_\chi, \in\r$ such that
\begin{enumerate}
  \item $X \cap (V_\kappa)^{V[G_\kappa]}=(V_\delta)^{V[G_\kappa]},$ for some $\delta\in S$.
  \item $\kappa, P_\kappa,D\in X$.
  \end{enumerate}

 Note that it is possible to find such $X$ due to stationarity of $S$.
Note also that  $(V_\kappa)^{V[G_\kappa]}=V_\kappa[G_\kappa]$ and $(V_\delta)^{V[G_\kappa]}=V_\delta[G_\delta]$,
since the iteration $P_\kappa$ splits nicely at inaccessibles.

Let us argue that $D \cap V_\delta[G_\delta]$ is  a dense open subset of $ Cohen(\delta) \text{ in } V[G_\delta] $.
\\Just note that $$D \cap X =D \cap X \cap (V_\kappa)^{V[G_\kappa]}=D \cap(V_\delta)^{V[G_\kappa]}=D\cap V_\delta[G_\delta].$$
So let $q \in (Cohen(\delta))^{V_\delta[G_\delta]}$. Then $q \in X$. Remember $X \preceq H_\chi$. So,
$$X\models D \text{ is dense open },$$
hence there is $p\geq q, p\in D\cap X$. But then, $p \in D \cap V_\delta[G_\delta]$, and we are done.
\\Contradiction.
\\
\begin{flushright}
   $\blacksquare$ of the claim.
\end{flushright}

It follows now that $C \in W$. Hence there is $n^{**}\geq n^*$ such that for every $m,n^{**}\leq m<\omega$,
$$\eta_m \in C.$$
So, for every $m,n^{**}\leq m<\omega$,
$$f_{\eta_m} \in D,$$
since $D$ is open.
\\It is almost what we need, however $H\upr \eta_m$ need not be  $f_{\eta_m}$, since an initial segment may was changed.
\\In order to overcome this, let us note the following basic property of the Cohen forcing:\\
\\
\textbf{Claim 2}
Let $E$ be a dense open subset of $Cohen(\kappa,2)$, then there is a dense subset $E^*$ of $E$ such that
$\text{ for every } p \in E^* \text{ and every inaccessible cardinal } \tau\in \dom(p) \text{ for every } q:\delta \to 2,
p\upr [\delta,\kappa)\cup q \in E^*.$
\\
\\
The proof is an easy use of $\kappa-$completeness of the forcing.

Now we can finish just replacing $D$ by its dense subset which satisfies the conclusion of the claim.
Then, $H\upr \eta_m$ will belong to it as a bounded change of $f_{\eta_m}$.
\\So we are done.\\
\begin{flushright}
  \(\blacksquare\)

\end{flushright}

\newpage
\section{Further directions.}

One of possible further directions is to extend our results from the Magidor forcing to the Radin forcing.
Note that we cannot claim that every subforcing of the Radin forcing is equivalent to Radin forcing.
Thus, the negation of \(o^{\vec{U}}(\kappa_i)<\min(\nu\mid 0<o^{\vec{U}}(\nu))\) provides a counterexample.
However, it is still reasonable that  every  set in a Radin extension is equivalent to a subsequence of the Radin sequence.
We conjecture   that this is the case.

An other direction is to proceed further with the Prikry forcing with $P-$points ultrafilters and to prove that every subforcing of it is equivalent
to a Prikry forcing. The complications starts once a $P-$point has more than $\kappa$ many generators. In such situations it is easy to construct a subset of $\kappa^+$ which is not equivalent to any of its initial segments. The opposite was crucial for the arguments  of \cite{PrikryCase} with
a normal ultrafilter.
We conjecture that it is possible to overcome this problem and that  every subforcing of the Prikry forcing with $P-$points ultrafilters  is indeed equivalent
to Prikry forcing.

Let us conclude with few questions.

\textbf{Question 1.} \emph{Is every  set in a Radin extension  equivalent to a subsequence of the Radin sequence?}

\textbf{Question 2.}\emph{ Is every subforcing of the Prikry forcing with $P-$points ultrafilters   equivalent
to a Prikry forcing?}

\textbf{Question 3.} \emph{Characterize all $\kappa-$complete ultrafilters $U$ over $\kappa$ such that $\kappa$ changes its cofinality
in $V[A]$, for any new set $A$ in the Prikry extension with $U$.}

In section 6 a rather large class of such ultrafilters was presented. It includes $P-$points, their products and limits.
But are there other ultrafilters like this?

\end{document}